\documentclass[11pt]{amsart}
\usepackage[latin1]{inputenc}
\usepackage{amsfonts}
\usepackage[toc,page,title,titletoc,header]{appendix}
\usepackage{graphicx,psfrag,epsfig,multirow,caption,subcaption}
\usepackage{amssymb,amsmath,amscd,amsthm,amssymb,verbatim,setspace}
\numberwithin{equation}{section}
\usepackage{mathrsfs,wrapfig}
\usepackage{indentfirst}
\usepackage{extarrows}
\usepackage{float}
\usepackage[dvipsnames]{xcolor}
\definecolor{revone}{RGB}{0,55,145}

\usepackage[colorlinks=true, backref=page]{hyperref}

\newtheorem{theorem}{Theorem}[section]
\newtheorem{definition}[theorem]{Definition}
\newtheorem{conjecture}[theorem]{Conjecture}
\newtheorem{proposition}[theorem]{Proposition}
\newtheorem{lemma}[theorem]{Lemma}
\newtheorem{remark}[theorem]{Remark}
\newtheorem{corollary}[theorem]{Corollary}

\newtheorem*{remark*}{Remark}

\numberwithin{equation}{section}

\newcommand{\eps}{\varepsilon}

\newcommand{\dt}{\delta}
\newcommand{\al}{\alpha}

\newcommand{\bn}{\mathbf{n}}
\newcommand{\br}{\mathbf{r}}

\newcommand{\bM}{\mathbf{M}}
\newcommand{\bS}{\mathbb{S}}

\newcommand{\cA}{\mathcal{A}}

\newcommand{\cH}{\mathcal{H}}
\newcommand{\cE}{\mathcal{E}}
\newcommand{\cB}{\mathcal{B}}

\newcommand{\cC}{\mathcal{C}}
\newcommand{\cL}{\mathcal{L}}
\newcommand{\cK}{\mathcal{K}}

\newcommand{\cQ}{\mathcal{Q}}

\newcommand{\cG}{\mathcal{G}}

\newcommand{\cF}{\mathcal{F}}
\newcommand{\cI}{\mathcal{I}}

\newcommand{\bH}{\mathbf{H}}

\newcommand{\sL}{\mathscr{L}}

\newcommand{\R}{\mathbb{R}}

\newcommand{\A}{\mathbb{A}}
\newcommand{\pr}{\partial}

\DeclareMathOperator{\Ric}{Ric}
\DeclareMathOperator{\Hess}{Hess}

\DeclareMathOperator{\spt}{spt}
\DeclareMathOperator{\Graph}{Graph}

\DeclareMathOperator{\tr}{tr}

\DeclareMathOperator{\id}{\text{Id}}

\title[Cylindrical singularities I]{Mean curvature flows with cylindrical singularities I: the normal forms and nondegeneracy}
\author{Ao Sun}
\address{Lehigh University, Department of Mathematics, Chandler-Ullmann Hall, Bethlehem, PA 18015}
\email{aos223@lehigh.edu}
\author{Jinxin Xue}
\address{New Cornerstone Science Laboratory, Department of Mathematics, Rm A115, Tsinghua University, Haidian District, Beijing, 100084}
\email{jxue@tsinghua.edu.cn}
\date{\today}

\begin{document}

	\begin{abstract}
		This paper studies the dynamics of mean curvature flow as it approaches a cylindrical singularity. We prove that the rescaled mean curvature flow converging to a smooth generalized cylinder can be written as a graph over the cylinder in a ball of radius $K\sqrt{\tau}$, and a normal form of the asymptotics, where $\tau$ is the rescaled time parameter. Such normal forms provide more precise dynamical and geometric information for mean curvature flows with cylindrical singularities beyond {\L}ojasiewicz inequalities.  Using the normal form, we can define the nondegeneracy of cylindrical singularities, and we show that nondegenerate cylindrical singularities are isolated in space, have a mean convex neighborhood, and are type-I.  
	\end{abstract}

	\maketitle

	\section{Introduction}
	This paper focuses on mean curvature flows with cylindrical singularities. Mean curvature flow is a fundamental geometric flow that has attracted considerable attention in diverse fields, including geometry, partial differential equations, and applied mathematics.
	
	Mean curvature flow is defined as a family of hypersurfaces $\{\bM_t\}_{t\in I}$ evolving in $\mathbb{R}^{n+1}$ according to the equation $\pr_{t} x=\vec{H}(x)$, where $\vec H$ is the mean curvature vector. Singularities must appear in mean curvature flows that originate from smooth embedded closed hypersurfaces. Therefore, understanding singularities is crucial to understanding mean curvature flows. The analysis of singularities involves a blow-up procedure. Huisken \cite{Hu1} introduced the concept of rescaled mean curvature flow to study these singularities. The rescaled mean curvature flow is a family of hypersurfaces $\{M_\tau\}_{\tau\in[0,\infty)}$ satisfying the equation
	$
	\pr_\tau x=\vec{H}+\frac{x^\perp}{2}.
	$
	When the first singularity of the mean curvature flow $\{\bM_t\}_{t\in[-1,0)}$ occurs at the spacetime point $(0,0)\in \R^{n+1}\times \R$, the corresponding rescaled mean curvature flow $\{M_\tau\}_{\tau\in[0,\infty)}$ that captures this singularity is defined as follows:
	\begin{equation}\label{EqMCFRMCF}
		M_\tau=e^{\tau/2}\bM_{-e^{-\tau}}, \quad \tau\in [0,\infty).
	\end{equation}
	The singularity models, known as shrinkers, are obtained as a (subsequential) limit of the rescaled mean curvature flow, and they satisfy the equation $\vec H+\frac{x^\perp}{2}=0$.
	
	The rescaled mean curvature flow can be viewed as a dynamical system, a perspective introduced by Colding-Minicozzi in a series of papers \cite{CM1,CM2,ColdingMinicozzi19_Dynamics,CM7}. 
	To describe this system, let $\Sigma\subset\mathbb{R}^{n+1}$ be a hypersurface, and define its Gaussian area as follows:
	$$
	\cF(\Sigma):=(4\pi)^{-n/2}\int_{\Sigma}e^{-\frac{|x|^2}{4}}d\cH^n(x),
	$$ where $\cH^n$ represents the $n$-dimensional Hausdorff measure. In this context, a rescaled mean curvature flow $M_\tau$ is the negative gradient flow of the Gaussian area, and a shrinker is a fixed point of this flow. 
	
	In a pioneering work \cite{CM1}, Colding-Minicozzi proved that the only stable shrinkers from a variational point of view are the hyperplane, the sphere, and the generalized cylinders $\cC_{n,k}:=\mathbb{S}^{n-k}(\sqrt{2(n-k)})\times\R^k$. Based on their result, there has been much recent research on generic mean curvature flows that avoid other shrinkers as the singularity models; cf. \cite{ColdingMinicozzi19_Dynamics, SX1, SX2, CCMS, CCS, BamlerKleiner23_multiplicity}. On the other hand, although the generalized cylinders are stable from a variational point of view, the singular sets modeled by them can still be very complicated. One famous example is the ``marriage ring'', namely, a very thin rotationally symmetric torus in $\R^3$, and it shrinks to the circle as the singular set, showing that the singular set can have high dimensions. Another example is the ``peanut surface'' constructed by Altschuler-Angenent-Giga \cite{AAG}, showing that a compact mean curvature flow can shrink to a point while the tangent flow at the point is the noncompact cylinder.
	
	In the problem list \cite{I1}, Ilmanen proposed the following conjecture. This conjecture was also noted by Colding-Minicozzi-Pedersen in \cite[Conjecture 7.1]{CMP}.
	
	\begin{conjecture}[Generic isolatedness conjecture]\label{ConIsolate}
		Mean curvature flow with generic initial data only has isolated spacetime singularities. 
	\end{conjecture}

	In this paper, we study the cylindrical singularities to establish the necessary tools to study the dynamics of cylindrical singularities, and we present a perspective from dynamics to hint at why this conjecture should be true; in another forthcoming paper (which is the second part of the first version of this paper), we will discuss further dynamical properties of the generic cylindrical singularities. Our focus is on the mean curvature flow of closed embedded hypersurfaces in $\R^{n+1}$. The basic assumption is that 
	\begin{equation}\label{eq:StarCondition}\tag{$\star$}
		\parbox{5.5in}{%
			\emph{$\{\bM_t\}_{t\in[-1,0)}$ is a mean curvature flow, with bounded entropy, with a cylindrical singularity modeled on $\cC_{n,k}$ at the spacetime point $(0,0)$, and $\{M_\tau\}_{\tau\in [0,\infty)}$ is its associated rescaled mean curvature flow that converges to $\cC_{n,k}$ in the $C^\infty_{loc}$-sense.}%
		}
	\end{equation}
	
	The reader may interpret ``mean curvature flow'' as a smooth embedded mean curvature flow. Nevertheless, the discussion in this paper can be applied to weak flows which satisfy the pseudolocality theorem (see Section \ref{SSPL}), such as the unit regular cyclic mod $2$ Brakke flows, see \cite{White09_CurrentsVarifolds}.

	\subsection{The normal forms}
	To describe our main results, we start with some notations. We fix $1\leq k< n$. Let $\varrho:=\sqrt{2(n-k)}$, and we use coordinates denoted as $x=(\theta, y)\in \cC_{n,k}:=\mathbb{S}^{n-k}(\varrho)\times\R^k \subset \mathbb{R}^{n-k+1}\times\mathbb{R}^k$. Here, $\{y_i\}_{i=1}^k$ denotes the coordinates on $\mathbb{R}^k$, and $\{\theta_j\}_{j=1}^{n-k+1}$ denotes the restriction of the coordinate functions of $\mathbb{R}^{n-k+1}$ to $\mathbb{S}^{n-k}(\varrho)$. We use $B_R$ to denote the Euclidean open ball of radius $R$ centered at the origin.
	
	On $\cC_{n,k}$, it is natural to introduce the Gaussian weighted $L^2$-space, because the linearized operator $L_{\cC_{n,k}}$ of the generalized cylinder is self-adjoint with respect to the Gaussian weighted inner product. Throughout this paper, for any domain $\Omega\subset\cC_{n,k}$ and $f:\Omega\to\R$, we define
	\[
	\|f\|^2_{L^2(\Omega)}=\int_{\Omega}|f(x)|^2 e^{-\frac{|x|^2}{4}}d\cH^n(x),
	\quad 
	\|f\|^2_{H^1(\Omega)}=\int_{\Omega}(|f(x)|^2 +|\nabla f(x)|^2)e^{-\frac{|x|^2}{4}}d\cH^n(x).
	\]
	
	Our first main theorem is the asymptotic profile of the rescaled mean curvature flow converging to a generalized cylinder. Such an expression is called a {\bf normal form}, which is a terminology in dynamical systems after Birkhoff. Throughout this paper, we use the big $O$ notation as follows: given a positive function $f:\R\to\R$, $O(f(\tau))$ means a term satisfying $\limsup\limits_{\tau\to\infty}|O(f(\tau))|\cdot (f(\tau))^{-1}<+\infty$.
	
	\begin{theorem}[$H^1$-normal form theorem]\label{ThmNF-sqrtt} Assume \eqref{eq:StarCondition}. Then we have the following dichotomy:
		\begin{enumerate}
			\item either for any $K>0$ and $\vartheta\in(0,1)$, there exists $T>0$  such that for $\tau>T$, the rescaled mean curvature flow is a graph over $\cC_{n,k}\cap B_{K\sqrt{\tau}}$, and up to a rotation in $\R^k$, the graphical function $u(\cdot,\tau):\ \cC_{n,k}\cap B_{K\sqrt{\tau}}\to \R$ has the following asymptotic in (weighted) $H^1$-norm
			\begin{equation}
				\left\|u(\theta,y,\tau)-\sum_{i\in \mathcal I}\frac{\varrho }{4\tau}(y_i^2-2)\right\|_{H^1(\cC_{n,k}\cap B_{K\sqrt{\tau}})}=O(\tau^{-1-\vartheta})
			\end{equation}
			as $\tau\to\infty$, where $\emptyset\neq \mathcal I\subset \{1,2,\ldots,k\}$ and $O$ depend on $\vartheta$ and $K$.  
			\item or we have $\|u(\cdot,\tau)\|_{H^1(\cC_{n,k}\cap B_{K\sqrt{\tau}})}=O(e^{-c\tau})$ for some $K>0,c>0$, and $\tau$ sufficiently large, when $\mathcal I=\emptyset$.
		\end{enumerate}
	\end{theorem}
	
	We also prove the following $C^1$-normal form theorem.
	
	\begin{theorem}[$C^1$-normal form theorem]\label{ThmNF-sqrtt-C1}
		Assume \eqref{eq:StarCondition}. Then we have the following dichotomy:
		\begin{enumerate}
			\item either for any $K>0$ and $\vartheta\in(0,1)$, there exists $T>0$ such that for $\tau>T$, the rescaled mean curvature flow is a graph over $\cC_{n,k}\cap B_{K\sqrt{\tau}}$, and the graphical function $u(\cdot,\tau):\ \cC_{n,k}\cap B_{K\sqrt{\tau}}\to \R$
			has the following asymptotic, up to a rotation in $\R^k$, in $C^1$-norm in $B_{K\sqrt{\tau}}$
			\begin{equation}
				\left\|u(\theta,y,\tau)-\varrho\left(\sqrt{1+\frac{ \sum_{i\in \mathcal I}(y_i^2-2)}{2\tau}}-1\right)\right\|_{C^1(\cC_{n,k}\cap B_{K\sqrt{\tau}})}=O(\tau^{-\vartheta}),
			\end{equation}
			as $\tau\to\infty$, where $\emptyset\neq\mathcal I\subset \{1,2,\ldots,k\}$ and $O$ depends on $\vartheta$ and $K$. 
			\item  or    we have $\|u(\cdot,\tau)\|_{C^1(\cC_{n,k}\cap B_{K\sqrt{\tau}})}=O(e^{-c\tau})$ for some $K>0,c>0$, and $\tau$ sufficiently large, when $\mathcal I=\emptyset$.
		\end{enumerate}
	\end{theorem}
	Combining the above two normal forms with the classical parabolic Schauder theory, we have the following estimate: 
	
	\begin{corollary}[$C^2$-normal form on a bounded domain]\label{ThmNF-sqrtt-C1bdd} In the first alternative of the last two theorems, for any $R>2n$, $\vartheta\in(0,1)$, there exists $T>0$ such that for $\tau>T$, the rescaled mean curvature flow is a graph over $\cC_{n,k}\cap B_{R}$, and the graphical function $u(\cdot,\tau):\ \cC_{n,k}\cap B_{R}\to \R$
		has the following asymptotic, up to a rotation in $\R^k$, in $C^2$ norm
		\begin{equation}
			u(\theta,y,\tau)=\sum_{i\in \mathcal I}\frac{\varrho }{4\tau}(y_i^2-2)+O(\tau^{-1-\vartheta}),
		\end{equation}
		as $\tau\to\infty$, where $\emptyset\neq\mathcal I\subset \{1,2,\ldots,k\}$ and $O$ is independent of $R$ but may depend on $\vartheta$. 
	\end{corollary}
	
	In a smaller region, we have a better error estimate. This estimate may have independent interest and we record it here.
	
	\begin{corollary}[Normal form theorem in $B_{\tau^{\kappa}}$]\label{cor_main_small_radius}
		Suppose $\kappa\in(0,1/2)$ and $\eta\in(0,1-2\kappa)$. Then there exists $T>0$ such that for $\tau>T$, the rescaled mean curvature flow is a graph over $\cC_{n,k}\cap B_{\kappa^\tau}$, and the graphical function $u(\cdot,\tau):\ \cC_{n,k}\cap B_{\kappa^\tau}\to \R$
		has the following asymptotic, up to a rotation in $\R^k$, in $C^1$-norm in $B_{K\sqrt{\tau}}$
		\begin{equation}
			\left\|u(\theta,y,\tau)-\sum_{i\in \mathcal I}\frac{\varrho }{4\tau}(y_i^2-2)\right\|_{C^1(\cC_{n,k}\cap B_{\tau^\kappa})}=O(\tau^{2\kappa-2+\eta})
		\end{equation}
		as $\tau\to\infty$, where $\emptyset\neq \mathcal I\subset \{1,2,\ldots,k\}$ and $O$ depend on $\kappa$ and $\eta$.  
	\end{corollary}

	The study of the normal form for MCFs with cylindrical singularities was initiated by Gang Zhou in a series of works \cite{G1, G2} etc, but restricted only to the case of $n=4$ and $k=3$ for technical reasons. Besides the significant simplification of the proof, the normal forms in full generality, as well as the exponential estimate in the case of complete degeneracy $\mathcal I=\emptyset,$ we would like to highlight the geometry conveyed by the normal forms
	beyond technicality.   Indeed, the normal forms should be considered as providing more precise information than the {\L}ojasiewicz inequality, since the latter is a consequence of analyticity appearing ubiquitously in a wide range of fields such as Ricci flow, minimal surfaces, MCFs, etc, while the former provides specific dynamical information for the particular singularity under study.  
	
	This paper is the first one in a series of papers of ours, as well as jointly with Zhihan Wang, to explore the geometric information contained in the above normal forms, as well as sharper ones later. In this paper, from the normal forms stated above, we establish the following geometric information: nondegenerate singularities $(\mathcal I=\{1,2,\ldots,k\})$ are isolated, which provides a mechanism for the isolatedness of cylindrical singularities as one step towards the Conjecture \ref{ConIsolate}, nondegenerate singularities have also mean convex neighborhoods and are of type I. Moreover, we also provide an asymptotic profile for MCFs developing a nondegenerate singularity. In later papers, we shall show that the normal form controls the topology change for the MCFs passing through nondegenerate singularities and the regularity of singularity sets in the degenerate case. 
	
	The technical innovation of our work consists mainly of two ingredients: the analysis of the dynamical system for the rescaled mean curvature flow near the cylinder and the regularizing effect of the Ornstein-Uhlenbeck operator discovered by Vel\'azquez \cite{V1, V2}. We would like to remark that the $H^1$-normal form and $C^1$-normal form capture the behavior of the rescaled mean curvature flow in different regions: the $H^1$-normal form mainly captures the behavior inside the region of radius $O(\tau^{\kappa})$ with $\kappa$ close to $0$, while the $C^1$-normal form mainly captures the behavior near the boundary of the region of radius $O(\tau^{1/2})$.

	It is worth mentioning that in the investigation of ancient flows, Angenent-Daskalopoulos-Sesum \cite{ADS1} made the novel discovery that the asymptotic form can provide insights into the geometry of ancient flows. Furthermore, in recent research on the classification of ancient solutions of mean curvature flow \cite{ADS1, ADS2, BC1, BC2, CHH, CHHW, DuZhu22_spectral}, the asymptotic expansion of the rescaled mean curvature flow over the cylinders at $-\infty$ time plays a crucial role. However, it is important to note that the backward expansion differs from the forward expansion presented in this work, yielding distinct geometric information. For example, the linear asymptotic modes in the study of ancient flows form a finite-dimensional space, while in the study of forward flow, the linear asymptotic modes form an infinite-dimensional space.
	
	\subsection{Nondegeneracy and isolatedness}
	From the normal form that we derived, a natural definition arises.
	
	\begin{definition}\label{DefDeg}
		A cylindrical singularity as in Theorem \ref{ThmNF-sqrtt} $($equivalently as in Theorem \ref{ThmNF-sqrtt-C1}$)$ is \begin{itemize}
			\item \textbf{nondegenerate} if $\mathcal I=\{1,2,\ldots,k\}$;
			\item {\bf partially nondegenerate} if $\emptyset\neq\mathcal I\subsetneqq\{1,2,\ldots,k\}$;
			\item \textbf{degenerate}  if $\cI=\emptyset$.   
		\end{itemize}
	\end{definition}
	
	The concept of nondegeneracy for neckpinching singularities was initially introduced and studied by Angenent-Vel\'azquez \cite{AV} in the context of rotationally symmetric mean curvature flows. In the rotationally symmetric setting, the cylindrical singularities are modeled by $\cC_{n,1}$. Angenent-Vel\'azquez also constructed examples of mean curvature flows with degenerate singularities. It is worth noting that Schulze-Sesum \cite{SS} also defined a notion of nondegenerate neckpinch for a singularity modeled by $\cC_{n,1}$. Their nondegenerate singularities are characterized by the blow-up of limits and the disconnecting of the flow. It seems that our notion of nondegeneracy for $\cC_{n,1}$ implies theirs, but not vice versa: in fact, some of those degenerate singularities constructed by Angenent-Vel\'azquez in \cite{AV} seem to satisfy Schulze-Sesum's nondegenerate notion in \cite{SS}. 
	
	We discovered several geometric and dynamic properties of nondegenerate singularities. The first significant result from the normal form theorem is the following theorem concerning the isolatedness of nondegenerate cylindrical singularities.
	
	\begin{theorem}[Isolatedness theorem]\label{ThmIsolate}
		Assume \eqref{eq:StarCondition} and the singularity is nondegenerate, then the singularity is isolated. In other words, there exists $\dt_1>0$ such that $(0,0)$ is the only singularity of $\bM_t$ in the spacetime neighborhood $B_{\dt_1}\times(-\delta_1^{2},0]$.
	\end{theorem}

	Another relevant property is the mean convex neighborhood. It was conjectured that any cylindrical singularity has a mean convex neighborhood. This conjecture was proved by Choi-Haslhofer-Heshkovitz \cite{CHH} for dimension $n=2$, by Choi-Haslhofer-Heshkovitz-White \cite{CHHW} for dimension $\geq 3$ with a $2$-convexity assumption, and by Gang Zhou \cite{G1} for the singularity model $\mathbb{S}^1\times\mathbb{R}^3$ with a nondegeneracy assumption. 
	
	We prove that all the nondegenerate cylindrical singularities have a mean convex neighborhood.
	
	\begin{theorem}[Mean convex neighborhood of a nondegenerate singularity]\label{prop:Mconnbhd}
		Assume \eqref{eq:StarCondition} and that the singularity is nondegenerate. Then there exist constants $\dt_2>0$ and $t_0>0$, such that for any $t\in(-t_0,0)$, $\bM_t\cap B_{\dt_2}(0)$ has positive mean curvature and is diffeomorphic to $\cC_{n,k}\cap B_{4n}$. 
	\end{theorem}
	
	The mean convex neighborhood property is related to the nonfattening property of level set flow. In \cite{HershkovitsWhite20_Nonfattening}, Hershkovits-White proved that if all the singularities are of mean convex/mean concave type, then the level set flow does not fatten. The condition in Theorem \ref{prop:Mconnbhd} is sufficient to show that a nondegenerate cylindrical singularity is of such a type. Thus, if a mean curvature flow only has nondegenerate singularities, its level set flow does not fatten.

	Another crucial property of nondegenerate singularities is the type-I curvature condition. A singularity $(y,T)$ of a mean curvature flow is called \textbf{type-I} if the curvature blows up with a speed of at most $O((T-t)^{-1/2})$. This means that there exist constants $r>0$ and $C>0$ such that for any $t<T$, in $\bM_t\cap B_r(y)$, we have $|A|(x,t)=O((T-t)^{-1/2})$, where $A$ represents the second fundamental form. On the other hand, if the curvature blows up faster than this rate, the singularity is called \textbf{type-II}. There are examples of type-II singularities; for instance, in \cite{AAG}, the authors constructed a surface called a ``peanut," and the mean curvature flow starting from such a peanut develops type-II singularities known as ``degenerate neckpinching." In addition, Angenent-Vel\'azquez \cite{AV} constructed a large family of type-II singularities in their work.
	
	\begin{theorem}[Type-I curvature condition of a nondegenerate singularity]\label{ThmtypeI}
		Assume \eqref{eq:StarCondition} and that the singularity is nondegenerate. Then there exist constants $\dt_2>0$ and $t_0>0$, such that for any $t\in(-t_0,0)$, $\bM_t\cap B_{\dt_2}(0)$ has curvature bound $|A|\leq C(-t)^{-1/2}$ for some constant $C$.
	\end{theorem}
	
	\subsection{Subsequent work}
	On the dynamics side, in Part II of the paper in the sequel, we shall address the stability and genericity of nondegenerate singularities. On the geometry side, in a subsequent work with Zhihan Wang \cite{SunWangXue1_Passing}, we further explored some finer properties of nondegenerate singularities, based on the normal form theorems of this paper. Particularly, we proved that Theorems \ref{ThmIsolate}, \ref{prop:Mconnbhd}, \ref{ThmtypeI} hold not only in the backward parabolic neighborhood, but also hold in the forward parabolic neighborhood. Another geometric consequence of the above normal forms is the following profile of the flow at a nondegenerate singularity. 
	\begin{theorem}[\cite{SunWangXue1_Passing}, Section 5.1]
		Assume $(\star)$ and the singularity is nondegenerate. When $t=0$, there is an $r_0>0$ such that the MCF restricted to the ball $B_{r_0}(0)\subset \R^{n+1}$ is the cylindrical coordinate graph $\{y+f(\theta,y)\theta\}$ of a function $f$ such that \[
		f(\theta,y)=\frac{\sqrt{2(n-k)}\ |y|}{2\sqrt{-\log(|y|)}}(1+o_y(1)),  \] 
		where $\|o_y(1)\|_{C^0}\to 0$ as $y\to 0$.
	\end{theorem}
	We give the proof in \cite{SunWangXue1_Passing} in the joint work with Zhihan Wang. The result will not be used in this paper. The profile first appeared in \cite{AV} in the spherical symmetric setting and in \cite{gang2025precise} in the case of $k=3,n=4.$ The profile inspires us to obtain a complete description of the change of geometry and topology of the mean curvature flow passing through a nondegenerate singularity. In \cite{SunWangXue2_RegSing}, we further studied those singularities that are partially nondegenerate and degenerate, and a refined normal form enables us to show that the singular sets have $C^{2,\alpha}$ regularity.
	
	The preprint of this paper attracted lots of attention from experts. It initiated a program of establishing the isolatedness of singularities of generic mean curvature flows using normal form and nondegeneracy. In the later joint work with Zhihan Wang \cite{SunWangXue1_Passing}, we proved an important step in the program by establishing that the nondegenerate singularities are isolated in \emph{spacetime}, and that the mean curvature flow passes through a nondegenerate singularity with topology change the same as that in Morse theory when crossing a critical point. This program was later completed in $\R^3$ and for mean convex mean curvature flows in $\R^{n+1}$ with $n\geq 3$ \cite{szekelyhidi2026mean, szekelyhidi2026nondegenerate, daniels2026generic}. To complete the program for generic MCFs in the higher-dimensional case, one might need to establish the higher dimensional Huisken conjecture on generic tangent flow and the multiplicity $1$ conjecture first.

	\subsection{Strategy of the proof}
	Let us first comment on two important ingredients of the proof. The first is the \textbf{pseudolocality} theorem of mean curvature flow. It is a nonlinear effect that was first discovered by Ecker-Huisken in the setting of mean curvature flow. The pseudolocality shows that if a mean curvature flow is sufficiently close to a cylinder in a sufficiently large domain, then in a definite amount of time, it behaves similarly to the flow of a shrinking cylinder. In particular, when we work with the rescaled mean curvature flow, the pseudolocality shows that if a rescaled mean curvature flow is sufficiently close to a cylinder in a sufficiently large domain, then in a definite amount of time, it will be close to the cylinder inside a larger domain, and its higher order derivatives will have better control. This is the key to the nonlinear estimate to study the rescaled mean curvature flows, even if the nonlinear error is much more complicated than the semilinear equations.
	
	The second is Vel\'azquez's argument, which is also known as the \textbf{Ornstein-Uhlenbeck regularization}. For a nonnegative function $Z(\cdot,\tau)$ defined on $\cC_{n,k}$ satisfying the inequality
	\[
	\pr_\tau Z-LZ\leq \eps_0 Z+\text{terms decaying in time}+\text{terms from cut-off},
	\]
	this tool provides an estimate of $\|Z(\cdot,\tau')\|_{C^0}$ inside ball of radius $r(\tau)=O(\tau^{1/2})$ by $(r(\tau))^2\|Z(\cdot,\tau)\|_{L^2}$ for $\tau'>\tau$. 
	
	In applications, we will choose $Z=|\chi u|+|\nabla (\chi u)|$, where $u$ is the graph function of the rescaled mean curvature flow over $\cC_{n,k}$ and $\chi$ is a cut-off function; we will also choose $Z=|\chi w|+|\nabla (\chi w)|$, where $w$ is the difference of $u$ and the function given by the normal form. Then Vel\'azquez's argument provides, roughly speaking, a $C^1$-control of the graph by the $H^1$-norm. In particular, the $C^1$-control is the ingredient to apply the pseudolocality theorem. In order to get a very small $C^1$-bound, we need a decay of $\|Z(\cdot,\tau)\|_{L^2}$ faster than $(r(\tau))^{-2}$. Hence, the graphical radius relies on the $H^1$-norm of the graph function. That is the reason that we need to bootstrap between the $C^1$-norm and the $H^1$-norm.

	The proof is divided into the following steps. The initial input is the $L^2$-bound of the graphical function (decay like $\tau^{-\xi}$, $\xi>0$ that can be very small) and the graphical radius $\br(\tau)\geq C\sqrt{(\log \tau)}$ by \cite{CM2} obtained using the {\L}ojasiewicz inequality. 
	\begin{enumerate}
		\item We improve the graphical radius to $\br(\tau)\geq C\tau^{\kappa}$ for $2\kappa<\xi$ as follows: starting at some large time $T$, we use the pseudolocality of mean curvature flow to improve the rough graphical radius to $T+\delta$. The graph function $u$ can have very large $C^0$-norm, but $|\nabla^k u|$ is small for $k=1,2,3$. Then we use Vel\'azquez's argument to show that the graph function indeed has a small $C^1$-norm. Then we iteratively use the pseudolocality and Vel\'azquez's argument to extend the graphical radius to all time $\tau>T$.
		\item Using the methods of dynamical systems, we obtain an $H^1$-normal form over $B_{\br(\tau)}$ and improve the $H^1$-estimate to $\|u-\frac{\varrho}{4\tau}\sum (y_i^2-2)\|_{L^2}=O(\tau^{-1-\vartheta})$. The dynamical argument is similar to the authors' previous work \cite{SX1, SX2, SX3}. To derive the explicit normal form, we follow the idea of Vel\'azquez \cite{V1} in the study of semilinear PDE.
		\item Let $w=u-\varrho(\sqrt{1+\frac{\varrho}{2\tau}\sum (y_i^2-2)}-1)$. Then use the perturbed Ornstein-Uhlenbeck regularization to estimate $w$ (hence $u$). In particular, together with pseudolocality, we extend the graphical radius to $K_0\sqrt{\tau}$, for some small $K_0$. $K_0$ can only be a small number because the perturbed linearized operator should be close to the original linearized operator. 
		\item Once we get the $C^1$-normal form on $(1+\alpha) K_0\sqrt{\tau}$, we repeatedly using pseudolocality to prove $|\nabla^k u|$ bound on the ball of radius $(1+\alpha)^m K_0\sqrt{\tau}$ for $m\geq 1$ and $k=1,2,3$. Particularly, we have the $C^1$-normal form on the ball of radius $K\sqrt{\tau}$ for any $K>0$.
	\end{enumerate}

	\subsection{Organization of the paper}
	Section \ref{S_preliminaries} contains some preliminaries; in Section \ref{S_OU_Estimate}, we present the adaptation of Vel\'azquez's Ornstein-Uhlenbeck regularization; in Section \ref{S:Extending_t_kappa} we proved the first step, showing that the graphical radius is $O(\tau^\kappa)$ for some $\kappa\in(0,1/2)$; in Section \ref{SDynamics}, we use the dynamical argument to study the evolution of the $0$-eigenmodes, and prove the $H^1$-normal form; in Section \ref{S:C1-sqrtK}, we prove the $C^1$-normal form based on the $H^1$-normal form and the Ornstein-Uhlenbeck regularization; in Section \ref{S_Geometric} we study some geometric properties of the nondegenerate singularities. We also discuss some nonlinear estimates in the Appendix.

	\subsection*{Acknowledgment}
	We would like to thank Professor Bill Minicozzi for the stimulating discussions, Professor Natasa Sesum and Professor Toti Daskalopoulos for their interest and comments. We also thank Zhihan Wang for valuable suggestions and comments -- some of them led to new statements and arguments in the current version. J. X. is supported by NSFC grants (No. 12271285) in China, the New Cornerstone investigator program, and the Xiaomi endowed professorship of Tsinghua University. J.X. is a member of IAS this semester, who would like to thank IAS for the hospitality and the support from the Shiing-Shen Chern member fund. 
	
	LLM AI tools were used during revision but only for polishing language and checking mathematical correctness. The framework, ideas, and techniques are the same as the original version, which was purely written by the authors.
	
	\section{Preliminaries}\label{S_preliminaries}
	\subsection{Eigenvalues and Eigenfunctions of the $L$-operator}\label{SPre}
	In this section, we summarize some previously known results on cylindrical singularities. We will be working in the weighted Sobolev space. Given two functions $f,g$ defined on a hypersurface $\Sigma$, we define an inner product
	$
	\langle f,g\rangle_{L^2} =\int_{\Sigma}f(x)g(x)e^{-\frac{|x|^2}{4}}d\cH^n(x).
	$
	Then we define the weighted $L^2$-norm by
	$\|f\|_{L^2}=\langle f,f\rangle_{L^2}^{1/2}$, and the weighted $L^2$-space consists of function $f$ with $\|f\|_{L^2}<+\infty$. Similarly, we define the weighted higher Sobolev space $H^{k}$ with
	$
	H^{k}(\cC_{n,k}):=\left\{f: \sum_{i=0}^k \|\nabla^i f\|^2_{L^2}<+\infty\right\}.
	$
	Throughout the paper, for simplicity of notation, we use $\|\cdot\|$ to denote the $H^1$-norm if not otherwise mentioned. We also introduce the weighted area measure
	\[
	d\mu(x)=e^{-\frac{|x|^2}{4}}d\cH^n(x).
	\]
	
	Recall that the linearized operator on a shrinker is defined as
	$$L:=\cL+(|A|^2+1/2)=\Delta-\frac{1}{2}\langle x,\nabla\cdot\rangle+(|A|^2+1/2).$$
	In the special case that the shrinker is $\cC_{n,k}$, we have
	\begin{equation}\label{EqL}
		L_{\cC_{n,k}} u= \Delta_{\bS^{n-k}(\varrho)} u +L_{\R^k}u,\quad L_{\R^k}u=\Delta_{\R^{k}} u -\sum_{i=1}^{k}\frac{1}{2}y_i\pr_{y_i}u+u.	
	\end{equation}
	Throughout this paper, an eigenvalue $\lambda$ of an operator $-L$ is a number such that there exists a nonzero function $f$ satisfying $Lf+\lambda f=0$. The following fact was proved in \cite[Section 5.2]{SunWangZhou20_MinmaxShrinker}: \emph{Suppose the eigenvalues of $-\Delta_{\bS^{n-k}(\varrho)}$ are given by $\mu_0\leq \mu_1\leq \mu_2\leq \cdots$ with corresponding eigenfunctions $\phi_0,\phi_1,\phi_2,\cdots$, and suppose the eigenvalues of $L_{\R^k}$ on $\R^{k}$ are given by $\nu_0\leq \nu_1\leq \nu_2\leq \cdots$ with corresponding eigenfunctions $\psi_0,\psi_1,\psi_2, \cdots$, then
		the eigenvalues of $L_{\cC_{n,k}}$ are given by $\{\mu_i+\nu_j\}_{i,j=0}^\infty$ with corresponding eigenfunctions $\{\phi_i\psi_j\}_{i,j=0}^\infty$.}
	
	The spectrum of $-\Delta_{\bS^{n-k}(\varrho)}$ is $0,1/2,\frac{n-k+1}{n-k},\cdots$, and the eigenfunctions are known to be the restriction of homogeneous harmonic polynomials. The first several eigenfunctions are listed as follows: constant functions for eigenvalue $0$; $\theta_i$, the restriction of linear functions in $\R^{n-k+1}$ to $\bS^{n-k}(\varrho)$, for eigenvalue $1/2$; and $\theta_i^2-\theta^2_j,\cdots$ for eigenvalue $\frac{n-k+1}{n-k}$. The spectrum of $-L_{\R^k}$ on $\R^{k}$ is given by half integers $\frac{m}{2}-1,\ m=0,1,2,\ldots$, and the eigenfunctions for eigenvalue $\frac{m}{2}-1$ are given by $h_{m_1}(y_1)\cdots h_{m_k}(y_k)$ with $m_1+\ldots+m_k=m$, where $h_{m_i}(y_i)=c_{m_i}\tilde h_{m_i}(y_i/2)$, and $\tilde h_{m_i}$ are standard Hermite polynomials with $c_{m_i}=2^{-m_i/2}(4\pi)^{-1/4}(m_i!)^{-1/2}$ are normalizing factors such that $\|h_{m_i}\|_{L^2(\R)}=1$. In particular, we have $$\tilde h_0(x)=1,\ \tilde h_1(x)=2x,\ \tilde h_2(x)=4x^2-2 \mathrm{\ and\ } c_0=(4\pi)^{-1/4}, c_1=4^{-1/2}\pi^{-1/4}, c_2=2^{-2}\pi^{-1/4}.$$ We shall use the following fact (c.f. Appendix B of \cite{HV2}).
	\begin{lemma}\label{LmTriple}
		Let $A_{m,n,\ell}=\int_\R h_m(x)h_n(x)h_\ell(x)e^{-\frac{|x|^2}{4}}dx.$ Then $A_{m,n,\ell}=0$ unless 
		$m+n+\ell$ is even and $n\leq m+\ell,\ m\leq n+\ell,\ \ell\leq m+n$, in which case we have $$A_{m,n,\ell}=(4\pi)^{-1/4}(m!n!\ell!)^{1/2}\left(\left(\frac{m+n-\ell}{2}\right)!\left(\frac{n+\ell-m}{2}\right)!\left(\frac{m+\ell-n}{2}\right)!\right)^{-1}.$$ In particular, we have $A_{2,2,2}=2\pi^{-1/4}=8c_2.$
	\end{lemma}

	Combining the spectra together, we obtain the first three eigenvalues and their corresponding eigenfunctions of $-L_{\cC_{n,k}}$, see Table \ref{TableEigen}. These eigenfunctions have geometric meanings: 
	\begin{itemize}
		\item Constant $1$ is proportional to the mean curvature on the generalized cylinder, representing infinitesimal (spacetime) dilation.
		\item $\theta_i$ and $y_j$ are both infinitesimal translations. Specifically, $\theta_i$'s represent the translations in the directions of the spherical components, and $y_i$'s represent the translations in the axis directions.
		\item $\theta_i y_j$ represents the infinitesimal rotation.
		\item $h_2(y_i)=c_2(y_i^2-2)$ is known to be the non-integrable Jacobi field. It represents \emph{non-degenerate neckpinching} in the related direction on the axis. 
		\item $y_iy_j$'s show up if we rotate $h_2(y_i)$'s in the $\R^k$ space. e.g. $\left(\frac{y_i+y_j}{\sqrt2}\right)^2-2=\frac{y_i^2-2}{2}+\frac{y_j^2-2}{2}+y_iy_j$.
	\end{itemize}
	
	\begin{table}[H]
		\begin{tabular}{|l|l|}
			\hline
			eigenvalues of $-L_{\cC_{n,k}}$ & corresponding eigenfunctions \\ \hline
			$-1$ & $1$ \\ \hline
			$-1/2$ & $\theta_i,y_j,\ i=1,2,\ldots,n-k+1,\ j=1,2,\ldots,k$ \\ \hline
			$0$ & $\theta_iy_j,\ h_2(y_j)=c_2(y_j^2-2),\ y_{j_1}y_{j_2}$ \\ \hline
			$\min\{1/(n-k),1/2\}$ & $\ldots$ 
			\\ \hline
		\end{tabular}
		\caption{Eigenvalues and eigenfunctions of $-L_{\cC_{n,k}}$.}
		\label{TableEigen}
	\end{table}
	
	Note that while $h_{m_i}(y)$ is a normalized eigenfunction in $\|\cdot\|_{L^2(\R)}$, $h_{m_i}(y_i)$ is not a normalized eigenfunction in $\|\cdot\|_{L^2(\cC_{n,k})}$. Thus, we define $H_{m_1,m_2,\cdots,m_k}$ be a multiple of $h_{m_1}h_{m_2}\cdots h_{m_k}$, such that $H_{m_1,m_2,\cdots,m_k}$ is a normalized eigenfunction in $\|\cdot\|_{L^2(\cC_{n,k})}$. In particular, $H_2(y_i)=c_0^{k-1}\cG_{n,k}^{-1/2}h_2(y_i)$, $H_{1,1}(y_i,y_j)=c_0^{k-2}\cG_{n,k}^{-1/2}h_1(y_i)h_1(y_j)$. Here \begin{equation}\label{EqcG}
		\cG_{n,k}=\int_{\mathbb S^{n-k}(\varrho)}e^{-\frac{|x|^2}{4}}d\cH^{n-k}=\varrho^{n-k}e^{-\frac{\varrho^2}{4}}\omega_{n-k},        
	\end{equation}
	where $\omega_{n-k}$ is the area of the $(n-k)$-dimensional unit sphere.

	\subsection{Coarse graphical scale and $L^2$ estimate}\label{SSCM}
	
	By Colding-Minicozzi's {\L}ojasiewicz inequality \cite{CM2}, if one tangent flow of a bounded entropy mean curvature flow $\bM_{t}$ at $(0,0)$ is $\cC_{n,k}$, then \eqref{eq:StarCondition} holds. In \cite{CM2}, Colding-Minicozzi introduced the notion of \emph{cylindrical scale}.	
	
	\begin{definition}[Graphical radius; called cylindrical scale in \cite{CM2}]\label{DefGraphRad} Let $\eps_0>0,\ \ell\geq 2$ and $C_{\ell}$ be some fixed constants.
		Given a hypersurface $\Sigma$, the cylindrical scale $($throughout the paper, we call it the \emph{graphical radius}$)$ $\br(\Sigma)$ is the largest radius such that $\Sigma \cap B_{\br(\Sigma)}$ is the graph of a function $u:\ \cC_{n,k}\cap B_{\br(\Sigma)}\to \R$ with $\|u\|_{C^{2,\alpha}}\leq \eps_0$ and the part of hypersurface $\Sigma \cap B_{\br(\Sigma)}$ has curvature bound $|\nabla^\ell A|\leq C_{\ell}$. 
	\end{definition}
	
	Because we will work with various graphical regions, and the behavior of the rescaled mean curvature flow can be different in different regions, we may \emph{not} always choose $\br$ to be the largest radius.
	
	We need some results from \cite{CM2} to serve as initial input for our proofs. The result is summarized in the following proposition. 
	\begin{proposition}\label{PropR}
		Assume \eqref{eq:StarCondition} and let $\br(\tau)$ be the graphical radius of the rescaled mean curvature flow $M_\tau$. Then there exists $\tau_0$ such that we have the estimate $\br(\tau)\geq (\al\log \tau)^{1/2}$ for some small $\al>0$ and all $\tau>\tau_0$. Moreover, if $r(\tau)=  (\al\log \tau)^{1/2}$, the graphical function $u:\ \cC_{n,k}\cap B_{r(\tau)}\to \R$ has the estimate  $\|u(\cdot,\tau)\|_{H^1(B_{r(\tau)})}\leq C\tau^{-\frac{1}{3}}$ for $\tau>\tau_0$. 
	\end{proposition}
	\begin{proof} 
		In Theorem 6.1 of \cite{CM2} we choose the $\tau$ therein, denoted by $\tau_{CM}$, to be $1/(3-\epsilon)$, then we get the estimate $F(\Sigma_{\tau-1})-F(\Sigma_{\tau+1})\leq \tau^{-(3-\epsilon)}$ by Lemma 6.9 of \cite{CM2} (we remark that the statement of Theorem 6.1 of \cite{CM2} gives only some $\tau_{CM}\in (1/3,1)$. However, the proof there is constructive. It is not hard to see that $\tau_{CM}$ can be chosen close to $1/3$. See \cite[Footnote 10]{CM2} on page 244). 
		Defining $R(\tau)$ (shrinker scale in \cite{CM2}) through $e^{-\frac{R(\tau)^2}{2}}=F(\Sigma_{\tau-1})-F(\Sigma_{\tau+1})$,
		we get $R(\tau)\ge \sqrt 2\sqrt{\frac{1}{\tau_{CM}}\log \tau} .$
		By Theorem 5.3 of \cite{CM2} on the graphical scale $\br(\tau)\geq (1+\mu)R(\tau)$ for some $\mu>0$, we get the estimate $$\br(\tau)\geq \frac{\sqrt2}{\sqrt{\tau_{CM}}}(1+\mu)\sqrt{\log\tau}=\sqrt{2(3-\eps)} (1+\mu)\sqrt{\log \tau}.$$ 
		
		To estimate $\|u\|_{H^1(B_{r(\tau)})}$, we first consider the situation where $r(\tau)=  (\al\log \tau)^{1/2}$. We integrate the square of \cite[(2.51)]{CM2} with respect to the Gaussian area to yield (with $r(\tau)$ possibly getting smaller by a constant)
		$$\|u(\tau)\|^2_{H^1(B_{r(\tau)})}\leq C r(\tau)^\rho(\|\phi\|_{L^1(B_{r(\tau)})}^{1-\eps}+e^{-(1-\epsilon)\frac{r(\tau)^2}{4}})
		\leq C(\al\log \tau)^{\rho/2}
		(\|\phi\|^{1-\eps}_{L^1(B_{r(\tau)})}+
		\tau^{-\frac{3}{2}(1+\mu)^2(1-\eps)}
		)
		,$$
		for some uniform constants $C$ and $\rho$, and arbitrarily small $\epsilon$,
		where $\phi$ is a function with the estimate (see equation (6.6) of \cite{CM2})
		$$\|\phi\|^4_{L^1(B_{r(\tau)})}\leq C \|\phi\|^2_{L^2(B_{r(\tau)})}\leq C(F(\Sigma_{\tau-1})-F(\Sigma_{\tau+1}))\leq \tau^{-1/\tau_{CM}}.$$
		Combining the two estimates, we get the estimate of $\|u\|_{H^1(B_{r(\tau)})}\leq C(\log\tau)^{\rho/4}\tau^{-3/8}\leq  \tau^{-1/3}$. 	\end{proof}
	We remark that the graphical radius growth rate and the decay rate of $u$ are uniform in entropy for all cylindrical singularities (see footnote 5 of \cite{CM2} on page 223). 
	
	\subsection{Pseudolocality}\label{SSPL}
	The following pseudolocality argument was first studied by Ecker-Huisken in \cite{EH}. Later, Ilmanen-Neves-Schulze \cite[Section 9]{INS} gave a simple proof using Brakke-White regularity theorem \cite{Wh4}. There, they studied the mean curvature flow of graphs over a hyperplane, and we adapt their proof to study the mean curvature flow of graphs over a cylinder.
	
	\begin{theorem}\label{thm:pseudo-locality}
		For any $\eps\in(0,1)$, $\rho\in(0,1)$ and $\sigma>0$, there exist $R_0>0$ and $\eta_0>0$ with the following significance: suppose $M_\tau$ is a rescaled mean curvature flow and $M_{\tau_0}$ is the graph of a function $u(\cdot,\tau_0)$ over $\cC_{n,k}\cap B_R(p)$ for some $R>R_0$, $p\in\{0\}\times\R^k$, and $\|u(\cdot,\tau_0)\|_{C^1(B_R(p))}\leq\eta\leq \eta_0$. Then there exists $\bar\delta >0$ such that when $\tau\in[\tau_0+\sigma,\tau_0+\bar\delta]$, $M_\tau$ is the graph of a function $u(\cdot,\tau)$ over $\cC_{n,k}\cap B_{e^{(\tau-\tau_0)/2}\rho R}(e^{(\tau-\tau_0)/2}p)$ with $$\|\nabla^ku(\cdot,\tau)\|_{L^\infty(\cC_{n,k}\cap B_{e^{(\tau-\tau_0)/2}\rho R}(e^{(\tau-\tau_0)/2}p))}\leq e^{-(k-1)(\tau-\tau_0)/2}\eps, \quad k=0,1,2,3.$$ Moreover, $\bar\delta\to\infty$ as $\eta\to 0$.
		
	\end{theorem}
	
	In the theorem, the factor $e^{-(k-1)(\tau-\tau_0)/2}$ is from the exponential expansion property of the rescaled mean curvature flow. In particular, when $k=0$, the $L^\infty$ norm of $u$ increases very fast. We need to use the equation of the rescaled mean curvature flow $\partial_\tau u=Lu+\cQ$ to show that the growth rate is actually not that fast.
	
	Theorem \ref{thm:pseudo-locality} follows from the following theorem for mean curvature flow, considering the scaling between mean curvature flow and rescaled mean curvature flow. 
	
	\begin{theorem}\label{thm:pseudo-locality-MCF}
		For any $t_0\in(-1,0)$, $\sigma>0$, $\rho\in(0,1)$, and $\eps\in(0,1)$, there exists $\eta_0>0$ such that for any $\eta\in(0,\eta_0)$, if $\bM_{-1}\cap B_{\eta^{-1}}$ is the graph of a function $v$ over $\cC_{n,k}$ with $\|v\|_{C^1}\leq \eta$, then 
		\begin{enumerate}
			\item $\bM_{t}\cap B_{\rho\eta^{-1}}$ is smooth mean curvature flow for $t\in(-1,t_0)$;
			\item $\bM_{t}\cap B_{\rho\eta^{-1}}$ has a singularity before time $|t_0|$;
			\item $\bM_{t}\cap B_{\rho\eta^{-1}}$ can be written as the graph of a function $u(\cdot,\tau)$ over $\sqrt{-t}\cC_{n,k}$ for $t\in(-1+\sigma,t_0)$, with $\|u(\cdot,t)\|_{C^3}\leq \eps$.
		\end{enumerate}
	\end{theorem}
	
	\begin{proof}
		We prove by contradiction. Suppose there is no $\eta_0$ such that all of the three items hold, then we can find a sequence of hypersurfaces $\bM_{-1}^k$, such that each one of them is the graph of a function $v_k$ inside the ball of radius $k$, with $\|v_k\|_{C^1}\leq k^{-1}$, but at least one of the three items fails. Then we can use the compactness of unit regular cyclic Brakke flow (see \cite[Section 7]{I2} and \cite{White09_CurrentsVarifolds}) to pass to a limiting weak mean curvature flow $\bM_t^{\infty}$. By Brakke-White's regularity theorem \cite{Wh4}, and Ecker-Huisken's curvature estimate, $\bM_t^{\infty}$ is a smooth mean curvature flow when $\tau$ is close to $-1$. Also, the sequence $\bM_{-1}^k$ has varifold limit being the standard cylinder $\cC_{n,k}$. Thus, the limit flow $\bM_t^{\infty}$ is exactly the shrinking soliton, with only one unique singular time $0$. Then, by the upper semi-continuity of the Gaussian density, we know that the singularity of $\bM_{\tau}^k$ must have the singular time converging to $0$, and $\bM_{t}^k$ are graphs over $\sqrt{-t}\cC_{n,k}$ with $C^3$-norm converging to $0$. Then all three items hold, which is a contradiction.
	\end{proof}
	
	\begin{proof}[Proof of Theorem \ref{thm:pseudo-locality}]
		Let $\sigma'=\frac12(1-e^{-\sigma})$. For $j\in\mathbb Z_+$, let $\eta_j$ be the constant given by Theorem \ref{thm:pseudo-locality-MCF} with $(t_0,\sigma,\rho,\eps)$ replaced by $(-e^{-j},\sigma',\frac12,\eps)$. Decreasing $\eta_j$ if necessary, we may assume $\eta_j<\frac{1}{2\varrho}$ and $\eta_j\searrow 0$. Fix an integer $j_0>\sigma+1$ and let $\eta_0=\eta_{j_0}$, $R_0=\frac{2}{(1-\rho)\eta_{j_0}}$.
		
		Given $R>R_0$ and $\eta\leq\eta_0$, let $\eta_1=\max\{\eta,\frac{2}{(1-\rho)R}\}\leq\eta_{j_0}$, and let $\bar\delta+1$ be the largest $j$ with $\eta_1\leq\eta_j$. Then $\bar\delta>\sigma$, and $\bar\delta\to\infty$ as $\eta_1\to0$. Define the mean curvature flow
		\[
		\bM_t:=\sqrt{-t}\,M_{\tau_0-\log(-t)}-p,\quad t\in[-1,0),
		\]
		so that $M_\tau=e^{(\tau-\tau_0)/2}(\bM_{-e^{-(\tau-\tau_0)}}+p)$, and $\bM_{-1}$ is a graph over $\cC_{n,k}\cap B_R$ with $C^1$-norm at most $\eta$.
		
		Since $\cC_{n,k}$ is invariant under translations in $\{0\}\times\R^k$, for any $q\in\{0\}\times\R^k$ with $|q|<\rho R$, the translate $\bM_{-1}-q$ is a graph over $\cC_{n,k}\cap B_{\eta_1^{-1}}$ with $C^1$-norm at most $\eta_{\bar\delta+1}$, because $\eta_1^{-1}+|q|<R$. Theorem \ref{thm:pseudo-locality-MCF} applied to $\bM_t-q$ shows that for $t\in(-1+\sigma',-e^{-\bar\delta-1})$, $\bM_t\cap B_{\eta_1^{-1}/2}(q)$ is a graph over $\sqrt{-t}\,\cC_{n,k}$ with $C^3$-norm at most $\eps$. As $\eta_1^{-1}/2>\varrho$, taking $q=(0,y)$ for $(\theta,y)\in \sqrt{-t}\,\cC_{n,k}\cap B_{\rho R}$, these balls cover $\sqrt{-t}\,\cC_{n,k}\cap B_{\rho R}$. Graph functions agree on overlaps, so $\bM_t\cap B_{\rho R}$ is the graph of some $\hat u(\cdot,t)$ with $\|\hat u(\cdot,t)\|_{C^3}\leq\eps$.
		
		Now let $\tau\in[\tau_0+\sigma,\tau_0+\bar\delta]$ and $\lambda=e^{(\tau-\tau_0)/2}$. Then $t=-\lambda^{-2}\in[-e^{-\sigma},-e^{-\bar\delta}]\subset(-1+\sigma',-e^{-\bar\delta-1})$. Since $\lambda\sqrt{-t}=1$, the dilation by $\lambda$ maps $\sqrt{-t}\,\cC_{n,k}$ to $\cC_{n,k}$ and preserves the unit normal. Hence $M_\tau=\lambda(\bM_t+p)$ is the graph over $\cC_{n,k}\cap B_{\lambda\rho R}(\lambda p)$ of
		\[
		u(z,\tau)=\lambda\,\hat u\big(\lambda^{-1}(z-\lambda p),t\big),
		\]
		and $\nabla^ku(z,\tau)=\lambda^{1-k}(\nabla^k\hat u)(\lambda^{-1}(z-\lambda p),t)$ gives the estimate for $k=0,1,2,3$.
	\end{proof}

	We shall also need a version of pseudolocality near cylinders whose radius is evolving, as this is the situation for the RMCF at the nondegenerate singularities away from the origin. For $b>0$ we say that a hypersurface is a normal graph of a function $w:\cC_{n,k}\to\R$ over $\bS^{n-k}(b)\times\R^k$ if it is given by $\{(\frac{b+w(\theta,y)}{\varrho}\theta,y):(\theta,y)\in\cC_{n,k}\}$ locally; the norms of $w$ are taken on $\cC_{n,k}$. Recall that the cylinder $\bS^{n-k}(b)\times\R^k$ moves by mean curvature flow with radius $b(t)=\sqrt{b^2-\varrho^2t}$.
	
	\begin{lemma}\label{lem:pseudo-radius}
		Let $0<b_-\leq b_+$, $0<\sigma\leq t_1$ with $\varrho^2t_1<b_-^2$, and $\eps\in(0,1)$. There exist $R_1,\eta_1>0$, depending only on $n,k,b_\pm,\sigma,t_1,\eps$ and the entropy bound, with the following significance: suppose $\{\bM_t\}_{t\in[0,t_1]}$ is a mean curvature flow, $b\in[b_-,b_+]$, $p\in\{0\}\times\R^k$, $R>R_1$, and $\bM_0\cap B_R(p)$ is a normal graph over $\bS^{n-k}(b)\times\R^k$ of a function $w_0$ with $\|w_0\|_{C^1}\leq\eta_1$. Then for every $t\in[\sigma,t_1]$, $\bM_t\cap B_{R/2}(p)$ is a normal graph over $\bS^{n-k}(b(t))\times\R^k$ of a function $w_t$ with $\|w_t\|_{C^3}\leq\eps$.
	\end{lemma}
	
	\begin{proof}
		Let $\lambda=\varrho/b\in[\varrho/b_+,\varrho/b_-]$ and $\hat\bM_t:=\lambda(\bM_{\lambda^{-2}(1+t)}-p)$, which is a mean curvature flow. Since $\lambda(\bS^{n-k}(b)\times\R^k)=\cC_{n,k}$, $\hat\bM_{-1}\cap B_{\lambda R}$ is the graph over $\cC_{n,k}$ of $\hat w_0(\theta,y)=\lambda w_0(\theta,p+\lambda^{-1}y)$, and $\|\hat w_0\|_{C^1}\leq C\eta_1$. Here and below $C$ depends only on $n,k,b_\pm,\sigma,t_1$. Set $\sigma'=\frac12\varrho^2\sigma/b_+^2$ and $t_0=\frac12(\varrho^2t_1/b_-^2-1)<0$, so that $[\lambda^2\sigma-1,\lambda^2t_1-1]\subset(-1+\sigma',t_0)$. We apply Theorem \ref{thm:pseudo-locality-MCF} with $\rho=\frac12$, $\sigma'$, $t_0$ and $\eps'=\eps/C$. If $R_1$ is large and $\eta_1$ is small, then $\hat\bM_t\cap B_{\lambda R/2}$ is a graph over $\sqrt{-t}\,\cC_{n,k}$ with $C^3$-norm at most $\eps'$ for $t\in[\lambda^2\sigma-1,\lambda^2t_1-1]$. Scaling back, $\bM_t\cap B_{R/2}(p)$ is a normal graph for $t\in[\sigma,t_1]$ over the cylinder of radius $\lambda^{-1}\varrho\sqrt{1-\lambda^2t}=b(t)$. The change of normalization from $\sqrt{-t}\,\cC_{n,k}$ to $\cC_{n,k}$ changes the $C^3$-norm by a factor at most $C$, since $-t$ is bounded above and below on this range. This gives $\|w_t\|_{C^3}\leq\eps$.
	\end{proof}

	\subsection{RMCF graph equation and cut-off}	\label{SSCone}
	To study the asymptotics of the rescaled mean curvature flow, we need to study the evolution of the graphical function. In Proposition \ref{PropQ} in Appendix \ref{appendix:Nonlinear estimate}, we write the graph of the rescaled mean curvature flow over a shrinker given by
	\[
	\pr_\tau u=Lu+\cQ(J^2u)
	\]
	where $u(\cdot, \tau):\ \cC_{n,k}\cap B_{\mathbf r(\tau)}\to \R$ whose graph is the rescaled mean curvature flow $M_\tau$ restricted to the ball $B_{\mathbf r(\tau)},$ and $\cQ$ is at least quadratic in $u$ given in \eqref{EqQ} and $J^2u:=(u,Du,D^2u)$ means the $2$-jet of $u$. We shall study the evolution of $u$ under the differential equation. 
	
	Since $M_\tau$ cannot be written as a global graph over $\cC_{n,k}$, we introduce a smooth cutoff function $\tilde \chi(\tau):\ \R_{\geq 0}\to \R$ that is $1$ over $[0,1]$ and vanishes outside $[0,2]$, so that $|\nabla^\ell \tilde \chi|\leq C_{\ell}$ for any $\ell\in\{0,1,2\}$. Then for any $r>1$ we define $\chi_r:\cC_{n,k}\to\R$ by 
	\[
	\chi_r(\theta,y)=
	\begin{cases}
		1,& |y|\leq r,
		\\
		\tilde \chi(|y|-r+1), & |y|\in[r,r+1],
		\\
		0,& |y|\geq r+1.
	\end{cases}
	\]
	In particular, given a differentiable function $f:\R\to\R_{\geq 1}$, we have $\pr_\tau \chi_{f(\tau)}(\theta,y)$ is only nonzero in $[f(\tau),f(\tau)+1]$, and $|\pr_\tau \chi_{f(\tau)}(\theta,y)|\leq C_1|f'(\tau)|$. In particular, when $f(\tau)=C(\log \tau)^{\alpha}$ or $\tau^\kappa$ for $\alpha>0$ and $\kappa<1$, we have $|\pr_\tau \chi_{f(\tau)}(\theta,y)|\to 0$ as $\tau\to\infty$.
	
	Given a function $\mathbf{r}:\R\to\R_{\geq 0}$, we define $\mathbb A_{\mathbf r(\tau)}$, or simply $\A_\tau$ if $\mathbf{r}$ is the graphical function and is fixed, to be the annulus region $\cC_{n,k}\cap(B_{\br(\tau)}\backslash B_{\br(\tau)-1})$. Then we derive the equation for $\chi u$:
	\begin{equation}\label{EqCutOff}
		\pr_\tau (\chi u)
		=
		L(\chi u)+\chi\cQ(J^2u)+(u(-\Delta\chi+\pr_\tau\chi)
		-2\langle \nabla\chi,\nabla u\rangle 
		+\frac{1}{2}\langle x,\nabla \chi\rangle u),
	\end{equation}
	where the last term on the RHS is supported on the annulus $\mathbb A_{\mathbf r(\tau)}$. For simplicity, we write
	\begin{equation}\label{Eq}
		\pr_\tau (\chi u)
		=
		L(\chi u)+\cB(J^2u).
	\end{equation}
	
	\section{Vel\'azquez's Ornstein-Uhlenbeck regularization estimate}\label{S_OU_Estimate}
	
	One of the key technical ingredients in this paper is the Ornstein--Uhlenbeck regularization studied by Vel\'azquez in \cite{V1}; see also \cite{HV1,HV2}. 
	
	Fix $\kappa\in(0,1/2]$ and $K>0$, and, for each target time $\tau$, define
	\begin{equation}\label{eq:OU-sliding-window}
		R(\tau)=K\tau^\kappa,\qquad
		H(\tau)=h_*+2\log(2+R(\tau)),\qquad
		s_\tau=\tau-H(\tau),
	\end{equation}
	where the fixed constant $h_*>4$ is chosen larger than the waiting-time constant in Lemma \ref{lem:perturbed-kernel} below. Thus
	\begin{equation}\label{eq:OU-window-properties}
		e^{-H(\tau)/2}R(\tau)\leq e^{-h_*/2},
		\qquad
		e^{H(\tau)}=e^{h_*}(2+R(\tau))^2.
	\end{equation}
	
	\begin{proposition}[OU regularization]
		\label{PropVelazquez}
		Let $0<\xi<2$, $\eta>0$, and $C_0,C_F>0$.  There are
		$\eps_*>0$, $c_*=c_*(\kappa)>0$, and $T_0>0$ such that the
		following holds.  Suppose $0\leq\eps_0\leq
		\min\{\eps_*,c_*\eta\}$, and $Z\geq0$ satisfies the following  on the interval 
		$[s_\tau,\tau]\subset[T_0,\infty)$: 
		\begin{equation}\label{EqVelazInequality}
			\partial_rZ-LZ
			\leq\eps_0Z+C_F\left(
			\frac{1+|y|^2}{r^2}
			+|y|\check\chi(y,r)\right),
		\end{equation}
		where \begin{equation}\label{eq:OU-cutoff-buffer}
			0\leq\check\chi\leq1,\qquad
			\check\chi(y,r)=0\quad\hbox{when }|y|\leq2R(r).
		\end{equation}
		Assume
		\begin{equation}\label{eq:OU-L2-assumption}
			\|Z(\cdot,r)\|_{L^2}
			\leq C_0r^{-\xi},
			\qquad s_\tau\leq r\leq\tau.
		\end{equation}
		Then, for every sufficiently large $\tau$,
		\begin{equation}\label{eq:OU-loss-estimate}
			\begin{split}
				\|Z(\cdot,\tau)\|_{C^0(B_{R(\tau)})}
				\leq C_\eta C_0\tau^{2\kappa-\xi+\eta}.
			\end{split}
		\end{equation}
	\end{proposition}

	While the linear effect of the drift Laplacian is the same, our problem is very different from the semilinear heat equation $\partial_\tau u=\Delta u+u^p$ considered by Vel\'azquez due to the nonlinear part. In particular, our nonlinearity contains second-order derivatives. In later proofs, we shall recast our problem into the form of the differential inequality \eqref{EqVelazInequality} that differs from the differential inequality considered in \cite{V1}. Moreover, our $L$-operator contains a spherical Laplacian summand.

	Let $S(h)$ be the positive semigroup generated by $L$. Since the coefficient of $Z$ in \eqref{EqVelazInequality} is the constant
	$\eps_0$, variation of constants absorbs it exactly:
	\begin{equation}\label{eq:OU-Duhamel-comparison}
		\begin{split}
			Z(\tau)\leq{}&
			e^{\eps_0H(\tau)}S(H(\tau))Z(s_\tau)+C_F\int_{s_\tau}^{\tau}
			e^{\eps_0(\tau-r)}S(\tau-r)
			\left(\frac{1+|y|^2}{r^2}
			+|y|\check\chi(y,r)\right)\,dr .
		\end{split}
	\end{equation}
	We have the following estimate bounding the growth given by the $\eps_0Z$ term
	\[
	e^{\eps_0H(\tau)}
	\leq C(2+R(\tau))^{2\eps_0}
	\leq C_\eta\tau^\eta
	\]
	when $\eps_0\leq c_*(\kappa)\eta$.  
	
	A particular application of this section is the following $C^0$-estimate of the linearized rescaled mean curvature flow equation, Proposition \ref{PropBoundary-global}. Although we will prove stronger estimates for nonlinear equations, the linear analysis may have independent interest, and we include it here.
	
	\begin{proposition}[\cite{V1}; see also Proposition 2.13 of \cite{KZ}]
		\label{PropBoundary-global}
		Let $v_{s_\tau}\in L^2$, and let
		\[
		v(r)=S(r-s_\tau)v_{s_\tau},
		\qquad s_\tau\leq r\leq\tau,
		\]
		be a  solution to $\partial_rv=Lv$.  Then
		\begin{equation}\label{eq:linear-OU-correct}
			\|v(\cdot,\tau)\|_{C^0(B_{K\tau^\kappa})}
			\leq C(1+K\tau^\kappa)^2
			\|v_{s_\tau}\|_{L^2}.
		\end{equation}
	\end{proposition}

	\subsection{The Mehler kernel and moving Gaussian norms}\label{SKZ}
	
	We now implement the heat-kernel argument of Vel\'azquez on the interval $[s_\tau,\tau]$. The semigroup $S(h)$ generated by $L=L_{\cC_{n,k}}$ has kernel
	\begin{equation}\label{EqS}
		S(h,(\theta,y),(\eta,z))
		=G(h,\theta,\eta)\frac{e^h}
		{(4\pi(1-e^{-h}))^{k/2}}
		\exp\left(-\frac{|e^{-h/2}y-z|^2}
		{4(1-e^{-h})}\right),
	\end{equation}
	where $G$ is the heat kernel of $\Delta_{\bS^{n-k}(\varrho)}$. Following \cite{V1}, we introduce the shifted center $L^p$-norm
	\begin{equation*}\label{eq:moving-Gaussian-norm}
		N_r^q(\psi)=\sup_{|\zeta|\leq r}
		\left(\int_{\cC_{n,k}}|\psi(\theta,y)|^q
		e^{-|y-\zeta|^2/4}\,d\theta\,dy\right)^{1/q},
		\qquad N_r=N_r^2.
	\end{equation*}
	In particular, $N_0(\psi)=\|\psi\|_{L^2}$.
	
	\begin{lemma}\label{lem:estNr}
		For $h>0$ and $\psi\in L^2$,
		\begin{equation}\label{eq:Mehler-L2-pointwise}
			|S(h)\psi(\theta,y)|
			\leq C_h e^h
			\exp\left(\frac{|y|^2}{4(e^h+1)}\right)\|\psi\|_{L^2},
		\end{equation}
		where $C_h$ is bounded for $h\geq1$. Consequently, for every $A<\infty$,
		\begin{equation}\label{eq:estNr}
			N_r(S(h)\psi)\leq C_Ae^h\|\psi\|_{L^2}
		\end{equation}
		whenever $h\geq1$ and $e^{-h/2}r\leq A$. Moreover, for every fixed $d>0$,
		\begin{equation}\label{eq:fixed-time-N-to-C0}
			\|S(d)\psi\|_{C^0(B_R)}
			\leq C_dN_{e^{-d/2}R}(\psi).
		\end{equation}
	\end{lemma}
	
	\begin{proof}
		Applying Cauchy--Schwarz in \eqref{EqS}, with the factor $e^{-|z|^2/4}$ assigned to $\psi^2$, gives
		\[
		|S(h)\psi(\theta,y)|
		\leq e^h\|\psi\|_{L^2}
		\left(\int_{\cC_{n,k}}G(h,\theta,\eta)^2
		\left[\frac{1} {(4\pi (1-e^{-h}))^k} 
		\exp\left(-\frac{|z-e^{-h/2}y|^2}{2(1-e^{-h})}+\frac{|z|^2}{4}\right)\right]
		\,d\eta\,dz\right)^{1/2}.
		\]
		Completing the square in $z$ shows that the term in the bracket equals
		$
		C_h\exp\left(\frac{|y|^2}{4(e^h+1)}\right).
		$
		The $L^2$-norm of the spherical heat-kernel $G$ is bounded uniformly for $h\geq1$. This proves
		\eqref{eq:Mehler-L2-pointwise}.
		
		Squaring \eqref{eq:Mehler-L2-pointwise}, multiplying by
		$e^{-|y-\zeta|^2/4}$, and completing the square once more gives
		\[
		\int_{\cC_{n,k}}|S(h)\psi|^2e^{-|y-\zeta|^2/4}
		\leq C e^{2h}
		\exp\left(\frac{|\zeta|^2}{2(e^h-1)}\right)N_0(\psi)^2
		\]
		for $h\geq1$. Taking the supremum over $|\zeta|\leq r$
		proves \eqref{eq:estNr}.
		
		Finally, repeat the first Cauchy--Schwarz argument with the Gaussian weight centered at $m=e^{-d/2}y$. Since $d>0$ is fixed, the remaining Gaussian integral is bounded by $C_d$, and hence
		\[
		|S(d)\psi(\theta,y)|\leq C_dN_{|e^{-d/2}y|}(\psi).
		\]
		Taking the supremum over $|y|\leq R$ proves
		\eqref{eq:fixed-time-N-to-C0}.
	\end{proof}
	
	Applying the last lemma to the time interval given by \eqref{eq:OU-sliding-window}, we get the following estimate. 
	\begin{corollary}\label{LmC0byL2}
		Let $H=H(\tau)$ be defined by \eqref{eq:OU-sliding-window}. Then
		\begin{equation}\label{EqC0byL2}
			\|S(H)\psi\|_{C^0(B_{R(\tau)})}
			\leq Ce^H\|\psi\|_{L^2}
			\leq C(1+R(\tau))^2\|\psi\|_{L^2}.
		\end{equation}
	\end{corollary}

	\subsection{Polynomial moments and the separated cutoff tail}\label{SSOU}
	We have understood the behavior of the propagator $S(h)$ in the last section. In this section, we show how to estimate the RHS of \eqref{eq:OU-Duhamel-comparison}. It is convenient to write $S(h)=e^hP_h$, where $P_h$ is the semigroup operator generated by
	$\Delta_{\cC_{n,k}}-\frac12y\cdot\nabla_y$.
	
	\begin{lemma}\label{lem:OU-source-estimates}
		For all $h\geq0$, there is $C>0$ such that
		\begin{equation}\label{eq:OU-polynomial-moment}
			S(h)(1+|z|^2)(\theta,y)\leq C(e^h+|y|^2).
		\end{equation}
		Moreover, there are $C,c,R_0>0$ such that, if $R\geq R_0$,
		$|y|\leq R$, and $D\geq3R/2$, then
		\begin{equation}\label{eq:OU-separated-tail}
			S(h)\bigl(|z|\mathbf1_{\{|z|\geq D\}}\bigr)(\theta,y)
			\leq Ce^he^{-cR^2}.
		\end{equation}
	\end{lemma}
	
	\begin{proof}
		The Euclidean part of $P_h$ is a Gaussian with mean $e^{-h/2}y$ and covariance $2(1-e^{-h})I$. Therefore
		\[
		P_h|z|^2=e^{-h}|y|^2+2k(1-e^{-h}),
		\]
		which proves \eqref{eq:OU-polynomial-moment}.
		
		For \eqref{eq:OU-separated-tail}, fix $\mu=1/32$. Direct Gaussian integration gives
		\begin{equation}\label{eq:OU-exp-moment}
			P_he^{\mu|z|^2}(\theta,y)
			=(1-4\mu(1-e^{-h}))^{-k/2}
			\exp\left(
			\frac{\mu e^{-h}|y|^2}{1-4\mu(1-e^{-h})}\right).
		\end{equation}
		For $\nu=3\mu/4$,
		\[
		P_h (r\mathbf1_{\{r\geq D\}})\leq C_\mu
		e^{-\nu D^2}e^{\mu r^2}.
		\]
		Since $D\geq3R/2$ and $|y|\leq R$, the exponent obtained from \eqref{eq:OU-exp-moment} is at most $-cR^2$, uniformly in $h$. Multiplication by $e^h$ proves \eqref{eq:OU-separated-tail}.
	\end{proof}
	
	\subsection{Proof of the unperturbed estimates}\label{SSVelaz}

	\begin{proof}[Proof of Proposition \ref{PropVelazquez}]
		Write $R=R(\tau)$, $H=H(\tau)$, and $s=s_\tau$. Since $H=O(\log\tau)$, after increasing $T_0$ we have
		\begin{equation}\label{eq:OU-window-comparability}
			s\geq\tau/2,\qquad
			R(r)\geq\frac34R(\tau)\quad(s\leq r\leq\tau).
		\end{equation}
		
		If $|y|\leq R$, Lemma \ref{LmC0byL2} and \eqref{eq:OU-L2-assumption} give
		\begin{equation}\label{eq:OU-initial-term}
			e^{\eps_0H}S(H)Z(s)(\theta,y)
			\leq CC_0e^{\eps_0H}(1+R)^2\tau^{-\xi}.
		\end{equation}
		For the polynomial source, let $h=\tau-r$. By \eqref{eq:OU-polynomial-moment} and \eqref{eq:OU-window-comparability},
		\begin{align}
			&C_F\int_s^\tau e^{\eps_0(\tau-r)}
			S(\tau-r)\frac{1+|z|^2}{r^2}\,dr \leq
			C C_F\tau^{-2}
			\int_0^H e^{\eps_0h}(e^h+R^2)\,dh \leq
			C C_Fe^{\eps_0H}
			(1+R)^2\tau^{-2}(1+H).
			\label{eq:OU-polynomial-source-bound}
		\end{align}
		
		On the support of $\check\chi(\cdot,r)$, \eqref{eq:OU-cutoff-buffer} and \eqref{eq:OU-window-comparability} imply $|z|\geq2R(r)\geq\frac32R.$ Thus Lemma \ref{lem:OU-source-estimates} yields
		\begin{align}
			&C_F\int_s^\tau e^{\eps_0(\tau-r)}
			S(\tau-r)(|z|\check\chi(z,r))\,dr \notag\\
			&\qquad\leq
			C C_Fe^{-cR^2}\int_0^He^{(1+\eps_0)h}\,dh
			\leq C C_Fe^{\eps_0H}(1+R)^2e^{-cR^2}.
			\label{eq:OU-cutoff-source-bound}
		\end{align}
		
		By the choice $\eps_0\leq c_*(\kappa)\eta$, we have $e^{\eps_0H}\leq C_\eta\tau^\eta$. Use $H\leq C(1+\log\tau)$ and $(1+R)^2\leq C(1+\tau^{2\kappa})$ in \eqref{eq:OU-initial-term}--\eqref{eq:OU-cutoff-source-bound}. The factor $\tau^\eta(1+R)^2$ multiplying the last Gaussian is absorbed by decreasing $c$. By the variation of constants expression of $Z$, \eqref{eq:OU-Duhamel-comparison}, this proves
		\eqref{eq:OU-loss-estimate}.
		
	\end{proof}
	
	\begin{proof}[Proof of Proposition \ref{PropBoundary-global}]
		Apply Lemma \ref{LmC0byL2} to $\psi=v_{s_\tau}$. Since $v(\tau)=S(H(\tau))v_{s_\tau}$, this gives \eqref{eq:linear-OU-correct}.
	\end{proof}
	
	\begin{remark}
		The increase of $h$ is necessary for the estimate. Indeed, for data constant on the spherical factor, evaluating at $y$ as a functional on $L^2$ is $C_he^h\exp\left(\frac{|y|^2}{4(e^h+1)}\right),$ by the same calculation as in \eqref{eq:Mehler-L2-pointwise}. For fixed $h$, this is not polynomially bounded on $B_R$ as $R\to\infty$. Therefore, to get a polynomial bound, we require $R^2/(e^h+1)=O(\log R)$, and hence at least $h=2\log R-O(\log\log R)$.  The choice \eqref{eq:OU-sliding-window} is therefore sufficient and has the correct leading order.
	\end{remark}
	
	\subsection{Perturbed Ornstein--Uhlenbeck regularization}
	\label{SS:perturbedOU}
	
	In this section, we state a version of the OU regularization estimate for the perturbed OU operator that will be used later. Let
	\begin{equation}\label{eq:variable-L}
		L_rv=a^{AB}(\cdot,r)\nabla_{AB}v
		-\frac12y\cdot\nabla_yv+b^A(\cdot,r)\nabla_Av+v.
	\end{equation}
	The indices $A,B$ range over both factors of the cylinder.  We assume that $a$ is symmetric and
	\begin{equation}\label{eq:coefficient-smallness-OU}
		\|a-g^{-1}\|_{C^1}+\|b\|_{C^1}\leq\delta,
		\qquad 0<\delta\leq\bar\delta.
	\end{equation}
	All norms here are spatial and uniform in $r$. The following proposition gives similar estimates as the nonperturbed case. 
	
	\begin{proposition}[Perturbed OU estimate]\label{prop:perturbed-OU}
		Fix $0<\kappa\leq1/2$, $0<\xi<2$, $K>0$, and $\eta>0$. There are $\eps_*,c_*,T_0>0$ such that the following holds. Suppose \eqref{eq:coefficient-smallness-OU} holds with $\delta\leq c_*\eta$, and $Z\geq0$ is a Gaussian-energy distributional subsolution on $[s_\tau,\tau]$ of
		\begin{equation}\label{eq:perturbed-OU-inequality}
			\partial_rZ-L_rZ
			\leq\eps_0Z+C_F\left(
			\frac{1+|y|^2}{r^2}
			+|y|\check\chi(y,r)\right),
		\end{equation}
		where $\eps_0\leq\min\{\eps_*,c_*\eta\}$,
		\[
		\|Z(\cdot,r)\|_{L^2}\leq C_0r^{-\xi},
		\qquad
		\check\chi(y,r)=0\quad\text{for }|y|\leq2Kr^\kappa.
		\]
		Then
		\begin{equation}\label{eq:perturbed-OU-conclusion}
			\begin{split}
				\|Z(\cdot,\tau)\|_{C^0(B_{K\tau^\kappa})}
				\leq C_\eta C_0\tau^{2\kappa-\xi+\eta}.
			\end{split}
		\end{equation}
	\end{proposition}
	The proof is parallel to the above unperturbed case with the above explicit Mehler kernel replaced by a subsolution estimate. We postpone it to Appendix \ref{App:Perturbed Ornstein--Uhlenbeck regularization}. 
	
	\section{Extending the graphical radius to $O(\tau^{\kappa})$}\label{S:Extending_t_kappa}
	In this section, we first extend the graphical radius to a polynomial rate $O(\tau^{\kappa})$ for some small $\kappa>0$ using Proposition \ref{PropVelazquez}. We shall use the $L^2$-decay rate in Proposition \ref{PropR} as the input, which gives $\|u\|_{L^2}=O(\tau^{-\xi})=O(\tau^{-\frac{1}{3}})$ and we expect the graphical radius to be $O(\tau^{\kappa})$ with $\kappa=1/6-\eps$ as the output of Proposition \ref{PropR}. Thus, in this section, we devote most of the effort to showing that $Z=|u|+|Du|$ satisfies the differential inequality \eqref{EqVelazInequality}. Note that the nonlinearity in our rescaled mean curvature flow equation depends on second-order derivatives of the graphical function $u$, which is the main difficulty for us to apply the machinery of Vel\'azquez. In this section, we show how to apply the pseudolocality to give a bound on higher-order derivatives of the graphical function over expanding domains.

	\subsection{Extending the graphical radius}
	In this section, we extend the graphical radius from $O(\sqrt{\log \tau})$ to $O(\tau^\kappa)$, where $0<2\kappa<\xi=1/3$. The proof uses Proposition \ref{PropVelazquez} with a fixed spatial buffer: the cutoff is placed at radius $2K\tau^\kappa$, while the estimate is taken on $B_{K\tau^\kappa}$.
	\begin{proposition}\label{thm:mainthmfor sqrtt}
		Suppose $M_\tau$ is a rescaled mean curvature flow converging to $\cC_{n,k}$ in the $C^\infty_{loc}$-sense. Then there exists $\eps_2>0$ such that for any $\eps\in(0,\eps_2)$, there exist $K>0$, and $T=T(\eps,\lambda(M_0))>0$, such that when $\tau>T$, $M_\tau$ can be written as the graph of a function $u(\cdot,\tau)$ over $\cC_{n,k}\cap B_{K \tau^\kappa}$, with $\|u(\cdot,\tau)\|_{C^2(B_{ K \tau^\kappa})}\leq \eps$.
	\end{proposition}
	\begin{lemma}\label{LmVelaz}
		Let $Z=|\chi_{2K\tau^\kappa}u|+|\nabla (\chi_{2K\tau^\kappa}u)|$. There exists $\eps_1>0$ such that, if $|D^ju|\leq\eps_1$ for $j=0,1,2,3$ on $B_{2K\tau^\kappa+1}$, then, for all sufficiently large $\tau$, $Z$ satisfies \eqref{EqVelazInequality} for a uniformly parabolic operator $L_\tau$ in \eqref{eq:variable-L}. Moreover, the cutoff error is supported in $\{|y|\geq2K\tau^\kappa\}$.
	\end{lemma}   
	
	\begin{proof}[Proof of Lemma \ref{LmVelaz}]
		The proof is based on an observation on the nonlinearity $\mathcal Q$ in Appendix \ref{appendix:Nonlinear estimate}. Write the quasilinear graph equation in nondivergence form as
		\[
		\partial_\tau u=a^{AB}(J^1u)\nabla_{AB}u
		-\frac12y\cdot\nabla_yu+b^A(J^2u)\nabla_Au+u.
		\]
		The estimates in Appendix \ref{appendix:Nonlinear estimate} and the assumed $C^3$ smallness imply uniform ellipticity, $\|a-g^{-1}\|_{C^1}+\|b\|_{C^1}\leq C\eps_1$. Thus, after decreasing $\eps_1$, the coefficient assumptions of Proposition \ref{PropVelazquez} hold.  Differentiating the equation once gives a uniformly parabolic system for $\nabla u$ with the same principal coefficients and lower-order coefficients bounded by $C\eps_1$.
		
		Next, let $\chi=\chi_{2K\tau^\kappa}$.  The commutators in the equations for $\chi u$ and $\nabla(\chi u)$ are supported in $B_{2K\tau^\kappa+1}\setminus B_{2K\tau^\kappa}$ and are bounded by $C|y|\check\chi$, where $\check\chi=0$ on $B_{2K\tau^\kappa}$. Applying the scalar and vector-valued Kato inequalities to $\chi u$ and $\nabla(\chi u)$, and then adding the two inequalities, gives \eqref{EqVelazInequality}.
	\end{proof}
	
	Now we come back to the proof of Proposition \ref{thm:mainthmfor sqrtt}. The proof uses the ``extension-improvement'' that Colding-Minicozzi \cite{CM2} used to prove the uniqueness of cylinder tangent flows.

	\begin{proof}[Proof of Proposition \ref{thm:mainthmfor sqrtt}]
		Set $\xi=1/3$. Fix $\kappa\in(0,1/2)$ with $2\kappa<\xi$, and choose $0<\eta<\xi-2\kappa$.
		
		Choose the constants in Lemma \ref{LmVelaz} sufficiently small that the perturbed OU estimate is applicable with the above value of $\eta$. By smooth local convergence, there are $\overline R\gg1$ and $T_1\gg1$ such that the flow is graphical on $\cC_{n,k}\cap B_{\overline R}$ for every $\tau\geq T_1$, with the smallness required by pseudolocality.
		
		We next choose $K>0$ sufficiently small that $2K\tau^\kappa+2<\overline R$ throughout an initial interval whose length is greater than $2H(T_1)+10$. This initializes the following extending argument.
		
		Let $T^*$ be the supremum of the times for which $\|u(\cdot,r)\|_{C^1(B_{K\tau^\kappa})}\leq\eta_0$ holds for all $\tau\in[T_1,T^*]$. Repeated application of Theorem \ref{thm:pseudo-locality} on fixed length $2$ intervals shows that, throughout the same time range and away from the initial endpoint, the flow is graphical on $B_{2Kr^\kappa+1}$ and
		\[
		\sum_{j=0}^3
		\|\nabla^ju(\cdot,r)\|_{L^\infty(
			B_{2Kr^\kappa+1})}
		\leq\eps_1.
		\]
		
		We claim that, for every sufficiently large $\tau\leq T^*$,
		\[
		\|u(\cdot,r)\|_{H^1(B_{2Kr^\kappa+1})}
		\leq
		Cr^{-\xi}.
		\]
		Indeed, Proposition \ref{PropR} gives the desired estimate on $B_{r_{\mathrm{CM}}(r)}$. On the complementary region, pseudolocality gives a uniform $C^1$ bound, and therefore
		\[
		\begin{aligned}
			&\int_{B_{2K\tau^\kappa+1}
				\setminus B_{(\alpha \log(\tau))^{1/2}}}
			\left(|u|^2+|\nabla u|^2\right)\,d\mu
			\leq
			C\eps_1^2
			\tau^{-\alpha/(4+\eps)}
			\leq
			C\tau^{-2\xi},
		\end{aligned}
		\]
		where the first inequality using the Gaussian decay of $d\mu$, and the last inequality follows from the choice $\alpha>8/3$ in Proposition \ref{PropR}.
		
		Fix a target time $\tau\leq T^*$ and let $s_\tau=\tau-H(\tau)$. For $\tau$ sufficiently large, $[s_\tau,\tau]\subset[T_1,T^*]$. The preceding estimates show that the cutoff quantity
		\[
		Z
		=
		|\chi_{2Kr^\kappa}u|
		+
		|\nabla(\chi_{2Kr^\kappa}u)|
		\]
		satisfies all hypotheses of Proposition \ref{prop:perturbed-OU} on the full sliding window $[s_\tau,\tau]$. Consequently,
		\[
		\begin{aligned}
			\|u(\cdot,\tau)\|_{C^1(B_{K\tau^\kappa})}
			&\leq
			C \tau^{2\kappa-\xi+\eta}
			.
		\end{aligned}
		\]
		After increasing $T_1$, the right-hand side is less than $\eta_0/2$. Thus the bootstrap estimate improves strictly at the possible terminal time $T^*$. Cylindrical pseudolocality then extends the graphical region beyond $T^*$, contradicting the definition of $T^*$. Hence $T^*=\infty$.
	\end{proof}
	
	From the proof, we can get an even better $C^2$ bound for $u$.
	
	\begin{corollary}\label{cor_C2_bound_for_kappa<1/2}
		Suppose the assumptions in Proposition \ref{thm:mainthmfor sqrtt}. Further assume $\xi=1/3$, $\eta>0$ and $\kappa\in(0,1/2)$ such that $2\kappa-\xi+\eta<0$. Then there exist $C>0$ and a possibly smaller $K>0$ such that $\|u(\cdot,\tau)\|_{C^2(B_{ K \tau^\kappa})}\leq C\tau^{2\kappa-\xi+\eta}$.
	\end{corollary}
	
	\begin{proof}
		From the proof of Proposition \ref{thm:mainthmfor sqrtt}, we have $\|u(\cdot,\tau)\|_{C^1(B_{ K \tau^\kappa})}\leq C\tau^{2\kappa-\xi+\eta}$. Therefore, it suffices to improve the estimate to $C^2$. First, we have already proved that $\|u\|_{C^2(B_{ K \tau^\kappa})}<\eps_0$, and hence, by the mean value theorem, $u$ satisfies a uniformly parabolic PDE in $B_{4n}$, and then we can apply the interior parabolic Schauder estimate to show that $$\|u(\cdot,\tau)\|_{C^{2,\alpha}(B_{3n})}\leq \sup_{s\in[\tau-1,\tau]}\|u(\cdot,s)\|_{C^0(B_{4n})}\leq C(\tau-1)^{2\kappa-\xi+\eta}\leq \tilde C \tau^{2\kappa-\xi+\eta},$$ if $t$ is sufficiently large. Next, we observe that for any $\mathbf{y}\in\{0\}\times\R^k$, we have $M_\tau-e^{(\tau-\tau_0)/2}\mathbf{y}$ is also a RMCF. In particular, this shows that $\|u(\cdot,\tau)\|_{C^{2,\alpha}(B_{3n}(e^2\mathbf{y}))}\leq \tilde C \tau^{2\kappa-\xi+\eta}$, whenever $B_{4n}(\mathbf{y})\subset B_{K(\tau-1)^{\kappa}}$. This implies that $\|u(\cdot,\tau)\|_{C^2(B_{ K' \tau^\kappa})}\leq C'\tau^{2\kappa-\xi+\eta}$, for some constant $C'$ larger than $C$ in the proof of Proposition \ref{thm:mainthmfor sqrtt} and $K'$ slightly smaller than $K$ in the proof of Proposition \ref{thm:mainthmfor sqrtt}.
	\end{proof}

	\section{The invariant cones and $H^1$-normal form}\label{SDynamics}
	In Section \ref{S:Extending_t_kappa}, we have extended the graphical scale to $K\tau^\kappa$ for some small $K>0$ and $\kappa=1/6-\eps$. However, this graphical radius is still not sufficient for us to derive the geometric consequence, such as the isolatedness of nondegenerate singularities, which requires a graphical radius of order $\tau^{1/2}$.  Note that by Proposition \ref{PropVelazquez}, the estimate of the graphical radius is mainly constrained by the $L^2$-estimate of the graphical $u$ provided in Proposition \ref{PropR}. To further extend the graphical scale, we need to improve the $L^2$-estimate of $u$. In this section, we prove that the linear equation $\partial_\tau u=L_{\cC_{n,k}} u$ controls the dynamics of the rescaled mean curvature flow, and prove the $H^1$-normal form theorem. 
	
	Recall that $\|\cdot\|$ without a subscript denotes the Gaussian weighted $H^1$-norm of a function over $\cC_{n,k}$. We consider the region of radius $O(\tau^\kappa)$ for some $\kappa\in(0,1/2)$, which is the improved graphical radius from Section \ref{S:Extending_t_kappa}. Suppose $u(\cdot,\tau)$ is the graph function of the rescaled mean curvature flow over $\cC_{n,k}\cap B_{C\tau^\kappa}$. We let $\bar \eps(\tau):=\tau^{2\kappa-\xi+\eta}$, where we fix a constant $\xi-\eta\in(2\kappa,1)$, and by Corollary \ref{cor_C2_bound_for_kappa<1/2}, $\|u\|_{C^2(\cC_{n,k}\cap B_{\br(\tau)})}\leq C\bar \eps(\tau)$.
	
	In the following, we define $2\vartheta:=\xi-2\kappa-\eta$, and we choose $\br(\tau)=\tau^{\vartheta/8}$, for the technical purpose in Section \ref{SS:BC-rotation}.
	
	\subsection{Brendle-Choi's argument on rotations}
	\label{SS:BC-rotation}
	The natural function space to study these equations is $H^1(\cC_{n,k})$, which admits a natural direct sum decomposition into eigenspaces of the $L$-operator. Lying in the kernel of $L$, those $\theta_\alpha y_j\in\mathfrak{so}(n+1)$ represent infinitesimal rotations. These terms do not influence the major decay behavior, see \cite{CM2}; on the other hand, they do show up in the nonlinear analysis. Following an idea of Brendle-Choi \cite{BC1, BC2}\footnote{While \cite{BC1, BC2} studied ancient RMCF, as the RMCF equations are the same, their arguments can be simply adapted to our setting by reversing the time range.}, we modulo these modes by considering the rotated rescaled mean curvature flow $\tilde M_\tau=S_\tau M_\tau,\ S_\tau\in \mathrm{SO}(n+1)$. Let us write $\tilde M_\tau$ over $\mathcal{C}_{n,k}\cap B_{\br(\tau)}$ as the graph of a function $\tilde u$, where the rotation $S_\tau$ is chosen to make sure that $\chi_{\br(\tau)} \tilde u$ does not have Fourier modes corresponding to rotations. The existence of $S_{\tau}$ can be proved as in \cite[Proposition 2.4]{BC2}.
	
	We give the details needed in the present forward-time setting.  Let $\mathcal R$ be the finite-dimensional space spanned by
	\[
	\phi_A(x)=\langle Ax,\mathbf n(x)\rangle,
	\qquad A\in\mathfrak{so}(n+1),
	\]
	modulo the stabilizer of $\cC_{n,k}$.  If $u_S$ is the graph function of $S\Graph(u)$, define
	\[
	\mathcal F(S,\tau)(A)=\int_{\cC_{n,k}}\chi_{\br(\tau)}u_S\phi_Ae^{-|x|^2/4}.
	\]
	At $(S,u)=(\id,0)$ the derivative in $S$ is the Gaussian Gram form of the $\phi_A$'s with cutoff.  It converges exponentially to the positive definite full Gram form and hence is uniformly invertible, on the complement of the stabilizer, for large $t$.  The parameter-dependent implicit function theorem produces a unique small $C^1$ family $S_t$ satisfying
	\begin{equation}\label{eq:rotation-orthogonality}
		\langle\chi_{\br(\tau)}\widetilde u,\phi_A\rangle_{L^2}=0,
		\qquad \phi_A\in\mathcal R.
	\end{equation}
	
	For a function $w$ with small $C^1$-norm, let $\Phi(A,w)$ be the graph function over $\cC_{n,k}$ of $e^{A}\Graph(w)$. Thus $\widetilde u=\Phi(A_\tau,u)$ and $u=\Phi(-A_\tau,\widetilde u)$. Let $\rho_A:=\Phi(A,0)$ be the graph function of the rotated cylinder $e^A\cC_{n,k}$. Because the rescaled mean curvature flow equation is invariant under $\mathrm{SO}(n+1)$, for a constant $A$, $\Phi(A,w(\cdot,\tau))$ solves $\partial_\tau\Phi=L\Phi+\cQ(J^2\Phi)$ whenever $w$ does. In particular, since $e^A\cC_{n,k}$ is a shrinker,
	\begin{equation}\label{eq:rotated-cylinder}
		L\rho_A+\cQ(J^2\rho_A)=0.
	\end{equation}
	Since $\phi_A=\langle Ax,\mathbf n\rangle$ is the normal component of $Ax$, we have $\rho_A=\phi_A+O(|A|^2)$. Differentiating \eqref{eq:rotated-cylinder} at $A=0$ (recall that $\cQ$ is at least quadratic) recovers $L\phi_A=0$.
	
	The nonlinear action of rotations is instead accounted for exactly by $\Phi$. Since $\Phi(0,w)=w$ and $\Phi(A,w)$ is smooth in $(A,w,\nabla w)$, Taylor expansion gives, on $B_{\br(\tau)-1}$ and for $|A|\br(\tau)^2\leq c$,
	\begin{equation}\label{eq:Phi-expansion}
		|E|+|\nabla E|\leq C|A|(1+|y|)^2\big(|w|+|\nabla w|+|\nabla^2w|\big)+C\big(|w|+|\nabla w|\big)\big(|w|+|\nabla w|+|\nabla^2w|\big),
	\end{equation}
	where $E:=\Phi(A,w)-\rho_A-w$. In particular, the equation of $\widetilde u=\Phi(A_\tau,u)$ is the same equation as that of $u$, with the only additional term $\mathcal A$ coming from $\partial_\tau A_\tau$ (this is \eqref{eq:bar equation}).

	Differentiability follows because the flow and the moving cutoff are differentiable in $t$.  The construction is pointwise in time and the subsequent equation is the same RMCF equation as in \cite{BC1,BC2}; therefore, no ancientness is used, and the estimates are integrated in the forward direction here.
	
	Let $\overline v=\chi_{\br(\tau)}\widetilde u$ and $\mathcal A=\langle S'_tS_t^{-1}x,\mathbf n(x)\rangle$.  hen
	\begin{equation}\label{eq:bar equation}
		\partial_t\overline v=L\overline v+\overline{\cQ}(J^2\overline v)+\chi_{\br(\tau)}\mathcal A,
	\end{equation}
	where $\overline{\cQ}$ contains the nonlinear terms and the commutators with the cutoff. We have the following lemma.
	\begin{lemma}\label{LmA}
		For sufficiently large $\tau$, 	we have $\|\chi_{\br(\tau)}\cA\|_{L^2}\leq 2\|\overline{\cQ}\|_{L^2}$.
	\end{lemma}
	
	\begin{proof}[Proof of Lemma \ref{LmA}] 
		From the definition of $\cA$, we have 
		$\int (\bar{\cQ}(J^2u)-\chi\cA(\tau))\cA(\tau)e^{-\frac{|x|^2}{4}}\,dx=0$
		and by Plancherel and Cauchy-Schwarz inequality,
		\begin{equation}\label{EqRotation}
			\|\chi_{\br(\tau)}\cA\|_{L^2}^2\leq \int \langle \chi_{\br(\tau)}\cA,\cA\rangle e^{-\frac{|x|^2}{4}} dx\leq \|\bar{\cQ}(J^2u)\|_{L^2}\|\cA\|_{L^2}. 	
		\end{equation}
		As $\cA$ is a concrete linear function in $x$, when $\tau$ is sufficiently large, $\|\cA\|_{L^2}\leq 2\|\chi_{\br(\tau)}\cA\|_{L^2}$. 
	\end{proof}

	As a consequence, we may also write the equation of $\bar v$ as 
	\begin{equation}\label{eq:bar equation after reduction}
		\pr_\tau \overline{v}=L\overline{v}+\overline{\cB}(J^2\overline{v}).
	\end{equation}
	In the rest of this section, we remove all the overline notations for notational simplicity.

	\subsection{The cone theorem}
	The invariant cone theorem is an important tool in dynamical systems to study the stable/unstable manifolds. Let $E$ be the orthogonal complement of $\text{span}\{\theta_\alpha y_j\}$ (recall we have modulo these rotations in Section \ref{SS:BC-rotation}) and we introduce the direct sum decomposition $E=E^+\oplus E^0\oplus E^-$ where $E^+$ (respectively $E^0$ and $E^-$) is spanned by eigenfunctions of $L$ with positive (respectively $0$ and negative) eigenvalues. We denote by $\Pi_+,\ \Pi_-,\ \Pi_0$ the $L^2$-projections to the spaces $E^+, E^-, E^0 $ respectively. 
	
	We introduce a double cone construction that is used to suppress the $E^+$-components and manifest the $E^0$-component. Let $\alpha>0$ be a positive number. We introduce two cones $\cK_{\geq 0}$ and $\cK_{0}$ as follows
	$$\cK_{\geq 0}(\alpha):=\left\{u=(u_+,u_0,u_-)\in E^+\oplus E^0\oplus E^-\ |\ \|u_++u_0\|\geq \alpha\|u_-\|\right\}$$
	is a $\alpha$-cone around $E^+\oplus E^0$
	and 
	$$\cK_{0}(\alpha):=\left\{u=(u_+,u_0,u_-)\in E^+\oplus E^0\oplus E^-\ |\ \|u_0\|\geq 
	\alpha\|u_++u_-\|\right\}$$
	is a $\alpha$-cone around $ E^0$, and 
	$$\cK_{+}(\alpha):=\left\{u=(u_+,u_0,u_-)\in E^+\oplus E^0\oplus E^-\ |\ \|u_+\|\geq 
	\alpha\|u_0+u_-\|\right\}$$
	is a $\alpha$-cone around $ E^+$. All are narrower when $\alpha$ is large. 
	
	We shall also need to compare the difference between the perturbed rescaled mean curvature flow and the unperturbed one under evolution. For this purpose, we introduce the following setting. Let $u_1,u_2$ be two graphical rescaled mean curvature flows over $\cC_{n,k}$, and $\|u_i(\cdot,\tau)\|_{C^{2,\alpha}}\leq \eps_0$ inside $\cC_{n,k}\cap B_{\br(\tau)}$. Then we write $v=\chi(u_1-u_2)$, where $\chi$ is a smooth cutoff function that is 0 outside the ball $B_{\br(\tau)}$ and is $1$ inside the ball $B_{\br(\tau)-1}$, and we can calculate
	\begin{equation}\label{EqDifference}
		\pr_\tau v=L v+\delta\cB. 
	\end{equation}
	where $\delta\cB=\delta\cB(J^2u_1,J^2u_2)=\cB(J^2u_1)-\cB(J^2u_2)$. In $\delta\cB$, we replace $u_i$ by $\chi u_i$ creating an error supported on the annulus $\mathbb A_{\br(\tau)}$ denoted by $\mathcal E(J^2u_1,J^2u_2). $ Thus we can write 
	\begin{equation}\label{EqPE}
		\delta\cB =P v+\mathcal E(J^2u_1,J^2u_2), \quad P = \int_0^1 D \cB(sJ^2(\chi u_1)+(1-s)J^2(\chi u_2))ds.		
	\end{equation}
	
	In this paper, we only consider the case that $u_1=u$ is a given rescaled mean curvature flow and $u_2=0$ is the rescaled mean curvature flow of the shrinking cylinder. While the idea of studying the difference of the graphs of two rescaled mean curvature flows is the same as this special case, the calculations are more involved, and it is only used in our forthcoming work. Therefore, we only discuss this special case in this paper.
	
	Throughout the rest of this section, we consider discrete-time $m\in\mathbb Z_{+}$ for presentation simplicity, although the proof for continuous time is almost verbatim. We also remind the readers that $\mathbb A_{\mathbf r(\tau)}$ is the annulus region $\cC_{n,k}\cap(B_{\br(\tau)}\backslash B_{\br(\tau)-1})$.
	
	In the following, we use $\omega(\tau)=e^{-\frac{\mathbf{r}(\tau)^2}{5}}$ to indicate the scale of the error from cut-off. The exponent $1/5$ is chosen below the Gaussian exponent $1/4$ so
	that every fixed polynomial factor on the cutoff annulus is absorbed by $\omega$.

	The following theorem shows that the linearized equation approximates the nonlinear equation nicely.

	\begin{proposition}\label{PropApproximation}
		Let $v(\tau)=\chi(\tau)(u_1( \tau)-u_2(\tau)):\ \cC_{n,k}\cap B_{\br(\tau)}\to \R$ be as above a solution to \eqref{EqDifference}. Let  $\eps_m
		:=
		\sup_{\tau\in[m,m+1]}
		\|v(\tau)\|_{C^2(B_{\mathbf r(\tau)})}$.
		For every sufficiently large integer $m$,
		\begin{equation}\label{eq:one-step-additive}
			\|v(m+1)-e^Lv(m)\|\leq C\eps_m\left(
			\|v(m)\|+\omega(m)^{1/2}
			\right).
		\end{equation}
		Here $e^Lv(m)$ means the heat semigroup generated by $L$ acting on $v(m)$ for time $1$.
	\end{proposition}
	
	In particular, if $\| v(m)\|\geq \omega(m)^{1/2}$, we have $\|v(m+1)-e^Lv(m)\|\leq 2C\eps_m \|v(m)\|$. This assumption either holds for all sufficiently large $m$, or fails for all sufficiently large $m$. This fact is verified in the following lemma, which is proved in Section \ref{SS:ConeBoundary}.

	\begin{proof}
		Let $w$ solve
		\[
		\partial_\tau w=Lw,
		\qquad
		w(m)=v(m).
		\]
		Then
		\[
		\partial_\tau(v-w)
		=
		L(v-w)+\delta\cB.
		\]
		
		Since we have $\dt\cB=Pv+\mathcal E$, we get $\|\dt \cB\|^2_{L^2}\leq 2\|Pv\|^2_{L^2}+2\|\mathcal E\|^2_{L^2}$. Notice that $\dt\cB$ is defined in equation \eqref{EqDifference} independent of $w$, and $\dt \mathcal B=Pv+\mathcal E$, where $\mathcal E$ is only supported on $\A_{\br(\tau)}$. The coefficients of $P$ are bounded by $C\eps$ because $\|u_i\|_{C^2}\leq\eps$. The cutoff commutator has the schematic form
		\[
		\mathcal E=(\nabla\chi)*\nabla(u_1-u_2)+(\nabla^2\chi+\partial_t\chi)*(u_1-u_2),
		\]
		with uniformly bounded coefficients on the unit-width annulus.  Since the Gaussian measure of that annulus is at most $C\omega(\tau)$ 
		\begin{equation}\label{eq:cB L_2}
			\|\dt\cB\|_{L^2}^2\leq \eps^2\|J^2 v\|_{L^2}^2+\omega(\tau) \|v\|^2_{C^2(\mathbb A_{\br(\tau)})},
		\end{equation}
		where we used Proposition \ref{PropR}.
		
		The term $\|J^2v\|_{L^2}^2$ is bounded by a multiple of $\|\mathcal L v\|_{L^2}^2+\|\nabla v\|_{L^2}^2$, and $\|\mathcal L v\|_{L^2}^2$ is further bounded by $\|\nabla^2 v\|^2+\|\nabla v\|^2$ by integrating the Bochner formula (details can be found in \cite[Section 5]{CM3} and \cite[Proposition 4.6]{SX3})
		\[
		\frac{1}{2}\cL|\nabla v|^2=|\nabla^2 v|^2+\langle \nabla\cL v,\nabla v\rangle +(\Ric+\Hess_{|x|^2/4})(\nabla v,\nabla v).
		\]
		
		In summary, we have
		\begin{equation}\label{eq:delta-B-additive}
			\|\delta\cB(\tau)\|_{L^2}^2
			\leq
			C\eps_m^2\|v(\tau)\|_{H^2}^2
			+
			C\eps_m^2\omega(m),
			\qquad
			\tau\in[m,m+1].
		\end{equation}
		Indeed, the cutoff terms are supported in
		$\mathbb A_{\mathbf r(\tau)}$, and
		\[
		\int_{\mathbb A_{\mathbf r(\tau)}}
		(1+|y|)^N e^{-|x|^2/4}
		\leq
		C_Ne^{-\mathbf r(\tau)^2/5}
		\leq
		C_N\omega(m).
		\]
		
		The $H^1$ energy estimate gives
		\begin{equation}\label{eq:H1-difference-additive}
			\frac d{d\tau}
			\|v-w\|^2
			\leq
			C\|v-w\|^2
			+
			C\eps_m^2\|v\|_{H^2}^2
			+
			C\eps_m^2\omega(m).
		\end{equation}
		
		We also claim that
		\begin{equation}\label{eq:H2-time-integral}
			\int_m^{m+1}
			\|v(\tau)\|_{H^2}^2\,d\tau
			\leq
			C\left(
			\|v(m)\|^2+\eps_m^2\omega(m)
			\right).
		\end{equation}
		In fact, after integrating by parts,
		\[
		\frac d{d\tau}
		\|\nabla v\|_{L^2}^2
		\leq
		-c\|\cL v\|_{L^2}^2
		+
		C\|v\|^2
		+
		C\eps_m^2\omega(m).
		\]
		The drift Bochner formula shows that $\|v\|_{H^2}^2
		\leq
		C\left(
		\|\cL v\|_{L^2}^2+\|v\|^2
		\right),$ then gives \eqref{eq:H2-time-integral}. Applying Gronwall's inequality to \eqref{eq:H1-difference-additive}, and using \eqref{eq:H2-time-integral}, yields
		\[
		\|v(m+1)-e^Lv(m)\|
		\leq
		C\eps_m
		\left(
		\|v(m)\|+\omega(m)^{1/2}
		\right).
		\]
	\end{proof}

	\begin{lemma}[Cone theorem]
		\label{lem:local-cone-estimates}
		There are constants $R_*>1$, $\alpha_*>0,$ $\beta_*>0,$ and $\nu>0$, such that the following holds for all sufficiently large $m$. Assume
		\begin{equation}\label{eq:active-boundary-scale}
			\|v(m)\|
			\geq
			R_*\omega(m)^{1/2}.
		\end{equation}
		
		\begin{enumerate}
			\item If $v(m)\in\cK_{\geq0}(\alpha)$, $C\eps_m\leq\alpha\leq C\eps_m^{-1}$, then $v(m+1)\in
			\cK_{\geq0}((1+\nu)\alpha)$.
			
			\item If $v(m)\in\cK_+(\beta)$, $C\eps_m\leq\beta\leq C\eps_m^{-1}$, then $v(m+1)\in \cK_+((1+\nu)\beta)$, and
			\[
			\|\Pi_+v(m+1)\|
			\geq
			(1+\nu)
			\|\Pi_+v(m)\|.
			\]
			
			\item If $v(m)\notin\cK_{\geq0}(\alpha_*)$, then
			\begin{equation}\label{eq:normalized-stable-contraction}
				\frac{\|v(m+1)\|}{\omega(m+1)^{1/2}}
				\leq
				q
				\frac{\|v(m)\|}{\omega(m)^{1/2}}
				+
				C\eps_m
			\end{equation}
			for some $q\in(0,1)$ independent of $m$.
		\end{enumerate}
	\end{lemma}
	
	\begin{proof}
		We use $\mu_-$ to denote the largest negative eigenvalue of $L$ and $\mu_+$ to denote the smallest positive eigenvalue of $L$. By Proposition \ref{PropApproximation},
		\[
		\|v(m+1)-e^L v(m)\|
		\leq
		C\eps_m
		\left(
		1+R_*^{-1}
		\right)
		\|v(m)\|.
		\]
		
		For item (1), we have
		\[
		\|\Pi_{\geq 0}v(m+1)\|
		\geq
		\|\Pi_{\geq 0}v(m)\|-C\eps_m(\|\Pi_{\geq 0}v(m)\|+\|\Pi_{<0}v(m)\|).
		\]
		Because
		\[
		\|\Pi_{< 0}v(m+1)\|
		\leq
		e^{-\mu_-}\|\Pi_{< 0}v(m)\|
		+
		C\eps_m(\|\Pi_{\geq 0}v(m)\|+\|\Pi_{< 0}v(m+1)\|).
		\]
		If $\|\Pi_{\geq 0}v(m)\|\geq\alpha \|\Pi_{<0}v(m)\|$, then
		\[
		\frac{\|\Pi_{\geq 0}v(m+1)\|}{\|\Pi_{< 0}v(m+1)\|}
		\geq
		\frac{
			1-C\eps_m(1+\alpha^{-1})
		}{
			e^{-\mu_-}\alpha^{-1}
			+
			C\eps_m(1+\alpha^{-1})
		}.
		\]
		For $C\eps_m\leq\alpha\leq c\eps_m^{-1}$, the last quantity is at least $(1+\nu)\alpha$. 
		
		The proof of item (2) is analogous. The positive component is expanded by at least $e^{\mu_+}$, while the operator norm of $e^L$ on $E^0\oplus E^-$ is at most $1$.
		
		For item (3), choose $\alpha_*$ sufficiently small. Outside
		$\cK_{\geq0}(\alpha_*)$,
		\[
		\|\Pi_{\geq0}v(m)\|
		<
		\alpha_*\|\Pi_-v(m)\|.
		\]
		Therefore
		\[
		\|v(m+1)\|
		\leq
		\left(
		e^{-\mu_-}
		+
		C\alpha_*
		+
		C\eps_m
		\right)
		\|v(m)\|
		+
		C\eps_m\omega(m)^{1/2}.
		\]
		Moreover,
		\[
		\frac{\omega(m+1)^{1/2}}{\omega(m)^{1/2}}
		=
		\exp\left(
		-\frac{\mathbf r(m+1)^2-\mathbf r(m)^2}{10}
		\right).
		\]
		For $\kappa<1/2$, this quotient converges to $1$, while for $\kappa=1/2$ it converges to $e^{-K_0^2/10}$. We choose $K_0$ so that $\frac{K_0^2}{10}<\mu_-.$ After decreasing $\alpha_*$ and increasing the initial time, division by $\omega(m+1)^{1/2}$ proves \eqref{eq:normalized-stable-contraction}.
	\end{proof}
	
	The next Lemma shows that if the cut-off error is very small relative to $\|v\|$, then the neutral mode will dominate after quite a long time. This estimate is coarse.
	
	\begin{lemma}
		\label{lem:coarse-neutral-reduction}
		Suppose $\|v(m)\|/\omega(m)^{1/2}$
		is sufficiently large. Then, after at most
		$C\log\mathbf r(m)$ times, and as long as $\|v\|\geq R_* \omega^{1/2}$,
		\begin{equation}\label{eq:coarse-neutral}
			\|\Pi_{\neq0}v(j)\|
			\leq
			C\delta_* \br(j)^{-2}\|\Pi_0v(j)\|.
		\end{equation}
		Once \eqref{eq:coarse-neutral} holds, it continues to hold.
	\end{lemma}
	
	\begin{proof}
		By Lemma \ref{lem:local-cone-estimates}, item 3, when $\|v(m)\|/\omega(m)^{1/2}$ becomes very large, $v(m)\in\cK_{\geq0}(\alpha_*)$. Then item (1) improves the cone until
		\[
		\|\Pi_{<0}v(j)\|
		\leq
		C\delta_* \br(j)^{-2}\|\Pi_{\geq0}v(j)\|.
		\]
		The number of required steps is $O(\log\mathbf r(m))$.
		
		We next control the positive component. Suppose that at a sufficiently large active time,
		\[
		\|\Pi_{>0}v(j)\|
		\geq
		\delta_*\br(j)^{-2}\|\Pi_{\leq0}v(j)\|.
		\]
		Since $\delta_*\br(j)^{-2}\approx \tau^{-\vartheta/4}\gg\tau^{-2\vartheta}\approx\eps_j$, item (2) of Lemma \ref{lem:local-cone-estimates} applies. The corresponding positive cone is preserved and becomes narrower, and $\|\Pi_{>0}v(j)\|$ grows exponentially, this is a contradiction to $\|v(j)\|\to0.$ Hence
		\[
		\|\Pi_{>0}v(j)\|
		<
		\delta_*\br(j)^{-2}\|\Pi_{\leq0}v(j)\|
		\]
		at every sufficiently large time. This shows \eqref{eq:coarse-neutral}.
	\end{proof}

	\subsection{Control the boundary term}\label{SS:ConeBoundary}
	
	Lemma \ref{lem:local-cone-estimates} relies on \eqref{eq:active-boundary-scale}, which is the smallness of boundary cut-off error relative to $v$. The main goal of this section is to show the following dichotomy: either $\|v\|$ is small relative to $\omega^{1/2}$ for all future time, or $\|v\|/\omega^{1/2}$ tends to infinity.
	
	\begin{lemma}
		\label{LmContinuation}
		Assume $K_0^2<10\mu_-.$ Then exactly one of the following alternatives occurs:
		\begin{enumerate}
			\item there is a constant $C$ such that
			\begin{equation}\label{eq:small-boundary-branch}
				\|v(\tau)\|
				\leq
				C\omega(\tau)^{1/2}
			\end{equation}
			for all sufficiently large $\tau$;
			\item
			\begin{equation}\label{eq:large-boundary-branch}
				\frac{\|v(\tau)\|^2}{\omega(\tau)}
				\longrightarrow\infty.
			\end{equation}
		\end{enumerate}
		Moreover, if $v(\tau)$ remains outside $\cK_{\geq0}(\alpha_*)$ for all sufficiently large $\tau$, then \eqref{eq:small-boundary-branch} holds; in the alternative \eqref{eq:large-boundary-branch},
		\begin{equation}\label{eq:sharp-nonneutral-control}
			\|\Pi_{\neq0}v(\tau)\|
			\leq
			C\bar\eps(\tau)
			\|\Pi_0v(\tau)\|
		\end{equation}
		for all sufficiently large $\tau$.
	\end{lemma}
	
	To prove this lemma, we need to understand how the neutral mode evolves when it dominates. The idea is that after modulo $\text{span}\{\theta_\alpha y_j\}$, the remaining neutral modes are $h_2(y_i)$'s and $h_1(y_i)h_1(y_j)$'s. Recall that we have defined the normalized (in $\|\cdot\|_{L^2(\cC_{n,k})}$) eigenfunctions $$H_2(y_i)=c_0^{k-1}\cG_{n,k}^{-1/2}h_2(y_i),\quad H_{1,1}(y_i,y_j)=c_0^{k-2}\cG_{n,k}^{-1/2}h_1(y_i)h_1(y_j),$$ where $\cG_{n,k}$ is given in \eqref{EqcG}. Define 
	$
	\begin{cases}
		m_{ii}(\tau)&:=\langle v,H_2(y_i)\rangle
		=
		c_0^{k-1}\cG_{n,k}^{-1/2}\langle v,h_2(y_i)\rangle,\\
		m_{ij}(\tau)&:=\langle v,H_{1,1}(y_i,y_j)\rangle
		= c_0^{k-2}\cG_{n,k}^{-1/2}\langle v, h_1(y_i)h_1(y_j)\rangle.        
	\end{cases}$
	Then we get $$
	\begin{aligned}
		v=&\Pi_{0}v+\Pi_{\neq 0}v=\sum_i m_{ii}H_2(y_i)+\sum_{i< j}m_{ij}H_{1,1}(y_i,y_j)+\Pi_{\neq 0}v.
	\end{aligned}
	$$
	We also define $\bar m_{ii}(\tau)=a\langle v(\tau),h_2(y_i)\rangle$ with $a=\sqrt{2}c_0$ and $\bar m_{ij}(\tau)=\langle v(\tau), h_1(y_i)h_1(y_j)\rangle$. The reason that we let $a=\sqrt{2}c_0$ in $m_{ii}$ is to make sure that later, when we compute the derivatives of $\bar{m}_{ii}$ and $\bar m_{ij}$, we will get uniform expressions. The choice of $a$ will become evident in later proofs. Let $\mathrm{M}(\tau)$ be a symmetric $k\times k$ matrix, whose entries are given by $\bar{m}_{ii}$ and $\bar{m}_{ij}$. 
	
	The following lemma provides a detailed evolution equation estimate for $\mathrm{M}$.
	
	\begin{lemma}\label{LmMatrix_0}
		Suppose at some time $\tau$, we have $v(\tau)\in \cK_0(\alpha(\tau))$, and $\|v(\cdot,\tau)\|_{C^2(B_{\br(\tau)})}\leq \bar\eps(\tau)$. Then there exists $C$ only depends on $n,k$ and $\alpha_0$ such that $\mathrm{M}(\tau)$ satisfies the following ODI:
		\begin{equation}
			\label{eq:ODE-0}
			|\mathrm M'
			+\gamma\mathrm M^2|
			\leq 
			C\left(
			1+\mathbf r(\tau)^2
			\right)
			\left(
			\alpha(\tau)^{-1}+\bar\eps(\tau)
			\right)
			|\mathrm M|^2
			+
			C\bar\eps(\tau)\omega(\tau).
		\end{equation}
		where $\gamma=\frac{1}{\varrho}c_0^{2(k-1)}\cG_{n,k}^{-1}$ is a fixed constant.
	\end{lemma}
	\begin{proof} In the following, $C$ may vary line to line, but only depend on $n,k$ and $\alpha_0$. By Section \ref{SS:BC-rotation}, we may write the equation of $v$ as $\pr_\tau v=Lv+\cQ(J^2 v)+\chi_{\br}\cA(\tau)$. From Proposition \ref{PropQ}, the remainder $\cQ(J^2u)$ of the rescaled mean curvature flow equation is given by
		\[
		\cQ(J^2u)=-\frac{1}{2\varrho}(u^2+4u\Delta_\theta u+2|\nabla_\theta u|^2)+\mathcal C(J^2u),
		\]
		where $\mathcal C(J^2u)=O(\|u\|_{C^{2}}(|u|^2+|\nabla u|^2))$. Then we take the time derivative of $\bar m_{ii}$ and $\bar m_{ij}$ to get
		\begin{equation}\label{eq:Aderivative-0}
			\begin{split}
				\bar m_{ii}'(\tau)
				&=
				a\int (Lv+\chi\cQ(v))h_2(y_i)d\mu
				\\
				&=
				a\int \left(-\frac{1}{2\varrho}(v^2+4v\Delta_\theta v+2|\nabla_\theta v|^2)+\mathcal C(J^2v)+\chi_{\br(\tau)}\cA(\tau)\right)h_2(y_i)d\mu
				,
				\\
				\bar m_{ij}'(\tau)
				&=
				\int \left(-\frac{1}{2\varrho}(v^2+4v\Delta_\theta v+2|\nabla_\theta v|^2)+\mathcal C(J^2v)+\chi_{\br(\tau)}\cA(\tau)
				\right)h_1(y_i)h_1(y_j)d\mu
				.
			\end{split}
		\end{equation}
		After integration by parts on the $\Delta_\theta$ term, we have 
		\[
		\begin{split}
			\bar m_{ii}'(\tau)
			&=
			a\int \left(-\frac{1}{2\varrho}(v^2-2|\nabla_\theta v|^2)+\mathcal C(J^2v)+\chi_{\br(\tau)}\cA(\tau)\right)h_2(y_i)d\mu
			,
			\\
			\bar m_{ij}'(\tau)
			&=
			\int \left(-\frac{1}{2\varrho}(v^2-2|\nabla_\theta v|^2)+\mathcal C(J^2v)+\chi_{\br(\tau)}\cA(\tau)
			\right)h_1(y_i)h_1(y_j)d\mu
			.
		\end{split}
		\] 
		We can decompose the integral involving $v^2$ as
		
		\[
		\begin{split}
			\int v^2 h_2(y_i)e^{-|x|^2/4}
			&=\int [(\Pi_{=0}v)^2+2v(\Pi_{\neq 0}v)-(\Pi_{\neq 0}v)^2 ]h_2(y_i)e^{-|x|^2/4},
			\\
			\int v^2 h_1(y_i)h_1(y_j)e^{-|x|^2/4}
			&=\int [(\Pi_{=0}v)^2+2v(\Pi_{\neq 0}v)-(\Pi_{\neq 0}v)^2 ]h_1(y_i)h_1(y_j)e^{-|x|^2/4}
			.
		\end{split}
		\]
		The integral involving $(\Pi_{=0}v)^2$ can be calculated explicitly by Lemma \ref{LmTriple}; Let $\phi$ be one of the normalized quadratic neutral eigenfunctions. On $B_{\mathbf r(\tau)+1}$, $|\phi|\leq
		C\left(
		1+\mathbf r(\tau)^2
		\right)$. Hence
		\[
		\left|
		\int_{B_{\mathbf r+1}}
		\left(
		v^2-(\Pi_{=0}v)^2
		\right)
		\phi\,d\gamma
		\right|
		\leq
		C(1+\mathbf r^2)
		\left(
		|\mathrm M|\|(\Pi_{\neq 0}v)\|_{L^2}
		+
		\|(\Pi_{\neq 0}v)\|_{L^2}^2
		\right),
		\]
		and
		\[
		\left|
		\int
		|\nabla_\theta (\Pi_{\neq 0}v)|^2\phi\,d\gamma
		\right|
		\leq
		C(1+\mathbf r^2)\|(\Pi_{\neq 0}v)\|^2.
		\]
		By Proposition \ref{PropQ},
		\[
		|\cC(J^2v)|
		\leq
		C\bar\eps
		\left(
		|v|^2+|\nabla v|^2
		\right).
		\]
		So the integral involving it is bounded by
		$
		C\bar\eps(1+\mathbf r^2)|\mathrm M|^2.
		$
		
		Finally, the term involving $\chi_{\br(\tau)}\cA(\tau)$ is bounded by $O(\bar\eps(\tau)\omega(\tau))$. In fact, $\cA(\tau)$ is the linear combination of infinitesimal rotations, so it is orthogonal to $h_2(y_i)$ and $h_1(y_i)h_1(y_j)$. Thus, the integral involving $\chi_{\br(\tau)}\cA(\tau)$ can be estimated by $\|\cA(\tau)\|_{L^2}\omega(\tau)$, and by Lemma \ref{LmA}, it is bounded by $O(\|v\|_{C^2(B_{\br(\tau)})}\omega(\tau))=O(\bar\eps(\tau)\omega(\tau))$. 
		
		This gives the desired inequality.
	\end{proof}

	\begin{proof}[Proof of Lemma \ref{LmContinuation}]
		We first work at integer times and suppose $\|v(m)\|=\beta(m)\omega(m)^{1/2}$. If $\sup_m \beta(m)<\infty$, then \eqref{eq:small-boundary-branch} holds at integer times. The energy estimate on every interval $[m,m+1]$ extends it to all times.
		
		Assume now that $\sup_m\beta(m)=\infty$. Choose a large number $A$, and let $m_0$ be the first sufficiently large integer such that $\beta(m_0)\geq A.$ Proposition \ref{PropApproximation} gives
		\[
		\beta(m_0)
		\leq
		C_0\beta(m_0-1)+C\eps_{m_0-1}.
		\]
		Thus, after choosing $A$ large,
		\[
		\beta(m_0-1)\geq C_0^{-1}A/2 \geq R_*.
		\]
		If $v(m_0-1)\notin\cK_{\geq0}(\alpha_*)$, then \eqref{eq:normalized-stable-contraction} yields
		\[
		\beta(m_0)
		\leq
		q\beta(m_0-1)+C\eps_{m_0-1}
		<
		A,
		\]
		a contradiction. Hence $v(m_0-1)\in\cK_{\geq0}(\alpha_*).$
		
		Lemma \ref{lem:coarse-neutral-reduction} gives an integer $m_1=m_0+O(\log\mathbf r(m_0))$ such that
		\begin{equation}\label{eq:coarse-neutral-at-m1}
			\|\Pi_{\neq0}v(m_1)\|
			\leq
			C\delta_*\br(\tau)^{-2}\|\Pi_0v(m_1)\|.
		\end{equation}
		
		Let $[m_1,T_*)$ be the maximal interval on which $\frac{\|v(\tau)\|}{\omega(\tau)^{1/2}}\geq R_*$. On this interval, \eqref{eq:coarse-neutral-at-m1}, Lemma \ref{lem:coarse-neutral-reduction}, and Lemma \ref{LmMatrix_0} imply
		\[
		|\mathrm M'
		+\gamma\mathrm M^2|\leq 
		\delta_0|\mathrm M|^2,
		\]
		where $\delta_0>0$ can be made arbitrarily small by first choosing $\delta_*$ small and then increasing $m_1$. Indeed, on the time range,
		\[
		(1+\mathbf r^2)\delta_*\br(\tau)^{-1}
		\leq
		C\delta_*,
		\qquad 
		(1+\mathbf r^2)\bar\eps\to 0,
		\qquad 
		\omega
		\leq
		C R_*^{-2}|\mathrm M|^2.
		\]
		
		Set
		$
		\Lambda(\tau)
		=
		|\mathrm M(\tau)|.
		$
		The ordered eigenvalues of a $C^1$ symmetric matrix are locally Lipschitz. At almost every time, if $\Lambda$ is realized by a positive eigenvalue $\lambda$, then 
		$
		\Lambda'
		=
		-\gamma\Lambda^2+O(\delta_0\Lambda^2).
		$
		If $\Lambda$ is realized by a negative eigenvalue $-\lambda$, then
		$
		\Lambda'
		=
		\gamma\Lambda^2+O(\delta_0\Lambda^2).
		$
		In either case, $\Lambda' \geq -C\Lambda^2$. Consequently,
		\begin{equation}\label{eq:operator-norm-lower-bound}
			\Lambda(\tau)
			\geq
			\frac{1}{
				\Lambda(m_1)^{-1}+C(\tau-m_1)
			}.
		\end{equation}
		
		We claim that $T_*=\infty$. Suppose otherwise, then
		$
		\|v(T_*)\|
		=
		R_*\omega(T_*)^{1/2}.
		$
		Lemma \ref{lem:coarse-neutral-reduction}  and \eqref{eq:coarse-neutral-at-m1} show that
		$
		|\mathrm{M}(T_*)|\approx 
		R_*\omega(T_*)^{1/2}.
		$
		Using Lemma \ref{LmMatrix_0} again,
		\[
		\pr_\tau(\Lambda/\omega^{1/2})
		\geq
		-C\Lambda
		+
		\frac{K_0^2\kappa}{5}
		\tau^{2\kappa-1}.
		\]
		At the possible first exit time $T^*$, $\Lambda=O(\omega^{1/2})\approx e^{-\tau^{\vartheta/4}/5}\ll\tau^{2\kappa-1}$. Hence the right-hand side is positive for all sufficiently large $\tau$. This contradicts the fact that $T_*$ is the first exit time. Therefore $T_*=\infty$.
		
		In addition, \eqref{eq:operator-norm-lower-bound} and the decay of $\omega(\tau)$ imply $\Lambda(\tau)/\omega(\tau)^{1/2}\to \infty$, which gives \eqref{eq:large-boundary-branch}. If $v$ remains outside $\cK_{\geq0}(\alpha_*)$, repeated using \eqref{eq:normalized-stable-contraction} shows that $R_m$ remains bounded, so \eqref{eq:small-boundary-branch} holds.
		
		Finally, if the second alternative holds, we estimate the positive, negative, and neutral parts separately. We first prove a crude bound for $\mathrm M$. Since $T_*=\infty$, the inequality $|\mathrm M'+\gamma\mathrm M^2|\leq\delta_0|\mathrm M|^2$ holds on $[m_1,\infty)$, and $\mathrm M(\tau)\to0$ because $\|v(\tau)\|\to0$. Taking the trace and using $\tr(\mathrm M^2)=|\mathrm M|^2$, we get $(\tr\mathrm M)'\leq-(\gamma-\delta_0)|\mathrm M|^2\leq0$. Hence $\tr\mathrm M(\tau)\geq\lim_{\sigma\to\infty}\tr\mathrm M(\sigma)=0$. Therefore the sum of the absolute values of the negative eigenvalues of $\mathrm M$ is at most the sum of the positive ones, and $|\mathrm M|\leq C_k\lambda_{\max}$, where $\lambda_{\max}\geq0$ is the largest eigenvalue of $\mathrm M$. By \eqref{eq:operator-norm-lower-bound}, $\lambda_{\max}>0$. At almost every time, $\lambda_{\max}'\leq-\gamma\lambda_{\max}^2+\delta_0|\mathrm M|^2\leq-\frac{\gamma}{2}\lambda_{\max}^2$ after decreasing $\delta_0$. Integrating this inequality gives
		\begin{equation}\label{eq:operator-norm-upper-bound}
			|\mathrm M(\tau)|\leq C_k\lambda_{\max}(\tau)\leq\frac{2C_k}{\gamma(\tau-m_1)}\leq\frac{C}{\tau},\qquad\tau\geq2m_1.
		\end{equation}
		Together with \eqref{eq:operator-norm-lower-bound}, we obtain $|\mathrm M(\tau)|\approx\tau^{-1}$, and by \eqref{eq:coarse-neutral}, also $\|v(\tau)\|\leq C\|\Pi_0v(\tau)\|\leq C|\mathrm M(\tau)|\leq C\tau^{-1}$.
		
		Now we iteratively apply Proposition \ref{PropApproximation}, 
		\[
		\begin{split}
			\|\Pi_{<0}v(m)\|
			&\leq
			e^{-\mu_-(m-m_0)}
			\|\Pi_{<0}v(m_0)\|+
			C\sum_{j=m_0}^{m-1}
			e^{-\mu_-(m-1-j)}
			\frac{\eps_j}{j}.
		\end{split}
		\]
		Since $\eps_j$ is eventually decreasing and has polynomial decay, 
		$
		\|\Pi_{<0}v(m)\|
		\leq
		C\frac{\eps_m}{m}.
		$
		
		For the positive component, since $\Pi_{>0}v(m)\to0$ and $\|e^{-L}|_{E^+}\|\leq e^{-\mu_+}$, backward iteration by Proposition \ref{PropApproximation} from infinity gives
		\begin{equation}\label{eq:positive-backward-series}
			\Pi_+v(m)
			=
			-
			\sum_{\ell=0}^\infty
			e^{-(\ell+1)L}
			\Pi_{>0}\big(v(m+\ell+1)-e^{L}v(m+\ell)\big).
		\end{equation}
		This implies that
		\[
		\begin{split}
			\|\Pi_+v(m)\|
			\leq
			C\sum_{\ell=0}^\infty
			e^{-\mu_+(\ell+1)}
			\frac{\eps_{m+\ell}}{m+\ell}
			\leq
			C\frac{\eps_m}{m}.
		\end{split}
		\]
		Finally,
		$
		\|\Pi_0v(m)\|
		\approx 
		|\mathrm M(m)|
		\approx
		\frac1m.
		$
		Combining all these estimates gives \eqref{eq:sharp-nonneutral-control}.
	\end{proof}

	\subsection{Behaviors of the neutral modes}
	In order to derive the refined $H^1$-normal form, we need a delicate analysis of the neutral mode.
	
	\begin{corollary}
		\label{cor:refined-neutral-matrix}
		In the alternative \eqref{eq:large-boundary-branch},
		\begin{equation}\label{eq:refined-Riccati}
			|\mathrm M'
			+\gamma\mathrm M^2|
			\leq
			C\tau^{-\vartheta}|\mathrm M|^2.
		\end{equation}
	\end{corollary}
	
	\begin{proof}
		By Lemma \ref{LmContinuation}, we may take $\alpha(\tau)=(C\bar\eps(\tau))^{-1}$ in Lemma \ref{LmMatrix_0}. Since $\bar\eps(\tau)=O(\tau^{-2\vartheta})$ and $\mathbf r(\tau)=K\tau^{\vartheta/8}$, we have
		\[
		(1+\mathbf r(\tau)^2)\bar\eps(\tau)
		=
		O(\tau^{-7\vartheta/4})
		=
		O(\tau^{-\vartheta}).
		\]
		Moreover, every cutoff term in \eqref{eq:ODE-0} is bounded by $C\tau^{-\vartheta}|\mathrm M(\tau)|^2$. This proves \eqref{eq:refined-Riccati}.
	\end{proof}
	
	\begin{lemma}
		\label{lem:matrix-Riccati-asymptotics}
		Suppose $\mathrm{M}$ is a $k\times k$-symmetric matrix-valued function defined on $[T,\infty)$, $\mathrm M(\tau)\to 0$, and assume for some $\vartheta\in(0,1)$,
		\begin{equation}\label{eq:abstract-matrix-Riccati}
			|\mathrm M'+\gamma\mathrm M^2|\leq
			C\tau^{-\vartheta}|\mathrm M|^2.
		\end{equation}
		If $|\mathrm M(\tau)|\geq c/\tau$, then there is a nonzero orthogonal projection $P_\infty$ such that
		\begin{equation}\label{eq:matrix-asymptotic}
			\mathrm M(\tau)
			=
			\frac{P_\infty}{\gamma\tau}
			+
			O(\tau^{-1-\vartheta}).
		\end{equation}
	\end{lemma}
	
	\begin{proof}
		Taking the trace gives
		\[
		(\tr\mathrm M)'
		=
		-\gamma|\mathrm M|^2+O(\tau^{-\vartheta}|\mathrm{M}|^2)
		\leq
		-\frac{\gamma}{2}|\mathrm M|^2.
		\]
		Since $\mathrm M(\tau)\to0$, $\tr\mathrm M(\tau)\to 0$. Therefore $\tr\mathrm M(\tau)\geq0$ for all sufficiently large $\tau$. Let $\lambda_{\max}$ be the largest eigenvalue of $\mathrm M$. The nonnegativity of the trace implies $|\mathrm M|\leq C_k\lambda_{\max}$. At almost every time,
		\[
		\left|
		\lambda_{\max}'
		+
		\gamma\lambda_{\max}^2
		\right|
		\leq
		C_k\tau^{-\vartheta}\lambda_{\max}^2.
		\]
		Choosing $\tau$ sufficiently large gives
		$
		-\frac{3\gamma}{2}\lambda_{\max}^2
		\leq
		\lambda_{\max}'
		\leq
		-\frac{\gamma}{2}\lambda_{\max}^2.
		$
		Integration proves $|\mathrm M(\tau)|\approx\tau^{-1}$. Consequently,
		\begin{equation}\label{eq:absolute-matrix-error}
			|\mathrm M'+\gamma \mathrm M^2|
			\leq
			C\tau^{-2-\vartheta}.
		\end{equation}

		By classical perturbation theory for linear operators, saying that for symmetric matrices, the eigenvalues may not be ordered from large to small, but they are differentiable with respect to given parameters if the matrices are differentiable; see \cite[Chapter Two, Section 6, Theorem 6.8]{Kato95_PertubLinearOp} and the proof of \cite[Lemma 3.1]{V2} in the semilinear PDE setting. We denote the eigenvalues of $\mathrm M(\tau)$ by $\lambda_1,\cdots,\lambda_k$. They are locally Lipschitz. At almost every differentiability time,
		\begin{equation}\label{eq:eigenvalue-ae}
			\left|
			\lambda_j'
			+
			\gamma\lambda_j^2
			\right|
			\leq
			C\tau^{-2-\vartheta}.
		\end{equation}
		This statement remains valid when these eigenvalues cross each other.
		
		After a change of variable $s=\log\tau$, $q_j(s)=\gamma e^s \lambda_j(e^s)$, we have 
		\begin{equation}\label{eq:scalar-logistic}
			q_j'
			=
			q_j-q_j^2+O(e^{-\vartheta s})
		\end{equation}
		for almost every $s$. The functions $q_j$ are bounded by $|\lambda_j(\tau)|\leq |\mathrm{M}(\tau)|\leq C\tau^{-1}$. As a consequence, $q_j(\tau)\to 0$ or $1$ as $\tau\to\infty$.
		
		If $q_j\to1$, then
		\[
		(q_j-1)'
		=
		-(q_j-1)-(q_j-1)^2+O(e^{-\vartheta s}),
		\]
		and variation of constants gives $(q_j-1)=O(e^{-\vartheta s}).$
		
		If $q_j\to0$, then
		\[
		q_j'-q_j
		=
		-q_j^2+O(e^{-\vartheta s}).
		\]
		Since $q_j(s)$ is bounded and tends to zero,
		\[
		q_j(s)
		=
		e^s
		\int_s^\infty
		e^{-\sigma}
		\left(
		q_j(\sigma)^2+O(e^{-\vartheta\sigma})
		\right)
		d\sigma,
		\]
		and for large $s$, $\sup_{\sigma>s}|q_j(\sigma)|\leq1/2$. Therefore
		$
		|q_j(s)|\leq \sup_{\sigma\geq s}|q_j(\sigma)|\leq Ce^{-\vartheta s}.
		$
		
		Thus every eigenvalue satisfies either
		\[
		\lambda_j(\tau)
		=
		\frac1{\gamma\tau}
		+
		O(\tau^{-1-\vartheta})
		,
		\qquad\text{or}\qquad
		\lambda_j(\tau)
		=
		O(\tau^{-1-\vartheta}).
		\]
		
		It remains to prove convergence of the associated eigenspaces. Following the idea in classical perturbation theory of linear operators (see \cite[Lemma 3.6]{V2} and the reference therein), define
		$
		\mathrm{N}(s)
		=
		\gamma e^s\mathrm M(e^s).
		$
		Then
		\begin{equation}\label{eq:matrix-logistic}
			\mathrm{N}'
			=
			\mathrm{N}-\mathrm{N}^2
			+
			O(e^{-\vartheta s}).
		\end{equation}
		For large $s$, the spectrum of $N(s)$ is contained in disjoint small neighborhoods of $0$ and $1$. Define
		\[
		P(s)
		=
		\frac1{2\pi i}
		\int_\Gamma
		(zI-\mathrm{N}(s))^{-1}\,dz,
		\qquad 
		\Gamma
		=
		\left\{
		z\in\mathbb C:
		|z-1|=\frac13
		\right\}.
		\]
		Then $P(s)$ is the orthogonal projection onto the space with eigenvalues converging to $1$.
		
		Differentiating,
		\[
		P'
		=
		\frac1{2\pi i}
		\int_\Gamma
		(zI-\mathrm{N})^{-1}
		\left(
		\mathrm{N}-\mathrm{N}^2+O(e^{-\vartheta s})
		\right)
		(zI-\mathrm{N})^{-1}\,dz.
		\]
		The contribution of $\mathrm{N}-\mathrm{N}^2$ is zero because it is a polynomial in $\mathrm{N}$. Thus
		$
		|P'|
		\leq
		Ce^{-\vartheta s}.
		$
		It follows that $P(s)\to P_\infty$ with $|P(s)-P_\infty|
		\leq
		Ce^{-\vartheta s}$. The eigenvalue estimates also imply
		$
		\mathrm{N}(s)=P(s)+O(e^{-\vartheta s}).
		$
		Therefore we get
		$
		\mathrm M(\tau)
		=
		\frac{P_\infty}{\gamma\tau}
		+
		O(\tau^{-1-\vartheta}).
		$
		The projection $P_\infty$ is nonzero because
		$|\mathrm M(\tau)|\geq c/\tau$.
	\end{proof}
	
	\subsection{The $H^1$-normal form, special case}\label{SNF}
	Let us first prove the $H^1$-normal form for a particularly chosen graphical radius $O(K \tau^{\vartheta/8})$.
	\begin{proof}[Proof of Theorem \ref{ThmNF-sqrtt} on $H^1$-normal form, with graphical radius $O(K\tau^{\kappa})$ and a particular $\vartheta$]
		We first consider the function $v$ after the Brendle-Choi rotation reduction. By Lemma \ref{LmContinuation}, we see that 
		\begin{enumerate}
			\item either $\|v(\tau)\|/\omega(\tau)^{1/2}$ remains bounded,
			\item or $\|v(\tau)\|/\omega(\tau)^{1/2}$ tends to infinity.
		\end{enumerate}
		In the former case, $\|v(\tau)\|\leq C\omega(\tau)^{1/2}\leq C\tau^{-1-\vartheta}$. In the latter case, we apply Corollary \ref{cor:refined-neutral-matrix} and Lemma \ref{lem:matrix-Riccati-asymptotics} to derive the $H^1$-normal form of $v$.
		
		It remains to go back from the graphical function modulo rotations in Section \ref{SS:BC-rotation} to the original function $u=\Phi(-A_\tau,\widetilde u)$. Throughout, $v=\chi_{\br(\tau)}\widetilde u$ is orthogonal to $\mathcal R$ for all $\tau$, and $M$, $\Pi_0$, $\Pi_{\neq0}$ refer to $v$. We claim that
		\begin{equation}\label{eq_ODE_of_Rot}
			\begin{split}
				&|A_\tau'|\leq C\omega(\tau)^{1/2}\ \text{ in the first alternative of Lemma \ref{eq:small-boundary-branch}},
				\\&|A'_\tau|\leq C\tau^{-2-\vartheta}\ \text{ in the second alternative of Lemma \ref{eq:small-boundary-branch}.}
			\end{split}
		\end{equation}
		Differentiating $\langle v,\phi_B\rangle_{L^2}=0$, $B\in\mathfrak{so}(n+1)$, in $\tau$, and using \eqref{eq:bar equation}, $L\phi_B=0$ and the self-adjointness of $L$, we get
		\begin{equation}\label{eq:rotation-projection}
			\langle\chi\mathcal A,\phi_B\rangle_{L^2}=-\langle\overline{\cQ},\phi_B\rangle_{L^2},\qquad B\in\mathfrak{so}(n+1).
		\end{equation}
		Moreover $\chi\mathcal A=\phi_{A_\tau'}+E_\tau$ to the leading order, with $\|E_\tau\|_{L^2}\leq C(\bar\eps+\omega^{1/2})|A'_\tau|$. As a consequence \eqref{eq:rotation-projection} gives
		\begin{equation}\label{eq:A'-by-Q}
			|A'_\tau|\leq C\max_{B\in \mathfrak{m},|B|=1}\big|\langle\overline{\cQ},\phi_B\rangle_{L^2}\big|.
		\end{equation}
		In the first alternative, \eqref{eq:A'-by-Q} and the Cauchy-Schwarz inequality give $|A_\tau'|\leq C\|\overline\cQ\|_{L^2}\leq C\omega^{1/2}$, as in Lemma \ref{LmA}.
		
		In the second alternative, we decompose $v=\Pi_0v+\Pi_{\neq0}v$. Since $v$ has no rotation modes, $\Pi_0v=\sum_im_{ii}H_2(y_i)+\sum_{i<j}m_{ij}H_{1,1}(y_i,y_j)$ is independent of $\theta$. Hence the quadratic part of $\cQ$ in \eqref{EqQ} satisfies $\cQ_2(\Pi_0v)=-\frac{1}{2\varrho}(\Pi_0v)^2$. As $\phi_B$ is odd in $\theta$ while $\chi$ and $d\mu$ are even in $\theta$,
		\[
		\langle\chi\cQ_2(\Pi_0v),\phi_B\rangle_{L^2}=-\frac{1}{2\varrho}\int\chi\,(\Pi_0v)^2\,\phi_B\,d\mu=0.
		\]
		Thus all the remaining nonzero terms in $\langle\overline\cQ,\phi_B\rangle$ contain a factor $\Pi_{\neq0}v$, are cubic, or are supported on $\mathbb A_\tau$. We estimate them as in the proof of Lemma \ref{LmMatrix_0}, integrating by parts once in $\theta$ and using $\Delta_\theta\phi_B=-\frac12\phi_B$ and $|\phi_B|\leq C(1+\br)$. By \eqref{eq:sharp-nonneutral-control} and $\|v\|\leq C|M|\leq C\tau^{-1}$,
		\[
		|\langle\overline\cQ,\phi_B\rangle_{L^2}|
		\leq C(1+\br^2)\big(|M|\|\Pi_{\neq0}v\|+\|\Pi_{\neq0}v\|^2+\bar\eps|M|^2\big)+C\bar\eps\omega
		\leq C(1+\br^2)\bar\eps\,\tau^{-2}+C\bar\eps\omega\leq C\tau^{-2-\vartheta},
		\]
		where we used $(1+\br^2)\bar\eps=O(\tau^{-\vartheta})$ as in the proof of Corollary \ref{cor:refined-neutral-matrix}. With \eqref{eq:A'-by-Q}, this proves \eqref{eq_ODE_of_Rot}.
		
		Integrating \eqref{eq_ODE_of_Rot} from $\tau$ to $\infty$ gives
		\[
		|A_\tau|\leq\int_\tau^\infty|A_s'|\,ds\leq C\tau^{-1-\vartheta}.
		\]
		Finally, $u=\Phi(-A_\tau,\widetilde u)=\rho_{-A_\tau}+\widetilde u+E$ (see the definition of $E$ in \eqref{eq:Phi-expansion}). Here $\|\chi\rho_{-A_\tau}\|_{H^1}\leq C|A_\tau|$, because $|\rho_A|+|\nabla\rho_A|\leq C|A|(1+|y|)$ on $\spt\chi$. We apply \eqref{eq:Phi-expansion} with $w=\widetilde u$, using $\|\widetilde u\|_{C^2(B_{\br})}\leq C\bar\eps$ and the weighted $L^2$-bound for the quadratic terms, while the annulus $\mathbb A_\tau$ contributes $O(\omega^{1/2})$. This gives
		\[
		\|\chi u-v\|_{H^1}=\|\chi (u-\widetilde u)\|_{H^1}\leq C|A_\tau|+C\omega^{1/2}\leq C\tau^{-1-\vartheta}.
		\]
		Combining this with the $H^1$-normal form of $v$ obtained above gives the $H^1$-normal form of $u$, with the error $O(\tau^{-1-\vartheta})$.
	\end{proof}
	
	Once we obtain the $H^1$-normal form of $u$ over the graphical radius $O(K\tau^{\vartheta/8})$, we can improve the error estimate to $O(\tau^{-1-\vartheta})$ for any $\vartheta\in(0,1)$ and some graphical radius $O(\tau^{\bar\kappa})$ for $\bar\kappa\in(0,1/2)$.
	
	\begin{proof}[Proof of Theorem \ref{ThmNF-sqrtt} on $H^1$-normal form, with any $\vartheta\in(0,1)$ and graphical radius $O(K\tau^{\bar\kappa})$]
		Given $\vartheta\in(0,1)$, we choose $\bar\kappa\in(0,1/2)$, $\bar\xi\in(0,1)$, and $\bar\eta>0$ so that $\bar\xi-2\bar\kappa-\bar\eta=\vartheta$. Because the special case of the $H^1$-normal form we have obtained so far showing that $\|u\|_{H^1(K\tau^{\kappa})}=O(\tau^{-1})$, we repeat the proof of Proposition \ref{thm:mainthmfor sqrtt} in Section \ref{S:Extending_t_kappa} and Corollary \ref{cor_C2_bound_for_kappa<1/2} to get that $\|u\|_{C^2(K\tau^{\bar\kappa})}=O(\tau^{2\bar\kappa-\bar\xi+\bar\eta})=O(\tau^{-\vartheta})$. Then by repeating all the proof of this section so far, we obtain the $H^1$-normal form, with any given $\vartheta$ and graphical radius $O(K\tau^{\bar\kappa})$.
	\end{proof}

	\begin{remark}\label{rmk_small_radius_better_H1}
		If one studies the region with radius $\sqrt{c\log \tau}$ with $c>4$, it seems plausible to obtain an $H^1$-normal form with error estimate $O(\tau^{-2}\log \tau)$.
	\end{remark}
	
	To obtain the $H^1$-normal form with the largest possible graphical radius $O(K\tau^{1/2})$, we need to extend the graphical radius to that size, and we will discuss that in the following sections.

	\section{The $C^1$-normal form}\label{S:C1-sqrtK}
	
	For the purpose of getting geometric properties such as the isolatedness of nondegenerate singularities, we need to upgrade the $H^1$-normal form to the $C^1$-normal form. Moreover, we can also prove that the graphical radius can be $K\sqrt{\tau}$ for any $K>0$. This larger graphical scale is not necessary in this paper, but we keep it for later purposes. 
	
	The main difficulty is as follows. In the equation of motion for $u$, we have $\pr_\tau u=Lu+\cQ(J^2u)$, where $\cQ$ is a quadratic error, (see Lemma \ref{PropQ}). However, in the radius $O(\sqrt{\tau})$, $u(\theta,y,\tau)\approx \frac{\varrho}{4\tau}\sum_{i\in\cI}h_2(y_i)$, is not small when $\cI\not=\emptyset$: when $|y|\approx K\sqrt \tau$, $u(\theta,y,\tau)\approx \frac{|\cI|\varrho}{4}K^2$. Therefore, in the $C^0$ norm, the nonlinear term $\cQ(J^2u)$ is at least $O(K^4)$, which leads to cumulative errors in the estimate of $u$ as $\tau\to\infty$. For large $K$, the error is nonperturbative. 
	The way to address the difficulty is to subtract from $u$ a term responsible for the unboundedness of $u$ near $\partial B_{K\sqrt{\tau}}$, which is mainly the zero eigenmode of $L,$ and apply Proposition \ref{PropVelazquez} to the remaining part of the graphical function. We also observe that in terms of the pseudolocality theorem, while $|u|$ is not necessarily small, $|\nabla^k u|$ is small for $k=1,2,3$, and hence the error introduced by higher-order derivatives can be nicely controlled.
	
	\subsection{The $C^1$-normal form with $K\ll1$} \label{SS_C1_K_0}
	In this section, we give the proof of Theorem \ref{ThmNF-sqrtt-C1} for $K\ll1$. Now we rewrite \eqref{EqRMCFFull} as
	$\pr_\tau u
	=
	\mathscr L( u)
	+\tilde m(u),$
	where we introduce the nonlinear operators $$\sL (u)= \frac{\varrho^2}{(\varrho+u)^2}\Delta_\theta u
	+\Delta_{\R^k} u
	-\frac12y_i\pr_i u
	+\frac{\varrho+u}{2}
	-\frac{\varrho^2}{2(\varrho+u)}.$$
	
	For the equation $\pr_\tau  f=\mathscr L (f)$, we have an approximate solution. 
	\begin{lemma}\label{Lmf} Suppose $1\leq \ell\leq k$. Let	$
		f(\theta,y,\tau)
		=
		\varrho(\Gamma(y,\tau)-1)$, where $\Gamma(y,\tau)=\sqrt{1+\frac{\sum_{i=1}^\ell (y_i^2-2)}{2\tau}}$ Then we have

		\begin{equation}\label{EqLf}
			\pr_\tau f - \sL (f)=
			\frac{\varrho}{4\tau^2\Gamma^{3}}\left(2\ell-\frac{\left(\sum_{i=1}^\ell(y_i^2-2)\right)^2}{2\tau}\right)
			=:P(y,\tau).
		\end{equation}
		Moreover,
		\begin{equation}\label{eq:f close to H2}
			\left\|f(\theta,y,\tau)-\frac{\varrho}{4\tau}\sum_{i=1}^\ell (y_i^2-2)\right\|_{H^1(\cC_{n,k})}
			=
			O(\tau^{-2}).
		\end{equation}
	\end{lemma}
	
	\begin{proof}[Proof of Lemma \ref{Lmf}]
		Direct differentiation gives
		\[
		\partial_\tau f
		=
		-\frac{\varrho \bigl(\sum_{i=1}^\ell(y_i^2-2)\bigr)}{4\tau^2\Gamma}.
		\]
		For $1\leq j\leq \ell$,
		$
		\partial_{y_j}f
		=
		\frac{\varrho y_j}{2\tau\Gamma},
		$
		whereas $\partial_{y_j}f=0$ for $j>\ell$. Similarly, for
		$1\leq j,m\leq \ell$,
		\[
		\partial_{y_jy_m}^2f
		=
		\frac{\varrho\delta_{jm}}{2\tau\Gamma}
		-
		\frac{\varrho y_jy_m}{4\tau^2\Gamma^{3}}.
		\]
		Consequently,
		\[
		\Delta_y f
		=
		\frac{\varrho\ell}{2\tau\Gamma}
		-
		\frac{\varrho \bigl(\sum_{i=1}^\ell(y_i^2-2)+2\ell\bigr)}{4\tau^2\Gamma^{3}}.
		\]
		Since $f$ is independent of $\theta$, we have $\Delta_\theta f=0$. Furthermore,
		\[
		-\frac{1}{2}y\cdot\nabla_y f
		+
		\frac{\varrho+f}{2}
		-
		\frac{\varrho^2}{2(\varrho+f)}
		=
		-\frac{\varrho\ell}{2\tau\Gamma}.
		\]
		It follows that
		\[
		\mathscr{L}(f)
		=
		-\frac{\varrho \bigl(
			\sum_{i=1}^{\ell}(y_i^2-2)+2\ell
			\bigr)}{4\tau^2\Gamma^{3}}.
		\]
		Therefore,
		\[
		\begin{aligned}
			\partial_\tau f-\mathscr{L}(f)
			&=
			\frac{\varrho}{4\tau^2\Gamma^{3}}\left(2\ell-\frac{\left(\sum_{i=1}^\ell(y_i^2-2)\right)^2}{2\tau}\right),
		\end{aligned}
		\]
		which proves \eqref{EqLf}. Finally, \eqref{eq:f close to H2} follows from the following direct computation
		\begin{equation}\label{eq:f-H2}
			f-\frac{\varrho}{4\tau}\sum_{i=1}^\ell (y_i^2-2)=-
			\frac{\varrho}{8\tau^2}(\sum_{i=1}^\ell (y_i^2-2))^2
			\left(\Gamma+1\right)^{-2}.
		\end{equation}
	\end{proof}
	
	\begin{lemma}\label{lem:profile-error}
		If we define $\mathcal E(y,\tau)
		:=\partial_\tau f-\sL(f)-\widetilde m(f)$, then 
		\begin{equation}\label{eq:complete-profile-error-bound}
			|\mathcal E(y,\tau)|+|\nabla\mathcal E(y,\tau)|
			\leq C_K\frac{1+|y|^2}{\tau^2}.
		\end{equation}
	\end{lemma}
	
	\begin{proof}
		By direct calculation, 
		\begin{equation}\label{eq:complete-profile-error}
			\mathcal E(y,\tau)
			=\frac{\varrho^4}{4\tau^2(\varrho+f)^3}
			\left[
			2|I|-\frac1{2\tau}
			\left(\sum_{i\in I}(y_i^2-2)\right)^2
			\right]+\frac{\varrho^8(\tau-\ell)\sum_{i\in I}y_i^2}
			{8\tau^4(\varrho+f)^5}
			\left(
			1+\frac{\varrho^4\sum_{i\in I}y_i^2}
			{4\tau^2(\varrho+f)^2}
			\right)^{-1}.
		\end{equation}
		Direct calculation shows that for every fixed $K>0$, on $|y|\leq K\sqrt\tau+2$ and for all sufficiently large $\tau$,
		\begin{equation}\label{eq:profile-derivative-bounds}
			\begin{aligned}
				|f|&\leq C_K\frac{1+|y|^2}{\tau},
				&|\nabla f|&\leq C_K\frac{|y|}{\tau},\\
				|\nabla^2f|&\leq C_K\left(\frac1\tau+
				\frac{|y|^2}{\tau^2}\right),
				&|\nabla^3f|&\leq C_K\left(\frac{|y|}{\tau^2}+
				\frac{|y|^3}{\tau^3}\right).
			\end{aligned}
		\end{equation}
		Plugging it into the expression of $\cE$ shows the desired estimate.
	\end{proof}
	
	\begin{proof}[Proof of Theorem \ref{ThmNF-sqrtt-C1} for $C^1$-normal form in $B_{K\sqrt{\tau}}$, with $K\ll1$]
		Fix $0<\vartheta<\vartheta_1<1$ and let $K_1>0$ and $\kappa<1/2$ be such that the graphical and $H^1$-normal-form estimates obtained in Sections \ref{S:Extending_t_kappa} and \ref{SDynamics} hold on $B_{K_1\tau^\kappa}$.
		
		Choose a sufficiently large initialization time $\tau_*$ and define $K_0=\frac{K_1}{12}\tau_*^{\kappa-1/2}$. Since the time interval used below have length $O(\log\tau_*)$, after increasing $\tau_*$ we have $6K_0\sqrt{\tau}\leq K_1\tau^\kappa.$ throughout the initial sliding window. Thus the flow is initially graphical on $B_{6K_0\sqrt{\tau}}$, with the derivative bounds required below. 
		
		We work in the region of $B_{6K\sqrt{\tau}}$. Let $ w=u-f.$ Then
		\begin{equation}\label{eq:w-exact-first}
			\partial_\tau w
			=\mathcal L(f+w)-\mathcal L(f)
			+\widetilde m(f+w)-\widetilde m(f)-\mathcal E.
		\end{equation}
		Apply the fundamental theorem of calculus along $f+s w$, $0\leq s\leq1$, 
		\begin{equation}\label{eq:w-variable-equation}
			\partial_\tau w
			=a^{AB}\nabla_{AB}w-\frac12y\cdot\nabla_yw
			+b^A\nabla_Aw+cw-\mathcal E.
		\end{equation}
		Here
		\begin{equation}\label{eq:a-average-metric}
			a^{AB}=\int_0^1\widetilde g^{AB}_{f+s w}\,ds,
		\end{equation}
		where $\widetilde g^{AB}_v$ is the inverse metric of the graph of $v$, as in the Appendix. More precisely,
		\begin{equation}\label{eq:a-angular-block}
			a^{\alpha\beta}=\int_0^1\Bigg[
			\left(1+\frac{f+s w}{\varrho}\right)^{-2}g^{\alpha\beta}
			-\frac{
				\left(1+\frac{f+s w}{\varrho}\right)^{-4}
				\nabla^\alpha(f+s w)\nabla^\beta(f+s w)}
			{1+
				\left(1+\frac{f+s w}{\varrho}\right)^{-2}
				|\nabla_\theta(f+s w)|^2+|\nabla_y(f+s w)|^2}
			\Bigg]\,ds,
		\end{equation}
		\begin{equation}\label{eq:a-Euclidean-block}
			\begin{split}
				a^{ij}=\int_0^1\Bigg[\delta^{ij}
				-\frac{\nabla_i(f+s w)\nabla_j(f+s w)}
				{1+
					\left(1+\frac{f+s w}{\varrho}\right)^{-2}
					|\nabla_\theta(f+s w)|^2+|\nabla_y(f+s w)|^2}
				\Bigg]\,ds
			\end{split}
		\end{equation}
		and
		\begin{equation}\label{eq:a-mixed-block}
			\begin{split}
				a^{\alpha i}=-\int_0^1
				\frac{
					\left(1+\frac{f+s w}{\varrho}\right)^{-2}
					\nabla^\alpha(f+s w)\nabla_i(f+s w)}
				{1+
					\left(1+\frac{f+s w}{\varrho}\right)^{-2}
					|\nabla_\theta(f+s w)|^2+|\nabla_y(f+s w)|^2}
				\,ds.
			\end{split}
		\end{equation}
		These formulas make explicit that every Hessian term involving $w$ is part of the uniformly elliptic operator. In particular, the portion coming only from $\mathcal L$ is
		\begin{equation}\label{eq:exact-L-difference}
			\mathcal L(f+w)-\mathcal L(f)
			=\frac{\varrho^2}{(\varrho+f+w)^2}\Delta_\theta w
			+\Delta_yw-\frac12y\cdot\nabla_yw
			+\left(
			\frac12+
			\frac{\varrho^2}
			{2(\varrho+f+w)(\varrho+f)}
			\right)w.
		\end{equation}
		The mean value theorem and one spatial derivative give the quantitative estimate
		\begin{equation}\label{eq:abc-quantitative}
			\|a-g^{-1}\|_{C^1(B_{6K\sqrt\tau})}
			+\|b\|_{C^1(B_{6K\sqrt\tau})}
			+\|c-1\|_{C^1(B_{6K\sqrt\tau})}
			\leq C\left(
			\|u\|_{C^3(B_{6K\sqrt\tau})}
			+\|f\|_{C^3(B_{6K\sqrt\tau})}
			\right).
		\end{equation}
		On $B_{6K\sqrt\tau}$, Lemma \ref{lem:profile-error} gives
		$\|f\|_{C^3}\leq C(K^2+\tau^{-1/2})$.
		
		Next we introduce the cut-off $\chi=\chi(y,\tau)$, with $\chi=1 \text{ on }B_{\frac92K\sqrt\tau},$ and $
		\chi=0 \text{ outside }B_{5K\sqrt\tau}$. We have
		\begin{equation}\label{eq:cutoff-equation}
			(\partial_\tau-L_\tau)(\chi w)
			=(c-1)\chi w-\chi\mathcal E
			+\left(
			\partial_\tau\chi-a^{AB}\nabla_{AB}\chi
			+\frac12y\cdot\nabla_y\chi-b^A\nabla_A\chi
			\right)w\\
			-2a^{AB}(\nabla_A\chi)(\nabla_Bw).
		\end{equation}
		Lemma \ref{lem:profile-error} shows that the right-hand side of \eqref{eq:cutoff-equation}, apart from $(c-1)\chi w$, is bounded by
		\begin{equation}\label{eq:cutoff-source-bound}
			C\left(
			\frac{1+|y|^2}{\tau^2}+|y|\check\chi(y,\tau)
			\right).
		\end{equation}
		Thus, if we define $Z=|\chi w|$, by Kato's inequality, we have 
		\begin{equation}\label{eq:Z-scalar-inequality}
			(\partial_\tau-L_\tau)Z
			\leq C\eps_0Z
			+C\left(
			\frac{1+|y|^2}{\tau^2}+|y|\check\chi(y,\tau)
			\right).
		\end{equation}
		Choose
		\begin{equation}\label{eq:exponent-choice}
			\vartheta<\vartheta_1<1,
			\qquad
			0<\eta<\min\{\vartheta_1-\vartheta,1-\vartheta\}.
		\end{equation}
		Take $K\ll1$, and hence we can take $\eps_0$ smaller than the constant in Proposition \ref{prop:perturbed-OU}. $L_\tau$ satisfies the condition in Section \ref{SS:perturbedOU}. Then we can use the perturbed Ornstein-Uhlenbeck regularization, together with the $H^1$-decay of $\chi w$, to show that 
		\begin{equation}\label{eq:C1-application-OU}
			\|w(\cdot,\tau)\|_{C^0(B_{2K_0\sqrt{\tau}})}
			\leq
			C\tau^{-\vartheta_1+\eta}\leq
			C\tau^{-\vartheta},
		\end{equation}
		Here the input is $\|\chi w\|_{L^2}=O(\tau^{-1-\vartheta_1})$, so Proposition \ref{prop:perturbed-OU} is used with $\xi=1+\vartheta_1<2$. Then the conclusion follows by choosing $\eta<\vartheta_1-\vartheta$.
		
		It remains to estimate the gradient $|\nabla w|$. Fix $(\theta,y)$ with $|y|\leq K\sqrt\tau$ and consider, for $s\in[\tau-1,\tau]$, a fixed ball centered at $(\theta,e^{-(\tau-s)/2}y)$. For large $\tau$, this entire parabolic neighborhood is contained in $B_{2K\sqrt s}$.  In coordinates moving with the center, the constant part of the drift cancels:
		\[
		-\frac12\left(e^{-(\tau-s)/2}y+z\right)\cdot\nabla_z
		+\frac12e^{-(\tau-s)/2}y\cdot\nabla_z
		=-\frac12z\cdot\nabla_z.
		\]
		Thus the drift is uniformly bounded on this fixed coordinate ball. The coefficients are uniformly parabolic and have uniformly bounded $C^1$ spatial norms by \eqref{eq:abc-quantitative} and pseudolocality; no time derivative of a coefficient is required.  Fix $p>n+2$.  The interior $W^{2,1}_p$-estimate\footnote{Here the Sobolev norm is the classical one, not the Gaussian weighted one.} for uniformly parabolic nondivergence equations with coefficients Lipschitz in space and measurable in time, followed by parabolic Sobolev embedding, gives
		\begin{equation}\label{eq:interior-gradient}
			|\nabla w(\theta,y,\tau)|
			\leq C\Bigg(
			\sup_{\substack{\tau-1\leq s\leq\tau\\
					|z-e^{-(\tau-s)/2}y|\leq1}}
			|w(\eta,z,s)|
			+
			\sup_{\substack{\tau-1\leq s\leq\tau\\
					|z-e^{-(\tau-s)/2}y|\leq1}}
			|\mathcal E(z,s)|
			\Bigg).
		\end{equation}
		We used the supremum of $|\mathcal E|$ in place of its $L^p$ norm because the parabolic cylinder has fixed volume.  The first supremum is $O(\tau^{-\vartheta})$ by the preceding $C^0$ estimate, applied at each time $s$.  By \eqref{eq:complete-profile-error-bound}, the second is $O_K(\tau^{-1})$.  Since $\vartheta<1$,
		\[
		|\nabla w(\theta,y,\tau)|
		\leq C\left(\tau^{-\vartheta}+\tau^{-1}\right)
		\leq C\tau^{-\vartheta}.
		\]
		Together with the $C^0$ estimate, this proves
		
		\begin{equation}\label{eq:C1-final}
			\left\|
			u(\theta,y,\tau)
			-\varrho\left(
			\sqrt{1+\frac1{2\tau}\sum_{i\in I}(y_i^2-2)}-1
			\right)
			\right\|_{C^1(B_{K\sqrt\tau})}
			\leq C\tau^{-\vartheta}.
		\end{equation}
		
		The above discussion not only established the $C^1$-normal form, but also extended the graphical radius to $K\sqrt{\tau}$ for a fixed time range. It remains to show that the estimate persists for all future times. Let $T^*$ be the supremum of the times for which the flow is graphical on $B_{6K_0\sqrt{\tau}}$, the coefficient bounds above hold, and $\|w(\cdot,\tau)\|_{C^1(B_{K_0\sqrt{\tau}})}\leq 2C\tau^{-\vartheta}$. The initialization above shows that $T^*>\tau_*$. Repeated finite-time pseudolocality supplies the buffered graphical region on every sliding window contained in $[\tau_*,T^*]$. The OU and interior parabolic estimates improve the bound to $\|w(\cdot,\tau)\|_{C^1(B_{K_0\sqrt{\tau}})}\leq C\tau^{-\vartheta}$. Thus all bootstrap inequalities are strict at a possible finite terminal time. Pseudolocality extends them beyond $T^*$, a contradiction. Hence $T^*=\infty$.
		
		This completes the proof of Theorem \ref{ThmNF-sqrtt-C1} for $C^1$-normal form, with $K\ll1$.
	\end{proof}
	
	\begin{proof}[Proof of Corollary \ref{cor_main_small_radius}]
		We simply repeat the above proof. When $\kappa\in(0,1/2)$, \eqref{eq:C1-application-OU} becomes
		\begin{equation}
			\|w(\cdot,\tau)\|_{C^0(B_{\tau^{\kappa}})}
			\leq
			C\tau^{2\kappa-2+\eta}.
		\end{equation}
		To get an improved gradient bound, in \eqref{eq:interior-gradient}, we use \eqref{eq:complete-profile-error-bound} with $|y|\leq \tau^{\kappa}$ to get 
		\begin{equation}
			|\nabla w(\theta,y,\tau)|\leq C\big(\tau^{2\kappa-2+\eta}+\tau^{2\kappa-2}\big)\leq C\tau^{2\kappa-2+\eta}.
		\end{equation}
		This concludes the proof.
	\end{proof}

	\subsection{Extension of the graphical scale to $K\sqrt{\tau}$ for all $K>0$}
	In this section, we show that once we obtain the $C^1$-normal form of the rescaled mean curvature flow inside the ball of radius $K_0\sqrt{\tau}$ for some $K_0>0$, we can extend the $C^1$-normal form to a larger scale.

	We fix $K_0\ll1$ from Section \ref{SS_C1_K_0}. The starting point is an extension of the graphical radius. Notice that in a nondegenerate direction, near the boundary of the ball of radius $K\sqrt{\tau}$, the $C^1$-normal form yields the graph function roughly $\varrho(\sqrt{1+K^2/2}-1)$, which is larger than any given $\eps$ if $K$ is sufficiently large. Therefore, we need to allow the graph function to have a possibly large $C^0$-norm. On the other hand, we still hope the higher-order derivatives are small, to ensure the nonlinear terms are small.
	
	We start with an application of the pseudolocality theorem.
	
	\begin{theorem}\label{thm:K0sqrt t to K sqrt t}
		For any $K>K_0$, $\bar\eps>0$, and $\vartheta\in(0,1)$, suppose for $\tau\geq T$, $M_\tau$ is a rescaled mean curvature flow that is the graph of a function $u(\cdot,\tau)$ over $ \cC_{n,k}\cap B_{K_0\sqrt{\tau}}$, and $\mathcal I\subset\{1,2,\cdots,k\}$, such that 
		\[
		\left\| u(\theta,y,\tau)-\left(\varrho\sqrt{1+\frac{\sum_{i\in\mathcal I}(y_i^2-2)}{2\tau}}-\varrho\right)
		\right\|_{C^1(\cC_{n,k}\cap B_{K_0\sqrt{\tau}})}=O(\tau^{-\vartheta}),\quad
		\tau\to\infty.
		\]
		Then there exists $\overline{T}=\overline{T}(\vartheta,K)>T$ such that $M_\tau$ is a graph of $u(\cdot,\tau)$ over $ \cC_{n,k}\cap B_{K\sqrt{\tau}}$ when $\tau>\bar T$, and 
		\begin{equation}\label{eq:Du and D2u bound}
			\sup_{ \cC_{n,k}\cap B_{K\sqrt{\tau}}}
			\sum_{j=0}^3|\nabla^j u(\cdot,\tau)|\leq C_K,
			\quad 
			\sup_{ \cC_{n,k}\cap B_{K\sqrt{\tau}}}
			\sum_{j=1}^3|\nabla^j u(\cdot,\tau)|\leq \bar\eps,
		\end{equation}
		where $C_K$ depends only on $n,k,K,K_0$.
	\end{theorem}

	\begin{proof}
		Let $f=\varrho(\Gamma-1)$ be the profile of Lemma \ref{Lmf} with $\{1,\ldots,\ell\}$ replaced by $\mathcal I$ ($f=0$ if $\mathcal I=\emptyset$), so that the assumption becomes $\|u-f\|_{C^1(B_{K_0\sqrt\tau})}=O(\tau^{-\vartheta})$. We treat the cases $\mathcal I=\emptyset$, $\mathcal I=\{1,\ldots,k\}$ and $\emptyset\neq\mathcal I\subsetneqq\{1,\ldots,k\}$ at once: near a point of the axis, the flow is close to the cylinder of radius $\varrho\Gamma$ given by the normal form, and we propagate it for a fixed rescaled time by Lemma \ref{lem:pseudo-radius}.
		
		Fix $K>K_0$ and $\bar\eps>0$, and we define the quantities
		\[
		h:=2\log\frac{2K}{K_0},\qquad t_1:=1-e^{-h},\qquad b_-:=\varrho\sqrt{1-\tfrac12e^{-h}},\qquad b_+:=\varrho\sqrt{1+K_0^2},
		\]
		so that $\varrho^2t_1<b_-^2$. Let $R_1,\eta_1$ be given by Lemma \ref{lem:pseudo-radius} for $b_\pm$, $\sigma=t_1$ and $\eps\in(0,1)$ to be chosen below, and fix $R>\max\{R_1,2(b_++1+\sqrt{2n})\}$. Let $s=\tau-h$. We work directly with $\bM_t$ on $[-e^{-s},-e^{-\tau}]$, using
		\[
		\bM_{-e^{-s}}=e^{-s/2}M_s,\qquad
		\bM_{-e^{-\tau}}=e^{-\tau/2}M_\tau,\qquad
		(-e^{-\tau})-(-e^{-s})=e^{-s}t_1.
		\]
		
		Let $(\theta,y)\in\cC_{n,k}$ with $|y|\leq K\sqrt\tau$, and put $y_0=e^{-h/2}y$ and $p=(0,y_0)$. Then $|y_0|\leq\frac{K_0}{2}\sqrt\tau$, so $B_R(p)\subset B_{K_0\sqrt s}$ and, by assumption, $M_s\cap B_R(p)$ is the graph of $u(\cdot,s)$. Let $b=\varrho\Gamma(y_0,s)$. Then $b_-\leq b\leq b_+$, because $\varrho^2(1-k/s)\leq b^2\leq\varrho^2(1+|y_0|^2/2s)$. Since $b-\varrho=\varrho(\Gamma(y_0,s)-1)=f(y_0,s)$, $M_s\cap B_R(p)$ is a normal graph over $\bS^{n-k}(b)\times\R^k$ of
		\[
		w_0:=(u-f)(\cdot,s)+f(\cdot,s)-f(y_0,s),
		\]
		and \eqref{eq:profile-derivative-bounds} gives $\|w_0\|_{C^1(B_R(p))}\leq C(s^{-\vartheta}+(1+R)s^{-1/2})\leq\eta_1$. Apply Lemma \ref{lem:pseudo-radius} to $\bM_t$ on $[-e^{-s},-e^{-\tau}]$, with lengths measured in units of $e^{-s/2}$ and elapsed time in units of $e^{-s}$. It follows that $\bM_{-e^{-\tau}}\cap B_{e^{-s/2}R/2}(e^{-s/2}p)$ is a normal graph over $\bS^{n-k}(e^{-s/2}b')\times\R^k$, where $b'=\sqrt{b^2-\varrho^2t_1}$. In these spatial units it has the parametrization
		\[
		\left(e^{-s/2}\frac{b'+w(\theta,z)}{\varrho}\theta,e^{-s/2}z\right),\qquad \|w\|_{C^3}\leq\eps.
		\]
		Dilation by $e^{\tau/2}$ sends $e^{-s/2}p$ to $(0,y)$ and this reference cylinder to $\bS^{n-k}(e^{h/2}b')\times\R^k$, where
		\[
		e^{h/2}b'=\sqrt{\varrho^2+e^h(b^2-\varrho^2)}\leq e^{h/2}b_+.
		\]
		Hence $M_\tau\cap B_{e^{h/2}R/2}((0,y))$ is the graph over $\cC_{n,k}$ of
		\[
		u(\theta,y',\tau):=e^{h/2}b'-\varrho+e^{h/2}w(\theta,e^{-h/2}y').
		\]
		Since $R/2>b_++1$, the point of this graph over $(\theta,y)$ lies in the ball, and 
		\[
		|u(\theta,y,\tau)|\leq e^{h/2}(b_++1),\qquad\sum_{j=1}^3|\nabla^ju(\theta,y,\tau)|\leq C\eps.
		\]
		Then choosing $\eps$ with $C\eps\leq\bar\eps$ gives \eqref{eq:Du and D2u bound} at $(\theta,y)$ with $C_K=e^{h/2}(b_++1)$. 
		
	\end{proof}

	Now we are ready to prove the $C^1$-normal form inside the graphical radius $K\sqrt{\tau}$ for arbitrary $K>0$. The key observation is that, although the linearized operator $L_\tau$ of the difference function $w=u-f$ can not be viewed as a perturbative operator of $L_{\cC_{n,k}}$ when $K$ is large, it is still a uniformly parabolic operator. Then, we can use a translating argument as in Corollary \ref{cor_C2_bound_for_kappa<1/2} with the parabolic estimate to extend the estimate.

	\begin{proof}[Proof of Theorem \ref{ThmNF-sqrtt-C1} ($C^1$-normal form in $B_{K\sqrt{\tau}}$)]
		From Theorem \ref{thm:K0sqrt t to K sqrt t}, $u$ and hence $w$ can be defined on $B_{2K\sqrt{\tau}}$. \eqref{eq:w-exact-first} gives the explicit PDE that $w$ satisfies, and 
		\begin{equation}\label{eq:large-K-coefficient-bounds}
			C_K^{-1}|\zeta|^2
			\leq a^{AB}\zeta_A\zeta_B
			\leq C_K|\zeta|^2,
			\qquad 
			\|a\|_{C^1}+\|b\|_{C^1}+\|c\|_{C^0}\leq C_K.
		\end{equation}
		Note that $C_K$ depends on $K$ but is a fixed constant. Also, from Lemma \ref{lem:profile-error}, when $|y|\leq 2K\sqrt{\tau}$, we have $|\cE(y,\tau)|\leq C_K/\tau$. We fix $\overline{T}$ from Theorem \ref{thm:K0sqrt t to K sqrt t} as the time range to run the parabolic PDE. 
		
		To estimate $w$ in a larger domain, we need the following translational argument. Define
		\begin{equation}\label{eq:large-K-moving-coordinates}
			v(\theta,z,\tau)
			:=w\left(
			\theta,e^{(\tau-s)/2}(y_0+z),\tau
			\right).
		\end{equation}
		We will choose a cutoff in the $z$ variable as follows. First, let $R_\tau$ be a sufficiently large fixed multiple of $\sqrt{\log \tau}$, which is increasing but also very small compared with $\sqrt{\tau}$. Then we choose the cut-off function $\chi(z,\tau)$ so that it is $1$ on $|z|\leq2R_\tau$, $0$ on $|z|\geq3R_\tau$, and satisfies
		\[
		|\nabla\chi|\leq C R_\tau^{-1},
		\quad
		|\nabla^2\chi|\leq C R_\tau^{-2}.
		\]
		Multiplying with the cut-off $\chi$, 
		\begin{equation}\label{eq:large-K-cutoff-equation}
			(\partial_\sigma-a^{AB}\nabla_{AB}
			-b^A\nabla_A-c)(\chi v)=-\chi\mathcal E
			-2a^{AB}(\nabla_A\chi)(\nabla_Bv)
			-\left(a^{AB}\nabla_{AB}\chi
			+b^A\nabla_A\chi\right)v.
		\end{equation}
		At the initial time, the support of $\chi$ is contained in $|y-y_0|<3R_\tau\subset B_{K_0\sqrt \tau}$, so the $C^1$-normal form for small $K_0$ bounds $|\chi v|$ by $C\tau^{-\vartheta}$. The contribution of $\cE$ is at most $C_K/\tau$. All cutoff commutators are supported where $2R_\tau\leq|z|\leq3R_\tau$ and, by \eqref{eq:Du and D2u bound}, are bounded by $C_K$. Evaluating at $z=0$, apply the variation of constants (e.g., use the heat kernel expression as in Section \ref{SS:perturbedOU}; the heat kernel expression does not rely on the closeness of $L_\tau$ and $L_{\cC_{n,k}}$) to the uniformly parabolic PDE \eqref{eq:large-K-cutoff-equation},
		\[
		|v(\theta,0,\tau)|
		\leq C_K\tau^{-\vartheta}+\frac{C_K}{\tau}+C_K(1+R_\tau)^k
		\exp\left(-\frac{R_\tau^2}{C_K\overline{T}}\right).
		\]
		The last term is the integral of the heat kernel over the cutoff annulus. Since $R_{\tau}$ is a fixed multiple of $\sqrt{\log\tau}$, choosing that multiple large makes the Gaussian tail smaller than any prescribed inverse power of $\tau$. Thus we have
		\begin{equation}\label{eq:finite-time-C0-calculation}
			|w(\theta,y,\tau)|
			\leq C_K\left(
			\tau^{-\vartheta}+\tau^{-1}
			+\exp\left(-\frac{R_\tau^2}{C_K\overline{T}}\right)
			\right).
		\end{equation}
		Choose the fixed multiple in the definition of $R_\tau$ sufficiently large that $\frac{R_\tau^2}{C_K\overline{T}}\geq3\log \tau$, we obtain, uniformly on $B_{2K\sqrt\tau}$,
		\begin{equation}\label{eq:finite-time-C0-final}
			|w(\theta,y,\tau)|\leq C_K\tau^{-\vartheta}.
		\end{equation}
		With the $C^0$-estimate, we can apply the interior nondivergence (classic) $W^{2,1}_p$-estimate, by shifting the time backwards a little bit. For $p>n+2$, the parabolic Schauder estimate shows that
		\begin{equation}
			|\nabla w(\theta,y,\tau)|
			\leq C_K\left(
			\|w\|_{L^\infty(Q_1)}
			+\|\mathcal E\|_{L^p(Q_1)}
			\right)
			\leq C_K\left(
			\tau^{-\vartheta}+\tau^{-1}
			\right).
		\end{equation}
		Here $Q_1$ is a fixed unit parabolic cylinder contained in the moving region. Together with \eqref{eq:finite-time-C0-final}, this proves Theorem \ref{ThmNF-sqrtt-C1} for arbitrarily large $K$.
	\end{proof}

	\begin{proof}[Proof of Theorem \ref{ThmNF-sqrtt}]
		We have proved the $H^1$-normal form in $B_{K_0\sqrt{\tau}}$ for $K_0\ll1$. For general $K>0$, we only notice that the $C^1$-normal form on the annulus region $B_{K\sqrt{\tau}}\backslash B_{K_0\sqrt{\tau}}$ contributes exponentially decay error, which is much smaller than $\tau^{-1-\vartheta}$. This proves the $H^1$-normal form. The moreover part is proved in Proposition \ref{prop:decay order zero normal form}.
	\end{proof}

	\section{Consequences for the $C^1$-normal form}\label{S_Geometric}
	\subsection{Isolatedness, Mean convex neighborhood, and Type-I}\label{SS:Mean convex neighborhood and Type-I}
	In this section, we shall apply our $C^1$-normal to study the local geometry of nondegenerate singularities. Throughout this section the singularity is nondegenerate, $\cI=\{1,2,\ldots,k\}$, and we write
	\[
	\Gamma=\sqrt{1+\frac{|y|^2-2k}{2\tau}},\qquad f=\varrho(\Gamma-1).
	\]
	We also define $g(s)=\sqrt{s}e^{-s/2}$. 
	
	We first compute the curvature of a graph over $\cC_{n,k}$ whose $C^0$-norm is not small. Let $\Sigma=\{(\frac{r\theta}{\varrho},y)\}$ be the graph of a function $r$ over a domain in $\cC_{n,k}$ with $0<r_-\leq r\leq r_+$ and $\sum_{j=1}^2|\nabla^jr|\leq\eps\leq1$. The computation in the proof of Proposition \ref{PropQ}, in particular \eqref{eq:H}, can be applied to $u=r-\varrho$ and only requires $r>0$. Hence
	\begin{equation}\label{eq:graph-curvature}
		\Big|H-\frac{n-k}{r}\Big|+\Big||A|-\frac{\sqrt{n-k}}{r}\Big|\leq C_4\eps,
	\end{equation}
	where $C_4$ depends only on $n,k,r_-,r_+$.
	
	Fix $\vartheta\in(0,1)$. By Theorem \ref{ThmNF-sqrtt-C1} and Theorem \ref{thm:K0sqrt t to K sqrt t} with $K=3$, there is a constant $C_3$, independent of $\bar\eps$ (see the proof of Theorem \ref{thm:K0sqrt t to K sqrt t}), such that for every $\eps\in(0,1)$ there is $\tau_1>1$ with the following property: for every $\tau\geq\tau_1$, $M_\tau\cap B_{3\sqrt\tau}$ is the graph of a function $u(\cdot,\tau)$ over $\cC_{n,k}\cap B_{3\sqrt\tau}$ with
	\begin{equation}\label{eq:geometric-input}
		\|u-f\|_{C^1}\leq C\tau^{-\vartheta},\qquad|u|\leq C_3,\qquad\sum_{j=1}^3|\nabla^ju|\leq\eps.
	\end{equation}
	Since $\Gamma\geq\sqrt{1-k/\tau}\geq\frac34$ for large $\tau$, we may also assume $\varrho+u\geq\frac\varrho2$.
	
	In the following, we use $\bH$ and $\mathbf{A}$ to denote the mean curvature and second fundamental form of the associated mean curvature flow $\bM_t$.
	
	\begin{proposition}\label{prop:geometric-region}
		Assume that the singularity is nondegenerate. Then there exist $\delta_2,t_0,C>0$ such that
		\begin{enumerate}
			\item for every $t\in(-t_0,0)$ and every $x\in\bM_t\cap B_{\delta_2}$, we have $\bH(x,t)>0$ and $|\mathbf A|(x,t)\leq C(-t)^{-1/2}$;
			\item for every $x\in B_{\delta_2}\setminus\{0\}$, $(x,0)$ is a regular point of the flow.
		\end{enumerate}
	\end{proposition}
	
	\begin{proof}
		By \eqref{eq:geometric-input}, the graph $M_s$ in
		$B_{3\sqrt s}$ has radius bounded above and bounded away from zero, while its first three derivatives tend uniformly to zero. Thus \eqref{eq:graph-curvature} gives $H>0$ and $|A|\leq C$ for
		all sufficiently large $s$.
		
		{\bf Claim.} For all sufficiently large $\tau$, every point on the sphere $\{|x|=\sqrt{\tau}e^{-\tau/2}\}$ has a neighborhood in which the flow is smooth, possibly empty, for $-\frac12e^{-\tau}\leq t\leq0$, and on its nonempty slices
		\begin{equation}\label{eq:shell-curvature}
			\bH>0,\qquad |\mathbf A|\leq Ce^{\tau/2}.
		\end{equation}
		The constant $C$ is independent of $x$.
		
		\begin{proof}[Proof of claim]
			We write $\sqrt{\tau}=e^{\tau/2}x=(\xi,y)$, so that $|\xi|^2+|y|^2=\tau$. In Lemma \ref{lem:pseudo-radius}, choose $\eps\ll1$, $b\in[\varrho\sqrt{5/4},\varrho\sqrt{3/2}]$ and $\sigma=1/2$ on the time interval $[0,1]$. Fix $R$ large enough in Lemma \ref{lem:pseudo-radius}so that
			\[
			R/4>\varrho+C_3+\sqrt{2n}+2.
			\]

			If $|\xi|\geq R/4$, we use the shrinking sphere as the comparison. The ball of radius $\sqrt{2n}+1$ centered at $(\xi,y)$ lies in $B_{3\sqrt\tau}$ and is disjoint from $M_\tau$, and its distance from the $y$-axis is at most $\varrho+C_3$. In the original variables, the initially empty ball centered at $x$ has radius $e^{-\tau/2}(\sqrt{2n}+1)$ at time $-e^{-\tau}$. At time $0$, the shrinking sphere's squared radius
			\[
			e^{-\tau}\big((\sqrt{2n}+1)^2-2n\big)>e^{-\tau}.
			\]
			The avoidance principle shows that a neighborhood of $x$ is empty throughout $[-e^{-\tau},0]$.
			
			If $|\xi|<R/4$, the fixed ball $B_R((0,y))$ lies in $B_{3\sqrt\tau}$ for large $\tau$. Let $b=\varrho\sqrt{1+\frac{|y|^2-2k}{2\tau}}$. Since $|y|^2=\tau-|\xi|^2$, we have $5\varrho^2/4\leq b^2\leq3\varrho^2/2$ for large $\tau$. The normal form and its derivative bounds give
			\[
			\|u(\cdot,\tau)-(b-\varrho)\|_{C^1(B_R((0,y)))}
			\leq C\tau^{-\vartheta}+C(1+R)\tau^{-1/2}\to0,
			\]
			uniformly in the chosen centers. Apply Lemma \ref{lem:pseudo-radius} directly to $\bM_t$ on $[-e^{-\tau},0]$, the comparison cylinder has radius
			\[
			e^{-\tau/2}\sqrt{b^2-\varrho^2(1+e^\tau t)},
			\qquad -e^{-\tau}\leq t\leq0.
			\]
			After rescaling back to the scale of the rescaled mean curvature flow, the radius of the cylinder lies between $\varrho/2$ and $\varrho\sqrt{3/2}$, and Lemma \ref{lem:pseudo-radius} gives a small $C^2$ perturbation of this cylinder on the second half of the time interval, in the ball $B_{e^{-\tau/2}R/2}(e^{-\tau/2}(0,y))$. This ball contains $B_{e^{-\tau/2}R/4}(x)$. Thus, \eqref{eq:graph-curvature} gives positive mean curvature and a uniform curvature bound at the scale of rescaled mean curvature flow, and scaling back proves \eqref{eq:shell-curvature}. This finishes the proof of the claim.
		\end{proof}

		With the claim, choose $\tau_1>1$ so that the claim holds for every $\tau\geq\tau_1$, $\delta_2=\sqrt{\tau_1}e^{-\tau_1/2}$, and let $t_0=e^{-\tau_1}$. Since $\sqrt{s}e^{-s/2}$ is strictly decreasing on $(1,\infty)$ and tends to zero, every $0<|x|<\sqrt{\tau_1}e^{-\tau_1/2}$ has the form $|x|=\sqrt{\tau}e^{-\tau/2}$ for a unique $\tau>\tau_1$. Then the claim proves item (2) by scaling, and also item (2).
	\end{proof}
	
	With the proposition, we can prove Theorems \ref{ThmIsolate}, \ref{prop:Mconnbhd} and \ref{ThmtypeI}.
	
	\begin{proof}[Proof of Theorem \ref{ThmIsolate}]
		Let $\delta_1=\delta_2$ in proposition \ref{prop:geometric-region}. By assumption the flow is smooth for $t<0$, and by Proposition \ref{prop:geometric-region} (2), $(x,0)$ is a regular point for every $x\in B_{\delta_1}\setminus\{0\}$. Hence $(0,0)$ is the only singularity of $\bM_t$ in $B_{\delta_1}\times(-\delta_1^2,0]$.
	\end{proof}

	\begin{proof}[Proof of Theorem \ref{prop:Mconnbhd}]
		By Proposition \ref{prop:geometric-region} (1), $\bM_t\cap B_{\delta_2}$ has positive mean curvature for every $t\in(-t_0,0)$. For the topology, we note that all the local descriptions in the proof of the proposition are graphs over $\bS^{n-k}(\varrho)\times\R^k$ in the same coordinates, with gradient at most $\eps\leq1$, and that each of them contains all points of $\bM_t$ in the corresponding ball. By the uniqueness of graph functions, they glue to a radial function $r_t$. After decreasing $t_0$, $r_t\leq|y|/10$ whenever $|y|\geq\delta_2/2$, and $r_t\leq\delta_2/10$ everywhere. Since $|\nabla r_t|\leq1$, the function $r_t^2+|y|^2$ is increasing in $|y|$ for $|y|\geq\delta_2/2$ and is smaller than $\delta_2^2$ for $|y|\leq\delta_2/2$. Therefore $\bM_t\cap B_{\delta_2}$ is the graph over $\{(\theta,y):|y|<\rho_t(\theta,y/|y|)\}$, for a smooth function $\rho_t\in(\delta_2/2,\delta_2)$. This is diffeomorphic to $\bS^{n-k}\times B^k$, and hence to $\cC_{n,k}\cap B_{4n}$.
	\end{proof}
	
	\begin{proof}[Proof of Theorem \ref{ThmtypeI}]
		By Proposition \ref{prop:geometric-region} (1), $|\mathbf A|(x,t)\leq C(-t)^{-1/2}$ for every $t\in(-t_0,0)$ and $x\in\bM_t\cap B_{\delta_2}$.
	\end{proof}
	
	\subsection{Vanishing of zero modes}
	The result of this subsection is not used anywhere in this paper, but will be used in \cite{SunWangXue2_RegSing}. We have seen that the $H^1$-normal form of the rescaled mean curvature flow is given by $\sum_{i\in\cI}\frac{\varrho}{4\tau}(y_i^2-2)+O(\tau^{1(1+\vartheta)})$, where $\cI\subset\{1,2,\cdots,k\}$. When $\cI=\emptyset$, we expect that the next eigenmodes of $L_{\cC_{n,k}}$ that have strictly positive eigenvalues should dominate the evolution. Then we should expect a faster decay of the graph function.
	
	\begin{proposition}[``Moreover'' part of Theorems \ref{ThmNF-sqrtt} and \ref{ThmNF-sqrtt-C1}]\label{prop:decay order zero normal form}
		Suppose $K>0$, the rescaled mean curvature flow is a graph of the function $u(\cdot,\tau)$ over $\cC_{n,k}\cap B_{K\sqrt{\tau}}$, and $\cI=\emptyset$ in the $H^1$-normal form. Let $K_0\in(0,K)$ be such that $K_0^2<10\mu_-$, that is, there is no eigenvalue of $-L_{\cC_{n,k}}$ in the interval $(0,K_0^2/10]$. Then for $\tau$ sufficiently large,
		\begin{equation}\label{eq:degenerate-H1-C1}
			\|u(\cdot,\tau)\|_{H^1(\cC_{n,k}\cap B_{K_0\sqrt{\tau}-1})}\leq Ce^{-K_0^2\tau/10},\qquad\|u(\cdot,\tau)\|_{C^1(\cC_{n,k}\cap B_{K_0\sqrt{\tau}/2})}\leq Ce^{-K_0^2\tau/16}.
		\end{equation}
	\end{proposition}
	
	\begin{proof}
		After we extend the graphical radius of $u$ to $K\sqrt{\tau}$, we may choose $\br(\tau)=K_0\sqrt{\tau}$ in Section \ref{SDynamics}, and hence $\omega(\tau)=e^{-K_0^2\tau/5}$. Since $\cI=\emptyset$, Theorem \ref{thm:K0sqrt t to K sqrt t} and the $C^1$-normal form give $\|u(\cdot,\tau)\|_{C^3(B_{K\sqrt\tau})}\to0$, so the analysis of Section \ref{SDynamics} applies with this choice of $\br$. The second alternative of Lemma \ref{LmContinuation} is excluded, because it would give $|\mathrm M(\tau)|\geq c/\tau$, and since the matrices $\mathrm M$ computed with the cutoffs at radii $\tau^{\vartheta/8}$ and $K_0\sqrt\tau$ differ by a quantity decaying faster than any power of $\tau$, Lemma \ref{lem:matrix-Riccati-asymptotics} would give $\cI\neq\emptyset$. Hence the first alternative holds, which is an estimate in the weighted $H^1$-norm: $\|\overline v(\tau)\|\leq C\omega(\tau)^{1/2}=Ce^{-K_0^2\tau/10}$.
		
		For the $C^1$-estimate, we use a local parabolic estimate. Fix $(\theta_*,y_*,\tau_*)$ with $|y_*|\leq K_0\sqrt{\tau_*}/2$, and let $Q'=\{(\theta,y,\tau):\tau\in[\tau_*-2,\tau_*],\ |y-e^{-(\tau_*-\tau)/2}y_*|\leq2\}$ and $Q\subset Q'$ the subregion with $[\tau_*-1,\tau_*]$ and radius $1$. For $\tau_*$ large, $Q'$ lies in $\{|y|\leq K_0\sqrt s-1\}$. In coordinates moving with the center as in \eqref{eq:interior-gradient}, $u$ satisfies the uniformly parabolic equation $\partial_\tau u=a^{AB}(J^1u)\nabla_{AB}u-\frac12y\cdot\nabla_yu+b^A(J^2u)\nabla_Au+u$ (see the proof of Lemma \ref{LmVelaz}) with bounded drift and coefficients whose spatial $C^1$-norms are uniformly bounded. The local maximum principle for such equations gives $\sup_Q|u|\leq C\|u\|_{L^2(Q')}$, where $L^2(Q')$ is the unweighted space-time $L^2$-norm, and the interior $W^{2,1}_p$-estimate gives $|\nabla u(\theta_*,y_*,\tau_*)|\leq C\sup_Q|u|$. On $Q'$ we have $|x|^2\leq\varrho^2+(K_0\sqrt{\tau_*}/2+2)^2$, so the unweighted $L^2$-norm on $Q'$ is bounded by $\exp\big(\frac{\varrho^2+(K_0\sqrt{\tau_*}/2+2)^2}{8}\big)$ times the weighted one. Therefore, 
		\[
		\begin{split}
			|u(\theta_*,y_*,\tau_*)|+|\nabla u(\theta_*,y_*,\tau_*)|&\leq Ce^{\frac{K_0^2\tau_*}{32}+C\sqrt{\tau_*}}\sup_{s\in[\tau_*-2,\tau_*]}\|u(\cdot,s)\|_{L^2(B_{K_0\sqrt s-1})}
			\\&\leq Ce^{-(\frac1{10}-\frac1{32})K_0^2\tau_*+C\sqrt{\tau_*}}\leq Ce^{-K_0^2\tau_*/16}.
		\end{split}
		\]
		This proves the second estimate in \eqref{eq:degenerate-H1-C1}, and hence the second alternatives of Theorems \ref{ThmNF-sqrtt} and \ref{ThmNF-sqrtt-C1}.
	\end{proof}
	
	\begin{remark}
		If $K_0$ is sufficiently large, $-K_0^2/10$ will be smaller than some eigenvalue of $-L$. Then the exponential decay is given by $\|u\|_{H^1(B_{K_0\sqrt{\tau}})}\leq Ce^{(-\lambda+\eps) \tau}$ where $\lambda$ is the smallest positive eigenvalue of $-L$ and $\eps>0$ is a small number. This is almost sharp, see \cite{SunWangXue2_RegSing}.
	\end{remark}

	\appendix
	
	\section{Estimates of the nonlinear term}\label{appendix:Nonlinear estimate}	
	In this appendix, we derive the equation of motion for the graphical function $u$ of a manifold evolving under rescaled mean curvature flow approaching a cylinder $\cC_{n,k}=\bS^{n-k}(\varrho)\times \R^{k}$ where $\varrho=\sqrt{2(n-k)}$ is the radius of the sphere. We will use coordinates $z=(\theta, y)\in \cC_{n,k}$ where $\theta\in \R^{n-k+1}$ denotes point on $\bS^{n-k}(\varrho)$ and $y$ denotes a point on $\R^{k}$. Note $|\theta|=\varrho$. The following computations are locally around a point $z_0=(\theta_0,y_0)$, and we choose local orthonormal frame $\{\theta_\alpha\}$ and $\{y_i\}$. Greek letters correspond to the spherical part and $i,j,k$'s correspond to the $\R^{k}$ part. We will use $\pr_\alpha$ and $\pr_i$ to simplify $\pr_{\theta_\alpha}$ and $\pr_{y_i}$ respectively. 
	\begin{proposition}\label{PropQ}Let $M_\tau$ be a rescaled mean curvature flow converging to a cylinder $\cC_{n,k}$ in the $C^\infty_{loc}$ sense as $\tau\to\infty$. Writing $M_\tau$ as a normal graph of a function $u$ over $\cC_{n,k}$ within the graphical radius, then we have  
		\begin{enumerate}
			\item 	we have
			$\pr_\tau u=L u+\cQ(J^2u),\ J^2u:=(u,\nabla u, \nabla^2u),$
			where we have 	\begin{equation}\label{EqQ<}
				|\cQ(J^2u)|\leq 
				C(|\nabla u|^4+|\nabla u|^2|\Hess_u|+|\nabla u|^2+u^2+|u||\Hess_u|)
			\end{equation}
			if $\|u\|_{C^2}\leq \eps_0$ for some $\eps_0$ small. 
			\item 	The leading terms in $\cQ$ is given explicitly as 
			\begin{equation}\label{EqQ}
				\cQ(J^2u)=-(2\varrho)^{-1}\left(u^2+4u\Delta_\theta u+2|\nabla_\theta u|^2\right) +\cC(J^2u),	
			\end{equation}
			where $\cC(J^2u)$ consists of terms cubic and higher power in $J^2u$, satisfying
			\begin{equation}
				|\cC(J^2u)|\leq C(|u|^3+|\nabla u|^3+|\nabla^2u|(u^2+|\nabla u|^2))\leq C\eps_0(u^2+|\nabla u|^2).
			\end{equation} 
			\item For three functions $v, u_1,u_2:\ \Sigma\cap B_R\to \R$ and the cutoff function $\chi$ that is 1 on $B_{R-1}$ and $0$ outside $B_R$, we have $($denoting $w=u_1-u_2)$\begin{equation*}\label{EqdC}
				\left|\int_{\Sigma\cap B_R} v 	(\cC(J^2(\chi u_1))-	\cC(J^2(\chi u_2)))d\mu\right|
				\leq C\max\{\|u_1\|_{C^2},\|u_2\|_{C^2}\}^2(\|w\|\cdot \|v\|+e^{-\frac{R^2}{4}}).
			\end{equation*}
		\end{enumerate}

	\end{proposition}
	
	\begin{proof}
		
		Note that the unit normal vector $\bn=\frac{\theta}{\varrho}$. Later we will use the following fact:
		\[\pr_i \theta=\pr_i\pr_j\theta=0,\ \ \ 
		\pr_\alpha\theta=\theta_\alpha,\ \ \ \pr_{\alpha}\pr_{\beta} \theta = -\frac{\delta_{\alpha\beta}}{\varrho^2}\theta.
		\]
		We consider a graph $\cC_{n,k}$ locally given by $\{\widetilde{F}(z)=(z)+u(z)\bn= (z)+u(z)\frac{\theta}{\varrho}\}$, which induces a frame on $\cC_{n,k}$ given by
		$
		\widetilde{\pr}_\alpha=(1+\frac{u}{\varrho})\theta_{\alpha}+\frac{\pr_\alpha u}{\varrho} \cdot\theta,\ \ \ 
		\widetilde{\pr}_{i}=\pr_i +\frac{\pr_iu}{\varrho} \cdot \theta.
		$
		The induced metric is given by
		\begin{equation*}
			\widetilde{g}_{\alpha\beta}=(1+\frac{u}{\varrho})^2 g_{\alpha\beta} + \pr_\alpha u \pr_\beta u
			,
			\ \ \ 
			\widetilde{g}_{i j}
			=
			g_{ij}+\pr_i u\pr_j u
			,
			\ \ \ 
			\widetilde{g}_{\alpha i}
			=
			\pr_{\alpha} u \pr_i u.
		\end{equation*}
		The inverse matrix is given by
		\begin{equation*}
			\begin{aligned}
				\widetilde{g}^{\alpha\beta}
				&=
				(1+\frac{u}{\varrho})^{-2}g_{\alpha\beta} - (1+\frac{u}{\varrho})^{-4}\pr_\alpha u \pr_\beta u+ m_{\alpha\beta},
				\quad	\\
				\widetilde{g}^{i j}
				&=
				g_{ij} - \pr_i u\pr_j u 
				+ m_{ij},\quad 
				\widetilde{g}^{\alpha i}
				=
				-\pr_\alpha u\pr_i u
				+ m_{\alpha i}.
			\end{aligned}
		\end{equation*}
		Here $m$ are functions of products of $\nabla u$, at least quadratically.	We find a unit normal vector field given by
		\begin{equation*}
			\widetilde{\bn}
			=
			\frac{\frac{\theta}{\varrho}-(1+\frac{u}{\varrho})^{-1}\pr_\alpha u\theta_{\alpha} -\pr_i u \pr_i}
			{\sqrt{
					1+\sum_{\alpha}(1+\frac{u}{\varrho})^{-2}(\pr_\alpha u)^2	
					+\sum_i (\pr_i u)^2.
			}}
			=:
			\frac{\frac{\theta}{\varrho}-(1+\frac{u}{\varrho})^{-1}\pr_\alpha u\theta_{\alpha} -\pr_i u \pr_i}{S}.
		\end{equation*}
		
		Now we calculate the 2nd fundamental form. We have
		\begin{equation*}
			\begin{aligned}
				\pr_{\alpha}\pr_{\beta} F
				&=
				-(1+\frac{u}{\varrho})\frac{\delta_{\alpha\beta}}{\varrho^2}\theta
				+
				\frac{\pr_{\beta} u}{\varrho}\theta_{\alpha}
				+
				\frac{\pr_{\alpha }u}{\varrho}\theta_\beta
				+
				\frac{\pr_{\alpha\beta}u}{\varrho}\theta,
				\\
				\pr_i\pr_{\alpha} F
				&	=
				\frac{\pr_{i\alpha} u}{\varrho}\theta
				+\frac{\pr_i u}{\varrho}\theta_{\alpha},
				\quad		\pr_i\pr_j F
				=
				\frac{\pr_{ij}u}{\varrho}\theta.
			\end{aligned}
		\end{equation*}
		
		Taking inner product with $\widetilde{\bn}$, we get the 2nd fundamental forms:
		\begin{equation*}
			\begin{aligned}
				\tilde{A}_{\alpha\beta}
				&=
				S^{-1}\left(
				-(1+\frac{u}{\varrho})\frac{\delta_{\alpha\beta}}{\varrho}
				-2(1+\frac{u}{\varrho})^{-1}\frac{\pr_\alpha u \pr_{\beta}u}{\varrho}
				+\pr_{\alpha\beta} u
				\right)	,\\
				\tilde{A}_{i\alpha}
				&=
				S^{-1}\left(
				\pr_{i\alpha} u 
				-
				(1+\frac{u}{\varrho})^{-1}\frac{\pr_{\alpha}u \pr_i u}{\varrho}
				\right)	,
				\quad 		\widetilde{A}_{ij}
				=
				S^{-1}\left(\pr_{ij}u
				\right)	.
			\end{aligned}	
		\end{equation*}
		
		In conclusion, we have (note that the convention in mean curvature flow is $H=-\tr A$)
		\begin{equation}\label{eq:H}
			\begin{split}
				-\widetilde{H}
				=&
				S^{-1}
				\left(
				(1+\frac{u}{\varrho})^{-2}\Delta_\theta u
				+\Delta_{\R^k} u
				-\pr_{ij}u\pr_i u\pr_j u
				-\pr_{i\alpha} u\pr_i u\pr_\alpha u
				\right.
				\\
				&
				\left.-(1+\frac{u}{\varrho})^{-4}\pr_\alpha u \pr_\beta u\pr_{\alpha\beta}u-(1+\frac{u}{\varrho})^{-1}\frac{n-k}{\varrho}
				-(1+\frac{u}{\varrho})^{-3}\frac{|\nabla_\theta u|^2}{\varrho}
				+m
				\right)
			\end{split}
		\end{equation}
		where $m$ is a remainder function of the linear combinations of products of $\nabla u$, with order at least $4$. In fact, we have
		\begin{equation}\label{Eqm1}
			|m|\leq C|\nabla u|^4,\quad  |\nabla m|\leq C|\nabla u|^3.
		\end{equation}
		
		Next we consider the term $\frac{1}{2}\langle \widetilde z,\widetilde \bn\rangle$. On $\cC_{n,k}$ near $z_0=(\theta_0,y_0)$, we have
		\begin{equation*}
			\begin{split}
				\langle \widetilde z,\widetilde\bn\rangle
				=&
				S^{-1}\left\langle (z)+\frac{u}{\varrho}\theta
				, \frac{\theta}{\varrho}-(1+\frac{u}{\varrho})^{-1}\pr_\alpha u\theta_{\alpha} -\pr_i u \pr_i
				\right\rangle
				=
				S^{-1}
				\left(
				\varrho+u
				-y_i\pr_i u
				\right)
			\end{split}
		\end{equation*}
		Therefore, we can obtain the equation of $u$ from the rescaled mean curvature flow equation (modulo diffeomorphism) $(\pr_\tau  \widetilde z)^{\bot}=-(\widetilde H-\frac{\langle \widetilde z,\widetilde\bn\rangle}{2})\widetilde \bn$. Take an inner product of the equation with $\widetilde \bn$, we get 
		\begin{equation*}
			\begin{split}
				S^{-1} \pr_\tau  u
				=&
				S^{-1}
				\left(
				(1+\frac{u}{\varrho})^{-2}\Delta_\theta u
				+\Delta_{\R^k} u
				-\pr_{ij}u\pr_i u\pr_j u
				-\pr_{i\alpha} u\pr_i u\pr_\alpha u
				-(1+\frac{u}{\varrho})^{-4}\pr_\alpha u \pr_\beta u\pr_{\alpha\beta}u
				\right.
				\\
				&
				\left.-(1+\frac{u}{\varrho})^{-1}\frac{n-k}{\varrho}
				-(1+\frac{u}{\varrho})^{-3}\frac{|\nabla_\theta u|^2}{\varrho}
				+m
				\right)+S^{-1}\frac{1}{2}(\varrho+u-y_i\pr_i u).
			\end{split}
		\end{equation*}
		Thus we obtain the equation of $u$ as follows:
		\begin{equation}\label{EqRMCFFull}
			\begin{split}
				&\pr_\tau  u
				=
				\left(
				(1+\frac{u}{\varrho})^{-2}\Delta_\theta u
				+\Delta_{\R^k} u
				-\pr_{ij}u\pr_i u\pr_j u
				-\pr_{i\alpha} u\pr_i u\pr_\alpha u	-(1+\frac{u}{\varrho})^{-1}\frac{n-k}{\varrho}\right.\\
				&\left.-(1+\frac{u}{\varrho})^{-4}\pr_\alpha u \pr_\beta u\pr_{\alpha\beta}u
				-(1+\frac{u}{\varrho})^{-3}\frac{|\nabla_\theta u|^2}{\varrho}
				+m
				\right)
				+\frac{1}{2}(\varrho+u-y_i\pr_i u).
			\end{split}
		\end{equation}
		When $|u|_{C^0}\leq \eps_0$, we can further rewrite the equation as 
		$\pr_\tau  u
		=
		Lu+\cQ(J^2u),$
		where $\cQ(J^2u)=m_1+m_2+m_3+m_4$, with $m_4$ consisting of all terms like $\pr u\pr u\pr^2 u$ and 
		\begin{equation*}
			m_1=m,
			\ \ \ m_2=-\frac{1}{2\varrho}\frac{u^2}{1+\frac{u}{\varrho}}-(1+\frac{u}{\varrho})^{-3}\frac{|\nabla_\theta u|^2}{\varrho},
			\ \ \ 
			m_3=-\frac{\frac{2u}{\varrho}+\frac{u^2}{\varrho^2}}{(1+\frac{u}{\varrho})^2}\Delta_\theta u.
		\end{equation*}
		From \eqref{EqRMCFFull}, we get \eqref{EqQ} immediately. Indeed, the leading terms come from $m_2$ and $m_3$ respectively. The terms $m_1$ (see \eqref{Eqm1}) and $m_4$ contribute only to $\cC(J^2u)$. Estimate \eqref{EqQ<} also follows immediately from \eqref{EqRMCFFull}. In particular, when $\|u\|_{C^{2}}\leq \eps_0$, we have the estimate
		$		|\cQ|\leq C(|\nabla u|^2+\eps_0|u|).
		$	
		For two different functions $u_1$ and $u_2$, the fundamental theorem of calculus shows that
		\begin{equation*}
			|\cQ(J^2u_1)-\cQ(J^2u_2)|\leq C\eps_0(|u_1-u_2|+|\nabla u_1-\nabla u_2|+|\Hess_{u_1-u_2}|).
		\end{equation*}
		We next consider item (3). To get the terms $\|u\|\|v\|$ on the RHS, we need to perform a step of integration by parts for terms of the form $v\mathrm{Hess}_u$, which also gives us the boundary term $e^{-R^2/4}$. Taking terms $\partial_iu\partial_j u\partial_{ij}u$ in $m_4$ as an example, we have $$\partial_iu\partial_j u\partial_{ij}u=(\nabla u)^T\nabla^2u\nabla u=\frac{1}{2}\nabla\cdot (|\nabla u|^2\nabla u)-\frac{1}{2}\Delta u \cdot |\nabla u|^2. $$ 
		Consider the first term on the RHS
		\begin{equation*}
			\begin{aligned}
				&\nabla\cdot (|\nabla (\chi u_1)|^2\nabla (\chi u_1))-\nabla\cdot (|\nabla (\chi u_2)|^2\nabla (\chi u_2))\\
				&=\nabla\cdot ((\nabla (\chi w)\cdot\nabla(\chi u_1+\chi u_2)) \nabla (\chi u_1))+\nabla\cdot (|\nabla (\chi u_2)|^2\nabla (\chi w))). 		
			\end{aligned}
		\end{equation*}
		When multiplied by $v$ and integrating by parts, we see that it can be bounded by the RHS of \eqref{EqdC}. 
		When the derivative during the integration by parts hits the Gaussian weight $e^{-\frac{|x|^2}{4}}$, we shall get $\int v|x|d\mu$ which is bounded by $\|v\|$ using the Ecker's inequality $\int v^2|x|^2d\mu\leq C\int (v^2+|\nabla v|^2) d\mu. $ 
		All other terms can be treated similarly and are easier, so we get \eqref{EqdC}. 
	\end{proof}

	\section{Perturbed Ornstein--Uhlenbeck regularization}\label{App:Perturbed Ornstein--Uhlenbeck regularization}
	In this section, we prove the perturbed OU regularization Proposition \ref{prop:perturbed-OU}.  The proof is parallel to what we have done for the unperturbed operator. The first step is to give an estimate of the propagator. Let $U(t,s)$ be the positive evolution family of operators generated by $L_r$, and write $U(t,s)=e^{t-s}P(t,s)$.
	
	\begin{lemma}
		\label{lem:perturbed-kernel}
		There are $\bar\delta,h_*,R_0,C,c>0$ such that the following holds. If \eqref{eq:coefficient-smallness-OU} holds, $h=t-s\geq h_*$, and $e^{-h/2}|y|\leq1$, then
		\begin{equation}\label{eq:perturbed-L2-pointwise}
			|U(t,s)\psi(\theta,y)|
			\leq Ce^h\|\psi\|_{L^2}.
		\end{equation}
		For all $t\geq s$,
		\begin{equation}\label{eq:perturbed-moment}
			U(t,s)(1+|z|^2)(\theta,y)
			\leq Ce^h\left(1+e^{-(1-C\delta)h}|y|^2\right).
		\end{equation}
		If $|y|\leq R$, $D\geq3R/2$, and $R\geq R_0$, then
		\begin{equation}\label{eq:perturbed-tail}
			U(t,s)(|z|\mathbf1_{\{|z|\geq D\}})(\theta,y)
			\leq Ce^he^{-cR^2}.
		\end{equation}
	\end{lemma}
	
	\begin{proof}
		Write
		\[
		U(t,s)=e^{h}P(t,s),\qquad h=t-s,
		\]
		where $P(t,s)$ is generated by $\mathscr L_r:=L_r-1$. Let $p(t,s;x,\zeta)$, $x=(\theta,y)$, $\zeta=(\eta,z)$, be its kernel with respect to the unweighted volume.  Since $\mathscr L_r1=0$, we have 
		\begin{equation}\label{eq:perturbed-conservative}
			P(t,s)1=1,\qquad
			\int p(t,s;x,\zeta)\,dV(\zeta)=1.
		\end{equation}
		
		The explicit Mehler formula used in Lemmas \ref{lem:estNr} and \ref{lem:OU-source-estimates} is replaced here by the following two kernel bounds. Let
		\[
		m_h=e^{-h/2}y,\qquad \rho_h=1-e^{-h}.
		\]
		{\bf Claim:} {\it there are $\mu>1/8$, $\beta_0>1/8$, and $h_*>0$ such that
			\begin{equation}\label{eq:perturbed-centered-kernel-exp}
				\int p(t,s;(\theta,y),(\eta,z))
				\exp\!\left(\frac{\mu|z-m_h|^2}{\rho_h}\right)
				\,d\eta\,dz\leq C,
				\qquad h>0,
			\end{equation}
			and
			\begin{equation}\label{eq:perturbed-long-kernel}
				p(t,s;(\theta,y),(\eta,z))
				\leq Ce^{-\beta_0|z-m_h|^2},
				\qquad h\geq h_*.
		\end{equation}}
		
		We can now argue exactly as in Lemma \ref{lem:estNr}.  If $|m_h|=e^{-h/2}|y|\leq1$, choose $\beta\in(1/8,\beta_0)$.  Since
		\[
		\beta_0|z-m_h|^2\geq\beta|z|^2-C_\beta,
		\]
		\eqref{eq:perturbed-long-kernel} implies $p\leq Ce^{-\beta|z|^2}$.  Therefore
		\[
		\begin{aligned}
			|P(t,s)\psi(\theta,y)|
			&\leq \|\psi\|_{L^2}
			\left(\int p(t,s;(\theta,y),\zeta)^2
			e^{|\zeta|^2/4}\,dV(\zeta)\right)^{1/2}  \leq C\|\psi\|_{L^2},
		\end{aligned}
		\]
		because $2\beta-1/4>0$.  Multiplication by $e^h$ proves \eqref{eq:perturbed-L2-pointwise}.
		
		Finally, we argue as in Lemma \ref{lem:OU-source-estimates}.  Directly from \eqref{eq:perturbed-centered-kernel-exp},
		\[
		\int p(t,s;(\theta,y),(\eta,z))
		|z-m_h|^2\,d\eta\,dz\leq C\rho_h.
		\]
		Hence
		\[
		P(t,s)(1+|z|^2)(\theta,y)
		\leq C(1+|m_h|^2+\rho_h)
		\leq C(1+e^{-h}|y|^2).
		\]
		After multiplication by $e^h$, this is stronger than
		\eqref{eq:perturbed-moment}.
		
		If $|y|\leq R$ and $D\geq3R/2$, then on $\{|z|\geq D\}$, we have $|z-m_h|\geq R/2$. Writing $w=z-m_h$, we have
		\[
		\mathbf1_{\{|w|\geq R/2\}}
		\leq e^{-\mu R^2/(4\rho_h)}e^{\mu|w|^2/\rho_h},
		\]
		and
		\[
		|w|\mathbf1_{\{|w|\geq R/2\}}
		\leq C\rho_h^{1/2}e^{-\mu R^2/(8\rho_h)}
		e^{\mu|w|^2/\rho_h}.
		\]
		Since $|z|\leq|w|+R$,
		\eqref{eq:perturbed-centered-kernel-exp} gives
		\[
		P(t,s)(|z|\mathbf1_{\{|z|\geq D\}})
		\leq C(1+R)e^{-cR^2/\rho_h}
		\leq Ce^{-c'R^2}
		\]
		for $R\geq R_0$.  The case $h=0$ is immediate.  Multiplication by $e^h$ proves \eqref{eq:perturbed-tail}.
	\end{proof}
	
	\begin{proof}[Proof of the claim]
		
		To prove the first bound, decrease $\bar\delta$ so that
		\[
		a^{ij}\nu_i\nu_j\leq\Lambda_*|\nu|^2,\qquad
		\Lambda_*<2,\qquad
		B_*:=\sup|b^y|\leq C\bar\delta .
		\]
		For $\lambda\in\mathbb R^k$, direct verification gives that the function
		\[
		\Phi_\lambda(q,y)=
		\exp\!\left(e^{-q/2}\lambda\cdot y
		+\Lambda_*(1-e^{-q})|\lambda|^2
		+2B_*(1-e^{-q/2})|\lambda|\right)
		\]
		satisfies (write $\ell_q=e^{-q/2}\lambda$)
		\[
		(\partial_q-\mathscr L_{s+q})\Phi_\lambda\geq(\Lambda_*e^{-q}|\lambda|^2-a^{ij}\ell_{q,i}\ell_{q,j}
		+B_*e^{-q/2}|\lambda|-b^i\ell_{q,i})\Phi_\lambda\geq 0,\qquad
		\Phi_\lambda(0,y)=e^{\lambda\cdot y}.
		\]

		This gives $ P(t,s)e^{\lambda\cdot z}(\theta,y)\leq\Phi_\lambda(h,y),$ rewritten as 
		\begin{equation}\label{EqB6}
			\int p(t,s;(\theta,y),(\eta,z))
			e^{\lambda\cdot(z-m_h)}\,d\eta\,dz
			\leq
			\exp\!\left(\Lambda_*\rho_h|\lambda|^2
			+2B_*(1-e^{-h/2})|\lambda|\right).
		\end{equation}
		Choose $\varepsilon>0$ so that $A:=\Lambda_*+\varepsilon<2$. Since $1-e^{-h/2}\leq\rho_h^{1/2}$, Young's inequality gives
		\begin{equation}\label{eq:perturbed-centered-linear-exp}
			\int p(t,s;(\theta,y),(\eta,z))
			e^{\lambda\cdot(z-m_h)}\,d\eta\,dz
			\leq C e^{A\rho_h|\lambda|^2}.
		\end{equation}
		Choosing $\frac18<\mu<\frac1{4A}$, multiplying \eqref{eq:perturbed-centered-linear-exp} by $\left(\frac{\rho_h}{4\pi\mu}\right)^{k/2}e^{-\rho_h|\lambda|^2/(4\mu)}$and integrating in $\lambda$ proves \eqref{eq:perturbed-centered-kernel-exp}, where we have used $e^{\mu|z|^2}=(4\pi\mu)^{-k/2}\int_{\R^k}e^{-|\lambda|^2/(4\mu)+\lambda\cdot z}\,d\lambda$.

		We next prove \eqref{eq:perturbed-long-kernel}. For $0<\ell\leq1$, the change $w(q,\theta,x)=v(r_0+q,\theta,e^{q/2}x)$ removes the drift $-\frac12y\cdot\nabla_y$. The transformed operator is uniformly parabolic with uniformly bounded coefficients. Its principal coefficients are spatially $C^1$, so it may be written in divergence form with bounded lower-order coefficients. The standard short-time Aronson estimate \cite[Chapter 1]{Friedman64_PDEparabolic} gives
		\begin{equation}\label{eq:perturbed-short-kernel}
			p(r_0+\ell,r_0;(\theta,y),(\eta,z))
			\leq C\ell^{-n/2}\exp\!\left[-\frac{c}{\ell}
			\bigl(d_\Sigma(\theta,\eta)^2+|z-e^{-\ell/2}y|^2\bigr)\right],
		\end{equation}
		uniformly in $r_0$ and $0<\ell\leq1$.
		
		Choose $\alpha\in(0,1)$ such that $\mu\alpha^2>1/8$, and then choose $\ell>0$ so small that $\frac{c(1-\alpha)^2}{\ell}>\frac18$. Take $h_*>2\ell$. For $h\geq h_*$, the semigroup property shows that
		\[
		p(t,s;x,\zeta)
		=\int p(t,s+\ell;x,\xi)p(s+\ell,s;\xi,\zeta)\,dV(\xi).
		\]
		Let $ m'=e^{-(h-\ell)/2}y=e^{\ell/2}m_h$, and split the integral into
		\[
		E_z=\{|y'-m'|\geq\alpha e^{\ell/2}|z-m_h|\}
		\]
		and $E_z^c$.  On $E_z$, the short-time kernel is bounded and \eqref{eq:perturbed-centered-kernel-exp}, applied to the long factor, gives
		\[
		\int_{E_z}\cdots
		\leq C\exp\!\left(
		-\frac{\mu\alpha^2e^\ell}{\rho_{h-\ell}}|z-m_h|^2\right)
		\leq Ce^{-\mu\alpha^2|z-m_h|^2}.
		\]
		On $E_z^c$,
		\[
		|z-e^{-\ell/2}y'|\geq(1-\alpha)|z-m_h|,
		\]
		so \eqref{eq:perturbed-short-kernel} and \eqref{eq:perturbed-conservative} give
		\[
		\int_{E_z^c}\cdots
		\leq C\exp\!\left(
		-\frac{c(1-\alpha)^2}{\ell}|z-m_h|^2\right).
		\]
		This proves \eqref{eq:perturbed-long-kernel} for some $\beta_0>1/8$.
		
	\end{proof}
	\begin{proof}[Proof of Proposition \ref{prop:perturbed-OU}]
		Variation of constants gives \eqref{eq:Mehler-L2-pointwise} with $S(\tau-r)$ replaced by $U(\tau,r)$. For $|y|\leq R(\tau)$, the choice of $H(\tau)$ and \eqref{eq:perturbed-L2-pointwise} bound the initial term by $CC_0e^{\eps_0H}e^H\tau^{-\xi}$. Equation \eqref{eq:perturbed-moment} bounds the term involving $\frac{1}{\tau^2}(1+|y|^2)$ by $C C_F\tau^{-2}e^{(\eps_0+C\delta)H}
		\bigl(e^H+H R(\tau)^2\bigr)$. Using \eqref{eq:perturbed-tail} we bound the term involving $|y|\check\chi$ by $C C_Fe^{\eps_0H}e^He^{-cR(\tau)^2}$. Since $\eps_0+\delta\leq Cc_*\eta$, take $c_*$ small enough that $e^{(\eps_0+C\delta)H}\leq C_\eta\tau^\eta$. Using $e^H\approx (1+R(\tau))^2$ proves \eqref{eq:perturbed-OU-conclusion}. 
	\end{proof}

	\bibliographystyle{alpha}
	\bibliography{GMT}

\begin{thebibliography}{CHHW22}

\bibitem[AAG95]{AAG}
Steven Altschuler, Sigurd~B. Angenent, and Yoshikazu Giga.
\newblock Mean curvature flow through singularities for surfaces of rotation.
\newblock {\em J. Geom. Anal.}, 5(3):293--358, 1995.

\bibitem[ADS19]{ADS1}
Sigurd Angenent, Panagiota Daskalopoulos, and Natasa Sesum.
\newblock Unique asymptotics of ancient convex mean curvature flow solutions.
\newblock {\em J. Differential Geom.}, 111(3):381--455, 2019.

\bibitem[ADS20]{ADS2}
Sigurd Angenent, Panagiota Daskalopoulos, and Natasa Sesum.
\newblock Uniqueness of two-convex closed ancient solutions to the mean
  curvature flow.
\newblock {\em Ann. of Math. (2)}, 192(2):353--436, 2020.

\bibitem[AV97]{AV}
S.~B. Angenent and J.~J.~L. Vel\'{a}zquez.
\newblock Degenerate neckpinches in mean curvature flow.
\newblock {\em J. Reine Angew. Math.}, 482:15--66, 1997.

\bibitem[BC19]{BC1}
Simon Brendle and Kyeongsu Choi.
\newblock Uniqueness of convex ancient solutions to mean curvature flow in
  {$\Bbb R^3$}.
\newblock {\em Invent. Math.}, 217(1):35--76, 2019.

\bibitem[BC21]{BC2}
Simon Brendle and Kyeongsu Choi.
\newblock Uniqueness of convex ancient solutions to mean curvature flow in
  higher dimensions.
\newblock {\em Geom. Topol.}, 25(5):2195--2234, 2021.

\bibitem[BK23]{BamlerKleiner23_multiplicity}
Richard~H Bamler and Bruce Kleiner.
\newblock On the multiplicity one conjecture for mean curvature flows of
  surfaces.
\newblock {\em arXiv preprint arXiv:2312.02106}, 2023.

\bibitem[CCMS24]{CCMS}
Otis Chodosh, Kyeongsu Choi, Christos Mantoulidis, and Felix Schulze.
\newblock Mean curvature flow with generic initial data.
\newblock {\em Invent. Math.}, 237(1):121--220, 2024.

\bibitem[CCS23]{CCS}
Otis Chodosh, Kyeongsu Choi, and Felix Schulze.
\newblock Mean curvature flow with generic initial data ii.
\newblock {\em arXiv preprint arXiv:2302.08409}, 2023.

\bibitem[CHH22]{CHH}
Kyeongsu Choi, Robert Haslhofer, and Or~Hershkovits.
\newblock Ancient low-entropy flows, mean-convex neighborhoods, and uniqueness.
\newblock {\em Acta Math.}, 228(2):217--301, 2022.

\bibitem[CHHW22]{CHHW}
Kyeongsu Choi, Robert Haslhofer, Or~Hershkovits, and Brian White.
\newblock Ancient asymptotically cylindrical flows and applications.
\newblock {\em Invent. Math.}, 229(1):139--241, 2022.

\bibitem[CM12]{CM1}
Tobias~H. Colding and William~P. Minicozzi, II.
\newblock Generic mean curvature flow {I}: generic singularities.
\newblock {\em Ann. of Math. (2)}, 175(2):755--833, 2012.

\bibitem[CM15]{CM2}
Tobias~Holck Colding and William~P. Minicozzi, II.
\newblock Uniqueness of blowups and \l ojasiewicz inequalities.
\newblock {\em Ann. of Math. (2)}, 182(1):221--285, 2015.

\bibitem[CM16]{CM3}
Tobias~Holck Colding and William~P Minicozzi.
\newblock The singular set of mean curvature flow with generic singularities.
\newblock {\em Inventiones mathematicae}, 204(2):443--471, 2016.

\bibitem[CM19]{ColdingMinicozzi19_Dynamics}
Tobias~Holck Colding and William~P. Minicozzi, II.
\newblock Dynamics of closed singularities.
\newblock {\em Ann. Inst. Fourier (Grenoble)}, 69(7):2973--3016, 2019.

\bibitem[CM21]{CM7}
Tobias~Holck Colding and William~P. Minicozzi, II.
\newblock Wandering singularities.
\newblock {\em J. Differential Geom.}, 119(3):403--420, 2021.

\bibitem[CMP15]{CMP}
Tobias~Holck Colding, William~P. Minicozzi, II, and Erik~Kj{\ae}r Pedersen.
\newblock Mean curvature flow.
\newblock {\em Bull. Amer. Math. Soc. (N.S.)}, 52(2):297--333, 2015.

\bibitem[DHHS26]{daniels2026generic}
JM~Daniels-Holgate, Or~Hershkovits, and G{\'a}bor Sz{\'e}kelyhidi.
\newblock Generic mean curvature flow has finitely many singularities.
\newblock {\em arXiv preprint arXiv:2610.12763}, 2026.

\bibitem[DZ22]{DuZhu22_spectral}
Wenkui Du and Jingze Zhu.
\newblock Spectral quantization for ancient asymptotically cylindrical flows.
\newblock {\em arXiv preprint arXiv:2211.02595}, 2022.

\bibitem[EH91]{EH}
Klaus Ecker and Gerhard Huisken.
\newblock Interior estimates for hypersurfaces moving by mean curvature.
\newblock {\em Inventiones mathematicae}, 105(1):547--569, 1991.

\bibitem[FKZ00]{KZ}
Clotilde Fermanian~Kammerer and Hatem Zaag.
\newblock Boundedness up to blow-up of the difference between two solutions to
  a semilinear heat equation.
\newblock {\em Nonlinearity}, 13(4):1189--1216, 2000.

\bibitem[Fri64]{Friedman64_PDEparabolic}
Avner Friedman.
\newblock {\em Partial differential equations of parabolic type}.
\newblock Prentice-Hall, Inc., Englewood Cliffs, NJ, 1964.

\bibitem[Gan21]{G1}
Zhou Gang.
\newblock On the mean convexity of a space-and-time neighborhood of generic
  singularities formed by mean curvature flow.
\newblock {\em J. Geom. Anal.}, 31(10):9819--9890, 2021.

\bibitem[Gan22]{G2}
Zhou Gang.
\newblock On the dynamics of formation of generic singularities of mean
  curvature flow.
\newblock {\em J. Funct. Anal.}, 282(12):Paper No. 109458, 73, 2022.

\bibitem[GW25]{gang2025precise}
Zhou Gang and Shengwen Wang.
\newblock Precise asymptotics near a generic s1$\times$ r3 singularity of mean
  curvature flow.
\newblock {\em Nonlinear Analysis}, 251:113679, 2025.

\bibitem[Hui90]{Hu1}
Gerhard Huisken.
\newblock Asymptotic behavior for singularities of the mean curvature flow.
\newblock {\em Journal of Differential Geometry}, 31(1):285--299, 1990.

\bibitem[HV92a]{HV2}
M.~A. Herrero and J.~J.~L. Vel\'{a}zquez.
\newblock Blow-up profiles in one-dimensional, semilinear parabolic problems.
\newblock {\em Comm. Partial Differential Equations}, 17(1-2):205--219, 1992.

\bibitem[HV92b]{HV1}
M.~A. Herrero and J.~J.~L. Vel\'{a}zquez.
\newblock Generic behaviour of one-dimensional blow up patterns.
\newblock {\em Ann. Scuola Norm. Sup. Pisa Cl. Sci. (4)}, 19(3):381--450, 1992.

\bibitem[HW20]{HershkovitsWhite20_Nonfattening}
Or~Hershkovits and Brian White.
\newblock Nonfattening of mean curvature flow at singularities of mean convex
  type.
\newblock {\em Communications on Pure and Applied Mathematics}, 73(3):558--580,
  2020.

\bibitem[Ilm94]{I2}
Tom Ilmanen.
\newblock {\em Elliptic regularization and partial regularity for motion by
  mean curvature}, volume 520.
\newblock American Mathematical Soc., 1994.

\bibitem[Ilm03]{I1}
Tom Ilmanen.
\newblock Problems in mean curvature flow.
\newblock {\em preprint}, 2003.

\bibitem[INS19]{INS}
Tom Ilmanen, Andr\'{e} Neves, and Felix Schulze.
\newblock On short time existence for the planar network flow.
\newblock {\em J. Differential Geom.}, 111(1):39--89, 2019.

\bibitem[Kat95]{Kato95_PertubLinearOp}
Tosio Kato.
\newblock {\em Perturbation theory for linear operators}.
\newblock Classics in Mathematics. Springer-Verlag, Berlin, 1995.
\newblock Reprint of the 1980 edition.

\bibitem[SS20]{SS}
Felix Schulze and Natasa Sesum.
\newblock Stability of neckpinch singularities.
\newblock {\em arXiv preprint arXiv:2006.06118}, 2020.

\bibitem[SWX25a]{SunWangXue1_Passing}
Ao~Sun, Zhihan Wang, and Jinxin Xue.
\newblock Passing through nondegenerate singularities in mean curvature flows.
\newblock {\em arXiv preprint arXiv:2501.16678}, 2025.

\bibitem[SWX25b]{SunWangXue2_RegSing}
Ao~Sun, Zhihan Wang, and Jinxin Xue.
\newblock Regularity of cylindrical singular sets of mean curvature flow.
\newblock {\em Preprint}, 2025.

\bibitem[SWZ24]{SunWangZhou20_MinmaxShrinker}
Ao~Sun, Zhichao Wang, and Xin Zhou.
\newblock Multiplicity one for min-max theory in compact manifolds with
  boundary and its applications.
\newblock {\em Calc. Var. Partial Differential Equations}, 63(3):Paper No. 70,
  52, 2024.

\bibitem[SX21a]{SX2}
Ao~Sun and Jinxin Xue.
\newblock Generic dynamics of mean curvature flows with asymptotically conical
  singularities.
\newblock {\em arXiv preprint arXiv:2107.05066, accepted by Science China
  Mathematics}, 2021.

\bibitem[SX21b]{SX1}
Ao~Sun and Jinxin Xue.
\newblock Generic dynamics of mean curvature flows with closed singularities.
\newblock {\em arXiv preprint arXiv:2104.03101}, 2021.

\bibitem[SX24]{SX3}
Ao~Sun and Jinxin Xue.
\newblock Generic regularity of level set flows with spherical singularity.
\newblock {\em Ann. PDE}, 10(1):Paper No. 7, 19, 2024.

\bibitem[Sz{\'e}26a]{szekelyhidi2026mean}
G{\'a}bor Sz{\'e}kelyhidi.
\newblock Mean convex flows with surgery.
\newblock {\em arXiv preprint arXiv:2606.00964}, 2026.

\bibitem[Sz{\'e}26b]{szekelyhidi2026nondegenerate}
G{\'a}bor Sz{\'e}kelyhidi.
\newblock Nondegenerate neck pinches along the mean curvature flow.
\newblock {\em arXiv preprint arXiv:2603.11022}, 2026.

\bibitem[Vel92]{V1}
J.~J.~L. Vel\'{a}zquez.
\newblock Higher-dimensional blow up for semilinear parabolic equations.
\newblock {\em Comm. Partial Differential Equations}, 17(9-10):1567--1596,
  1992.

\bibitem[Vel93]{V2}
J.~J.~L. Vel\'{a}zquez.
\newblock Classification of singularities for blowing up solutions in higher
  dimensions.
\newblock {\em Trans. Amer. Math. Soc.}, 338(1):441--464, 1993.

\bibitem[Whi05]{Wh4}
Brian White.
\newblock A local regularity theorem for mean curvature flow.
\newblock {\em Annals of mathematics}, pages 1487--1519, 2005.

\bibitem[Whi09]{White09_CurrentsVarifolds}
Brian White.
\newblock Currents and flat chains associated to varifolds, with an application
  to mean curvature flow.
\newblock {\em Duke Mathematical Journal}, 148(1):41--62, 2009.

\end{thebibliography}
\end{document}